\documentclass[twoside,11pt]{amsart}
\usepackage{mabliautoref}
\usepackage{amsmath,amsfonts,amssymb,amsthm,mathtools,multirow,tabls,mathrsfs}
 \usepackage[a4paper,margin=1.1in]{geometry}  
\usepackage{lmodern}
\usepackage[all]{xy}
\usepackage[utf8]{inputenc}
\usepackage[T1]{fontenc}
\usepackage[shortlabels]{enumitem}
  \usepackage{hyperref}

\let\counterwithin\relax
\usepackage{chngcntr}

\numberwithin{equation}{section}

\newcommand{\Sch}{\mathfrak{Sch}}
\newcommand{\Sets}{\mathfrak{Sets}}

\newcommand{\scr}{\mathcal}

\newcommand{\sE}{\scr{E}}
\newcommand{\sF}{\scr{F}}
\newcommand{\sG}{\scr{G}}

\newcommand{\sI}{\scr{I}}
\newcommand{\sJ}{\scr{J}}

\newcommand{\sL}{\scr{L}}

\newcommand{\sM}{\scr{M}}

\newcommand{\sO}{\scr{O}}

\newcommand{\sR}{\scr{R}}

\newcommand{\sT}{\scr{T}}

\newcommand{\sV}{\scr{V}}

\newcommand{\bb}[1]{\mathbb{#1}} 
\newcommand{\FF}{\bb{F}}
\newcommand{\CC}{\bb{C}}
\newcommand{\GG}{\bb{G}}
\newcommand{\NN}{\bb{N}}
\newcommand{\PP}{\bb{P}}
\newcommand{\QQ}{\bb{Q}}

\newcommand{\ZZ}{\bb{Z}}

\renewcommand{\phi}{\varphi}
\newcommand{\isom}{\simeq} 
\newcommand{\piet}{\pi^{\acute{e}t}_{1}}

\newcommand{\et}{{\rm\acute{e}t}}
\newcommand{\ra}{\longrightarrow}    

\newcommand{\cExt}{{\mathscr E}\kern -.5pt xt}
\newcommand{\cHom}{\mathscr{H}\kern -.5pt om}

\newcommand{\stacksproj}[1]{{\cite[Tag~\href{http://stacks.math.columbia.edu/tag/#1}{#1}]{stacks-project}}}

\newcommand{\wt}[1]{{\mathchoice%
  {\raisebox{1.5ex}{\resizebox{1.7ex}{!}{{}\hphantom{i}\ensuremath{{\sim}}}} \hspace{-1.7ex}{#1}}%
  {\smash{\raisebox{1.5ex}{\resizebox{1.7ex}{!}{{}\hphantom{i}\ensuremath{{\sim}}}}\hspace{-1.7ex}{#1 }}\vphantom{\tilde I}}%
  {\raisebox{1.1ex}{\resizebox{1.3ex}{!}{{}\hphantom{i}\ensuremath{{\sim}}}}\hspace{-1.3ex}{#1}}%
  {\raisebox{0.8ex}{\resizebox{1ex}{!}{{}\hphantom{i}\ensuremath{{\sim}}}}\hspace{-1ex}{#1}}%
}}

\newcommand{\cE}{\sE}

\newcommand{\cI}{\sI}
\newcommand{\cJ}{\sJ}

\newcommand{\cL}{\sL}

\newcommand{\cO}{\sO}

\newcommand{\cR}{\sR}

\newcommand{\cT}{\sT}

\newcommand{\cV}{\sV}

\RequirePackage[dvipsnames,usenames]{xcolor}
\usepackage{tikz}
\usepackage{tikz-cd}
\usetikzlibrary{decorations.markings}
\tikzset{degil/.style={
  decoration={markings,
  mark= at position 0.5 with {
  \node[transform shape] (tempnode) {$\backslash$};
  \draw[thick] (tempnode.north east) -- (tempnode.south west);
  }}, postaction={decorate}
}}
\usetikzlibrary{decorations.text, calc, matrix, arrows, decorations.pathreplacing}
\usepackage{enumitem}
\usepackage{chngcntr}
\usepackage{minibox}
\usepackage[framemethod=TikZ]{mdframed}
\usepackage{blkarray}

\setlist[itemize,enumerate]{leftmargin=0.9cm}

\mdfsetup{
 skipabove=6pt,
 skipbelow=3pt,
 roundcorner=10pt
}

\usepackage{xparse}

\NewDocumentCommand\DownArrow{O{2.0ex} O{black}}{%
   \mathrel{\tikz[baseline] \draw [<-, line width=0.5pt, #2] (0,0) -- ++(0,#1);}
}

\NewDocumentCommand\UpArrow{O{2.0ex} O{black}}{%
   \mathrel{\tikz[baseline] \draw [line width=0.5pt, decoration={markings,mark=at position 1 with {\arrow[scale=2, line width=0.25pt]{to}}},  postaction={decorate}, #2] (0,0) -- ++(0,#1);}
}

\newcommand{\expl}[2]{\underset{\mathclap{\minibox[c]{$\UpArrow[10pt]$\\ \fbox{\footnotesize #2}}}}{#1}}

\newcommand{\expllower}[3]{ 
\underset{\mathclap{\minibox[c]{$\UpArrow[#1]$\\ \fbox{\footnotesize #3}}}}{#2}}

\newcommand{\explshift}[3]{\underset{\mathclap{\minibox[c]{$\UpArrow[10pt]$\\ \hspace{#1} \fbox{\footnotesize #3}}}}{#2}}

\newcommand{\explpar}[3]{\underset{\mathclap{\minibox[c]{$\UpArrow[10pt]$\\  \fbox{\parbox{#1}{\footnotesize #3}}}}}{#2}}

\newcommand{\explparshift}[4]{\underset{\mathclap{\minibox[c]{$\UpArrow[10pt]$\\ \hspace{#2} \fbox{\parbox{#1}{\footnotesize #4}}}}}{#3}}

\hypersetup{
bookmarks,
bookmarksdepth=3,
bookmarksopen,
bookmarksnumbered,
pdfstartview=FitH,
colorlinks,backref,hyperindex,
linkcolor=Sepia,
anchorcolor=BurntOrange,
citecolor=MidnightBlue,
citecolor=OliveGreen,
filecolor=BlueViolet,
menucolor=Yellow,
urlcolor=OliveGreen
}

\DeclareMathAlphabet{\mathchanc}{OT1}{pzc}%
                                 {m}{it}

\newcommand{\mcH}{\mathchanc{H}}
\newcommand{\mcm}{\mathchanc{m}}
\newcommand{\mco}{\mathchanc{o}}

\newcommand{\bC}{\mathbb{C}}
\newcommand{\bF}{\mathbb{F}}
\newcommand{\bG}{\mathbb{G}}
\newcommand{\bN}{\mathbb{N}}
\newcommand{\bP}{\mathbb{P}}
\newcommand{\bQ}{\mathbb{Q}}

\newcommand{\bZ}{\mathbb{Z}}

\newcommand{\tH}{\tilde{H}}

\newcommand{\tW}{\tilde{W}}

\newcommand{\ok}{\overline{k}}

\newcommand{\ot}{\overline{t}}

\DeclareMathOperator{\pol}{pol}
\DeclareMathOperator{\Alb}{Alb}
\DeclareMathOperator{\Aut}{Aut}

\DeclareMathOperator{\Hom}{Hom}
\DeclareMathOperator{\Pic}{Pic}
\DeclareMathOperator{\Proj}{Proj}

\DeclareMathOperator{\Spec}{Spec}
\DeclareMathOperator{\Sym}{Sym}
\DeclareMathOperator{\gen}{{gen}}
\DeclareMathOperator{\pr}{{pr}}
\DeclareMathOperator{\charfield}{{char}}

\DeclareMathOperator{\fin}{{fin}}

\DeclareMathOperator{\characteristic}{{char}}

\DeclareMathOperator{\alb}{{alb}}

\DeclareMathOperator{\codim}{codim}

\DeclareMathOperator{\Ker}{{Ker}}
\DeclareMathOperator{\pre}{{pre}}

\DeclareMathOperator{\Exc}{Exc}
\DeclareMathOperator{\Hilb}{Hilb}
\newcommand{\sHom}[0]{{\mcH\mco\mcm}}    
\DeclareMathOperator{\id}{{id}}

\DeclareMathOperator{\Isom}{Isom}

\DeclareMathOperator{\specialpt}{special}

\DeclareMathOperator{\red}{red}

\DeclareMathOperator{\rk}{{rk}}
\DeclareMathOperator{\coeff}{{coeff}}
\DeclareMathOperator{\supp}{{supp}}
\DeclareMathOperator{\Supp}{{Supp}}

\newcommand{\factor}[2]{\left. \raise 2pt\hbox{\ensuremath{#1}} \right/
        \hskip -2pt\raise -2pt\hbox{\ensuremath{#2}}}

\counterwithin*{equation}{section}
\counterwithin*{equation}{subsection}

\makeatletter
\renewcommand\subsection{
  \renewcommand{\sfdefault}{pag}
  \@startsection{subsection}%
  {2}{0pt}{.8\baselineskip}{.4\baselineskip}{\raggedright
    \sffamily\itshape\small\bfseries
  }}
\renewcommand\section{
  \renewcommand{\sfdefault}{phv}
  \@startsection{section} %
  {1}{0pt}{\baselineskip}{.8\baselineskip}{\centering
    \sffamily
    \scshape
    \bfseries
}}
\makeatother

\renewcommand{\theenumi}{\arabic{enumi}}

\newcommand{\wh}{\widehat}

\newcommand{\obF}{\overline{\mathbb{F}}}

\usetikzlibrary{shapes.geometric, arrows, positioning, fit, calc,decorations.pathmorphing,shapes}
\tikzstyle{e} = [ellipse, minimum width=1cm, minimum height=0.5cm,text centered, draw=black]
\tikzstyle{arrow} = [thick,-,>=stealth]
\tikzstyle{connection}=[inner sep=0,outer sep=0]

\author[Zs. Patakfalvi]{Zsolt Patakfalvi}
\address{\'Ecole Polytechnique F\'ed\'erale de Lausanne, Chair of Algebraic Geometry \newline 
\indent MA C3 625 (Bâtiment MA), Station 8, CH-1015 Lausanne}
\email{zsolt.patakfalvi@epfl.ch}

  \author[M. Zdanowicz]{Maciej Zdanowicz}
  \address{\'Ecole Polytechnique F\'ed\'erale de Lausanne, Chair of Algebraic Geometry \newline 
    \indent MA C3 585 (Bâtiment MA), Station 8, CH-1015 Lausanne}
  \email{maciej.zdanowicz@epfl.ch}

  \title[The weak Beauville--Bogomolov decomposition in characteristic $p\geq 0$ ]{The weak Beauville--Bogomolov decomposition in characteristic $p\geq 0$ \\[10pt] {\normalfont \footnotesize (with an appendix joint with Giulio Codogni)}}
  \date{\today}
  \subjclass[2010]{Primary 14G17, Secondary 14M17, 14M25, 14J45} 
  \keywords{
    Beauville--Bogomolov decomposition, trivial canonical class, characteristic $p$, ordinary}
    
\begin{document}

\begin{abstract}

We prove a variant of the Beauville--Bogomolov decomposition for weakly ordinary, or generally globally $F$-split, varieties $X$ with $K_X \sim 0$, in characteristic $p>0$.  We also show that the weakly ordinary assumption in our statement cannot be dropped. Additionally, if the assumption $K_X \sim 0$ is replaced by  $-K_X$ being semi-ample, we show the weaker statement that all closed fibers of the Albanese morphism are isomorphic.  
Finally, we apply our main theorem to draw consequences to the behavior of rational points and fundamental groups of weakly ordinary $K$-trivial varieties in positive characteristic.

{\scshape French version:} On d\'emontre un variant de la d\'ecomposition  Beauville--Bogomolov pour les vari\'et\'es faiblement ordinaire, ou $F$-scind\'e en g\'en\'eral, avec $K_S \sim 0$ en caract\'eristic $p>0$. On d\'emontre aussi que l'hypoth\`ese faiblement ordinaire est n\'ecessaire. De plus, si l'hypoth\`se $K_X \sim 0$ est remplac\'e par $-K_X$ semi-ample, alors on d\'emontre l'affirmation moins forte que toutes fibres du morphisme Albanese sont isomorphes. A la fin, on applique notre th\'or\`eme principal pour arriver aux cons\'equences sur les points rationnels et les groupes foundamentaux des vari\'et\'es  faiblement ordinaires $K$-triviales en caract\'eristique positive. 
\end{abstract}

\maketitle

\tableofcontents


\section{Introduction}

\subsection{Historical motivation}
\label{sec:historical_motivation}
The classical Beauville--Bogomolov decomposition (see \cite{Bogomolov_The_decomposition_of_Kahler_manifolds_with_a_trivial_canonical_class,
Beauville_Varietes_Kahleriennes_dont_la_premiere_classe_de_Chern_est_nulle})  
states that up to an \'etale cover, a compact K\"ahler manifold can be split into a product of an abelian variety, a Calabi--Yau variety $U$ (satisfying $H^i(U,\cO_{U}) = 0$, for $0 < i < \dim U$) and a hyperk\"ahler manifold $W$ (simply connected and admitting a holomorphic symplectic form). 

The main goal of this article is to investigate what one can say about the existence of such decomposition over fields of positive characteristic. The above decomposition is based on the existence of K\"ahler--Einstein metrics, an analytic tool whose algebraic substitute is not known even in charateristic zero. Hence, in this article we analyze the weaker question that is approachable by algebraic tools, and to which we refer to as \emph{weak Beauville--Bogomolov} decomposition. 
Over $\bC$, it states that for every smooth projective variety $X$ with $K_X \sim 0$ there exists a finite \'etale cover 
\begin{equation}
\label{eq:weak_BB_classical}
V \times B \to X,
\end{equation} 
where $B$ is an abelian variety and $V$ is a smooth projective simply connected variety with $K_V \sim 0$.  Already this statement has substantial corollaries concerning geometry of varieties with trivial canonical class.  For example, it directly implies that the fundamental group of such variety is virtually abelian, that is, admits a finite index abelian subgroup.   We note that the weak Beauville--Bogomolov was already observed by Calabi for K\"ahler manifolds in \cite{Calabi_On_Kahler_manifolds_with_vanishing_canonical_class} (see also \cite{Matsushima_Holomorphic_vector_fields_and_the_first_Chern_class}).

Additionally, instead of requiring $V$ to be simply connected, the decomposition \autoref{eq:weak_BB_classical} can be alternatively also characterized up to an \'etale cover of $V$ by using the notion of \emph{augmented irregularity}. By definition, the augmented irregularity of $X$ is
\begin{equation}
\label{eq:augmented_irregularity}
\wh{q}(X) := \max \left\{ \left. \dim \Alb_{X'} \right|  X' \to X \text{ is a finite \'etale morphism}\right\},
\end{equation}
for which one needs to prove that the above maximum exists. This follows from Calabi's original statement, or from \cite[Thm 1]{Kawamata_Characterization_of_abelian_varieties} in a more general singular setting. Then, the above mentioned alternative characterization of \autoref{eq:weak_BB_classical} is by requiring that $\wh{q}(X)=\dim B$ and $\wh{q}(V)=0$. 

\subsection{List of the results}

In this subsection $X$ is a smooth projective variety of dimension $n$ over the stated base-field  $k$. If we work over a perfect field of charateristic $p>0$ and $K_X \sim 0$, then the only cohomology group that is guaranteed to be non-zero is $H^n(X, \sO_X) \cong k$. In particular, the only arithmetic invariant one can  extract is whether the Frobenius morphism acts on this cohomology as a  bijection or not. We say that $X$ is \emph{weakly ordinary} if this action is bijective. Our results are:

\begin{itemize}[leftmargin=18pt,labelwidth=0pt,itemindent=0pt,label={\huge\textbullet}]
\item 
We fully answer the question when a weak Beauville--Bogomolov decomposition exists. That is, 
\begin{enumerate}
\item In the weakly ordinary case we clarify the relaxations of the statement over $\bC$ (see \autoref{sec:historical_motivation}) needed to obtain a weak Beauville--Bogomolov decomposition:
\begin{enumerate}
\item $V$ has to be allowed mildly singular, and 
\item the morphism $V \times B \to X$ has to be allowed to be inseparable, i.e., a torsor under a possibly non-reduced group scheme.
\end{enumerate}
(\autoref{thm:smooth_BB_decomp} is the smooth case, and \autoref{thm:beauville_bogomolov_general_logcy} is the  singular and pair case.)\\[5pt]
Additionally, we note that the above two relaxations are unavoidable (see \cite[Thm 2.11]{Lam_On_a_question_of_Patakfalvi_Zdanowicz} and \autoref{ex:Lam}).

\item If $X$ is not weakly ordinary, then we find an example that cannot admit a a decomposition (\autoref{prop:bb_decomp_counterexample_moret_bailly}).

\end{enumerate}

\item We deduce corollaries to the number of rational points and to the behavior of the \'etale fundamental group (see \autoref{cor:fundamental_group_virtually_abelian} and \autoref{cor:rational_points}).

\item We solve the positive characteristic version of the Demailly--Peternell--Schneider conjecture explained in the following paragraph (\autoref{thm:smooth_isotriviality_Albanese} is the smooth case, and  \autoref{thm:decomposition_theorem} is the singular and pair case). This proof is a byproduct of our methods for solving the above Beauville-Bogomolov question, and it requires practically no additional effort. 

\end{itemize}

\subsection{The smooth cases of the characteristic $p>0$ results}
\label{sec:results_smooth_case}

\emph{Our base field $k$ is perfect and of characteristic $p>0$}, and we also need  the following positive characteristic notion:  a projective variety $X$ over $k$ is \emph{weakly ordinary} if the action of the absolute Frobenius morphism $F_X$ on $H^{\dim X} ( X, \sO_X)$ is bijective. This is a genericity notion. That is, being weakly ordinary is an open condition in positive equicharacteristic and it is typically dense in moduli. Additionally, it is conjectured to be dense over mixed characteristic bases that are finite type over $\bZ$ \cite{Mustata_Srinivas_Ordinary_varieties_and}.
For a smooth projective weakly ordinary variety $X$ with $K_X \sim 0$ 
we define the \emph{augmented irregularity} just as in \autoref{eq:augmented_irregularity}, where the existence of the maximum is guaranteed by \cite[Thm 1.1]{Ejiri_When_is_the_Albanese_morphism_an_algebraic_fiber_space_in_positive___characteristic?}.

Also, already when $X$ is smooth, our decomposition statement uses the notion of strongly $F$-regular singularities, which is a characteristic $p$ class of mild singularities. In particular, it is contained both in the class of klt and rational singularities.  We refer to the many surveys in the topic for a detailed introduction on the notions of $F$-singularities \cite{Schwede_Tucker_A_survey_of_test_ideals,Patakfalvi_Schwede_Tucker_Positive_characteristic_algebraic_geometry,Patakfalvi_Frobenius_techniques_in_birational_geometry}

\begin{theorem}[{\scshape Smooth \& weakly ordinary case of the Beauville--Bogomolov decomposition}, special case of \autoref{thm:beauville_bogomolov_general_logcy}]
\label{thm:smooth_BB_decomp}
\label{cor:beauville_bogomolov}
Let $X$ be smooth  projective variety over $k$, such that $K_X \sim 0$ and $X$ is weakly ordinary. 
Then there is a composition 
\begin{equation*}
B \times V \to Z \to X
\end{equation*}
 of two finite covers,  such that
\begin{enumerate}
\item $Z \to X$ is \'etale, $B \times V \to Z$ is a torsor under $\prod_{i=1}^{\wh{q}(X)} \mu_{p^{j_i}}$ for some integers $j_i \geq 0$, 
\item $B$ is an abelian variety with $\dim B = \wh{q}(X)$, and
\item $V$ is a weakly ordinary  projective Gorenstein variety over $k$ with strongly $F$-regular singularities, $K_V \sim 0$ and  $\wh{q}(V)=0$. 
\end{enumerate}
Additionally, the action of $G:=\prod_{i=1}^n \mu_{p^{j_i}}$ on $B \times V$ is the diagonal action induced by an action on $V$ and an action on $B$, respectively, such that
\begin{itemize}
\item $G$ acts freely and faithfully on $B$, and
\item $G$ acts faithfully on $V$. 
\end{itemize}
\end{theorem}

\begin{remark}
\label{rem:differences_to_char_zero}
There are two major differences between \autoref{thm:smooth_BB_decomp} and the original characteristic zero statement mentioned above:
\begin{enumerate}
\item \label{itm:differences_to_char_zero:inseparable} $B \times V \to X$ in \autoref{thm:smooth_BB_decomp} is not  \'etale as in characteristic zero, but has an infinitesimal part as well, and
\item \label{itm:differences_to_char_zero:singularities} $V$ is not necessarily smooth, but has only strongly $F$-regular singularities.
\end{enumerate}
These two phenomena are in fact interconnected, as both are caused by the possible presence of non-reduced polarized automorphism groups.  In \autoref{sec:examples}, we present examples showing that both of the above extensions of the characteristic zero statements are actually necessary.  We note that our examples are based on constructions by Matsumoto \cite[Example 10.5]{Matsumoto_mu_p_and_alpha_p_actions_on_K3_surfaces} and Lam \cite[Appendix A]{Lam_On_a_question_of_Patakfalvi_Zdanowicz}.
\end{remark}

The weak ordinarity assumption in \autoref{thm:smooth_BB_decomp} cannot be dropped. That is, we have counterexamples to the statement that one would obtain from \autoref{thm:smooth_BB_decomp} by erasing the weak ordinarity assumption. This is phrased in the following proposition, taking into account that in \autoref{thm:smooth_BB_decomp} the projection $B \times V \to B$ can be identified with the Albanese morphism of $B \times V$.

\begin{proposition}
(special case of \autoref{prop:bb_decomp_counterexample_moret_bailly})
\label{prop:bb_cecomp_counterexample_intro}
Set 
\begin{equation*}
d :=
\left\{
\begin{matrix}
3 & \textrm{if }p=2 \\
2p -3 & \textrm{otherwise}
\end{matrix}
\right.
\end{equation*}
Then, there exists a $d$ dimensional smooth projective variety $X$ over $k$ such that $K_X \sim 0$, and the following property is satisfied: whenever  $Y \to Z \to  X$ is a composition of finite morphisms such that
\begin{enumerate}
\item $Y$ is an integral variety, 
\item  $Z \to X$ is \'etale, and 
\item $Y \to Z$ is a  torsor under an infinitesimal   group scheme over $k$,
\end{enumerate}
then, the Albanese morphism $Y \to \Alb Y$ is not split. 
\end{proposition}

\begin{remark}
As explained before its statement, 
\autoref{prop:bb_cecomp_counterexample_intro} shows that weak ordinarity assumption cannot be dropped from \autoref{thm:smooth_BB_decomp}. However, it does not exclude the possibility that weaker versions of \autoref{thm:smooth_BB_decomp} still hold if the assumption is removed. Some options for such weakenings would be:
\begin{enumerate}
\item requiring $p \geq p(d)$, where $p(d)$ depends only on the dimension $d$ of $X$,
\item allowing also other types of morphisms into the cover $Y \to X$ (although it is hard to find other classes of relatively $K$-trivial finite morphisms, so probably in this case one would loose the $K$-triviality of $Y$),
\item relaxing, but not removing the weak ordinarity assumption, when $X$ comes from a special class of varieties with finer Hodge type data.
\end{enumerate}
The questions above are beyond the scope and the goal of the present article, hence we leave them as problems for future investigations.  As a first step, one could verify if $p(3) = 3$ or $5$ using current MMP-based methods.  For example results related to $p(3) \leq 7$ in the case $\dim \Alb X = 1$ can be found in \cite[Section 4.4]{Zhang_Abundance_for_3-folds_with_non-trivial_Albanese_maps_in_positive___characteristic}.
\end{remark}

\begin{remark}
\label{rem:weak_ordinary_global_F_split}
If $X$ is a normal, projective variety over $k$ such that $K_X \sim 0$, then $X$ is weakly ordinary if and only if it is globally $F$-split. That is if and only if the structure morphism $\sO_X \to F_* \sO_X$ splits as a homomorphism of $\sO_X$-modules.
This is the framework in which the statement of \autoref{thm:smooth_BB_decomp} generalizes  to the case when $K_X \equiv 0$. Additionally, in this framework, one can allow also singularities. Hence, in \autoref{sec:results_singular_case}, and in general in most parts of the article, the weakly ordinary condition will be replaced with global $F$-splitting. 

\end{remark}

In \autoref{thm:smooth_isotriviality_Albanese}, we are also able to prove,  the positive characteristic version of the main result of \cite{Cao_Albanese_maps_of_projective_manifolds_with_nef_anticanonical_bundles}, with the caveat that if we are not over a finite field or its algebraic closure, then we have to replace nef by semi-ample. 
This result, i.e, \autoref{thm:smooth_isotriviality_Albanese}, uses  \autoref{thm:smooth_isotriviality}, which might be of interest independently as well. 

These statements are also impossible to state completely without the use of the local and global $F$-singularity notions.  The notion of global $F$-splitting generalizes weak ordinarity, as explained in \autoref{rem:weak_ordinary_global_F_split}, and strong $F$-regularity generalizes smoothness, as mentioned before \autoref{thm:smooth_BB_decomp}. 

\begin{theorem}[{\scshape Isotriviality over curves}, special case of \autoref{cor:finite_field} and \autoref{cor:semi_ample}.]
\label{thm:smooth_isotriviality}
 Let  $f : X \to T$ be a surjective morphism from a  smooth projective variety to a smooth projective curve with strongly $F$-regular general fiber such that   
either
\begin{enumerate}[label=(\alph*)]
\item \label{situation:smooth_isotriviality:semi_ample} $-K_{X/T}$ is semi-ample, or
\item \label{situation:smooth_isotriviality:finite_field} $k \subseteq \obF_p$ and $-K_{X/T} $ is nef.
\end{enumerate}  
Then $f$ has isomorphic closed fibers over $\ok$.
\end{theorem}

Luckily, in the situation of \autoref{thm:smooth_isotriviality_Albanese}, we can prove that  the general fibers of the Albanese morphism are strongly $F$-regular (see \autoref{thm:general_Fiber_SFR2}). This, together with \autoref{thm:smooth_isotriviality}, leads then to the following statement:

\begin{theorem}[{\scshape Isotriviality of the Albanese morphism}, special case of \autoref{thm:decomposition_theorem}]
\label{thm:smooth_isotriviality_Albanese}
Let $X$ be a smooth projective globally $F$-split variety such that
either
\begin{enumerate}[label=(\alph*)]
\item \label{situation:smooth_isotriviality_Albanese:semi_ample} $-K_X$ is semi-ample, or
\item \label{situation:smooth_isotriviality_Albanese:finite_field} $k \subseteq \obF_p$ and $-K_X $ is nef.
\end{enumerate}  
Then the Albanese morphism $f \colon X \to A$ has isomorphic closed fibers over $\ok$.
\end{theorem}

\begin{remark}
\label{rem:trivializes_over_torsor_smooth}
The general version of \autoref{thm:smooth_isotriviality_Albanese} provided in \autoref{thm:decomposition_theorem}, makes a more precise statement: $f$ gets trivialized over a specific torsor under the polarized automorphism group of (any) fiber. 
\end{remark}

We conclude the present section with a few applications of our results. We start with a conjecture concerning existence of rational points on $K$-trivial varieties. 

\begin{conjecture}
\label{conj:rational_points}
For every positive integer $d$ there exists an integer $q_d$ such that if $X$ is a smooth projective  variety of dimension $d = \dim X$ over $\bF_q$ with $K_X \sim 0$ and $q \geq q_d$, then $X( \bF_q) \neq \emptyset$.
\end{conjecture}

Although, in this precise form \autoref{conj:rational_points} might not be written down in the literature, it is most likely well-known for the experts. It  is partially motivated by the following:
\begin{itemize}
\item If there was a  Beauville-Bogomolov decomposition for all $K$-trivial varieties over finite fields together with the boundedness of the irreducible summands, then \autoref{conj:rational_points} would follow via the Weil conjectures, see the proof of \autoref{cor:rational_points_intro}.
\item Calabi-Yau varieties in characteristic zero are conjectured to always have rational curves, as well as rational points over number fields \cite[Problem 10.2]{Miyaoka_Peternell_Geometry_of_higher-dimensional_algebraic_varieties}, c.f. \cite{Wilson_Calabi-Yau_manifolds_with_large_Picard_number}.
\item Symplectic varieties over $\bC$ also tend to have many rational curves, e.g, \cite{Charles_Mongardi_Pacienza_Families_of_rational_curves_on_holomorphic_symplectic_varieties_and___applications_to_zero-cycles,Amerik_Verbitsky_Rational_curves_on_hyperkahler_manifolds}.
\item Over finite fields, geometrically abelian varieties always have rational points \cite[Thm 3]{Lang_Abelian_varieties_over_finite_fields}.
\item On the other hand if $p \geq 5$ is a prime number, then $\left\{\sum_{i=0}^{p-2} x_i^{p-1} =0 \right\} \subseteq \bP^{p-2}_{\bF_p}$ is a smooth Calabi-Yau variety without an $\bF_p$-rational point. Hence, the bound $p_d$ is necessary in \autoref{conj:rational_points}.
\end{itemize}
\autoref{thm:smooth_BB_decomp} implies the following special case of  \autoref{conj:rational_points}.

\begin{corollary} [special case of \autoref{cor:rational_points}]
\label{cor:rational_points_intro}
If $X$ is a smooth projective weakly ordinary and geometrically irreducible $3$-fold over $\bF_q$ with $K_X \sim 0$, $\wh{q}(X)\neq 0$ and $q \geq 83$, then $X(\bF_q) \neq \emptyset$.  
\end{corollary}

Our other application of our main results is towards the conjecture, the characteristic zero counterpart of which is well known, that $K$-trivial smooth projective varieties have virtually abelian \'etale fundamental groups.

\begin{corollary}
[special case of \autoref{cor:fundamental_group_virtually_abelian}]
If $X$ is a smooth projective weakly ordinary variety with $K_X \sim 0$ and $\hat{q}(X)=\dim X -2$, then $\pi_1^{\'et}(X)$ is virtually abelian. 
\end{corollary}

\begin{remark}
We note that in the full statements, that is, in  \autoref{cor:rational_points} and in \autoref{cor:fundamental_group_virtually_abelian} we are able to also include the $K_X$ numerically trivial case by assuming $F$-purity instead of weak ordinarity.  
\end{remark}

\subsection{Statements for singular varieties in characteristic $p>0$}
\label{sec:results_singular_case}

As explained before, the main recipe of turning the statements concerning smooth varieties provided in \autoref{sec:results_smooth_case} into statements about singular varieties is to replace every occurrence of smooth by strongly $F$-regular, and every occurrence of weakly ordinary by globally $F$-split. Hence for \autoref{thm:smooth_BB_decomp} we obtain:

\begin{theorem}[{\scshape Singular Beauville--Bogomolov decomposition},  \autoref{thm:beauville_bogomolov_general_logcy}]
\label{thm:BB_decomp}
Let $(X,\Delta)$ be a globally $F$-split  projective pair over $k$ with strongly $F$-regular singularities, such that $K_X + \Delta \equiv 0$. 
Then there is a  composition $Y \to Z \to X$ of two finite covers  such that $Z \to X$ is quasi-\'etale, $Y \to Z$ is a torsor under $\prod_{i=1}^{\wh{q}(X)} \mu_{p^{j_i}}$  for some integers $j_i \geq 0$, and such that 
\begin{equation*}
(Y, \Delta_Y) \cong (V,\Delta_V) \times B,
\end{equation*}
where 
\begin{enumerate}
\item $B$ is an abelian variety with $\dim B = \wh{q}(X)$. 
\item $(V, \Delta_V)$ is a globally $F$-split  projective pair over $k$ with strongly $F$-regular singularities, $K_V + \Delta_V \equiv 0$ and  $\wh{q}(V)=0$. 
\end{enumerate}\end{theorem}

The proof of the theorem actually works in a more general setting.  We managed to show the following result under the assumption that $-K_X-\Delta$ is nef/semi-ample. 

\begin{theorem}[{\scshape Isotriviality of the Albanese morphism}, \autoref{thm:decomposition_theorem}]
\label{thm:isotriviality_Albanese}
Let $(X,\Delta)$ be a projective globally $F$-split pair with strongly $F$-regular singularities such that  $K_X + \Delta$ is $\bQ$-Cartier with index prime-to-$p$, and 
either
\begin{enumerate}[label=(\alph*)]
\item \label{situation:smooth_isotriviality_Albanese:semi_ample} $-(K_X + \Delta)$ is semi-ample, or
\item \label{situation:smooth_isotriviality_Albanese:finite_field} $k \subseteq \obF_p$ and $-(K_X +\Delta)$ is nef.
\end{enumerate}  
Then, $(X_t,\Delta_t) \cong (X_s, \Delta_s)$ for every  $s,t \in A\left(\ok\right)$.
\end{theorem}

While proving \autoref{thm:isotriviality_Albanese}, we also obtain the following theorem concerning fibrations over curves.

\begin{theorem}[{\scshape Isotriviality over curves},\autoref{cor:finite_field} and \autoref{cor:semi_ample}.]
\label{thm:isotriviality}
 Let  $f : (X,\Delta) \to T$ be a surjective morphism from a   projective normal pair to a smooth projective curve such that $\Delta$ is an effective $\bQ$-divisor, $K_X + \Delta$ is $\bQ-Cartier$, the general fiber $(X_t, \Delta_t)$ is strongly $F$-regular and  
either
\begin{enumerate}[label=(\alph*)]
\item \label{situation:smooth_isotriviality:semi_ample} $-(K_{X/T}+ \Delta)$ is semi-ample, or
\item \label{situation:smooth_isotriviality:finite_field} $k \subseteq \obF_p$ and $-(K_{X/T} + \Delta)$ is nef.
\end{enumerate}  
Then, $(X_t,\Delta_t) \cong (X_s, \Delta_s)$ for every  $s,t \in T\left(\ok\right)$.
\end{theorem}

\begin{remark}
\label{remark:intro_historical}
{\scshape Historical remarks (in characteristic zero):} As mentioned above, the original smooth version of the Beauville--Bogomolov decomposition was shown in \cite{Bogomolov_The_decomposition_of_Kahler_manifolds_with_a_trivial_canonical_class,Beauville_Varietes_Kahleriennes_dont_la_premiere_classe_de_Chern_est_nulle}.
The singular Beauville--Bogomolov decomposition has recently attracted a serious amount of attention which culminated in the series of papers
\cite{Greb_Kebekus_Petternel_Singular_Spaces_with_Trivial_Canonical_Class,Greb_Guenancia_Kebekus_Klt_varieties_with_trivial_canonical_class_-_Holonomy__differential___forms__and_fundamental_groups,Druel_A_decomposition_theorem_for_singular_spaces_with_trivial_canonical_class___of_dimension_at_most_five,
Horing_Peternell_Algebraic_integrability_of_foliations_with_numerically_trivial_canonical___bundle} leading to the full decomposition theorem for klt varieties with numerically trivial canonical class.  The weak Beauville--Bogomolov condition was shown even in the singular logarithmic setting in the papers
\cite{Kawamata_Characterization_of_abelian_varieties,Ambro_The_moduli_b-divisor_of_an_lc-trivial_fibration}.
\end{remark}

\begin{remark}
For the additional statement about trivialization over a flat torsor, mentioned in \autoref{rem:trivializes_over_torsor_smooth}, we refer to \autoref{cor:finite_field}, \autoref{cor:semi_ample} and \autoref{thm:decomposition_theorem}.
\end{remark}


\subsection{Outline of the proof}

In the following section, we give an outline of the proof of \autoref{thm:smooth_BB_decomp}.  This includes all the techniques necessary to get the most general results, i.e.,  the statements of   \autoref{sec:results_singular_case}.  However, for the purpose of clarity we avoid some technical difficulties.   Let $X$ be a smooth weakly ordinary projective variety defined over an algebraically closed field $k$ of characteristic $p>0$.  Suppose that the canonical divisor satisfies the condition $K_X \sim 0$.

\subsubsection{Albanese morphism and augmented irregularity in characteristic $p$}
\label{ss:intro_albanese_morphism_irregularity}

We begin our proof with the application of the results of Ejiri provided in \cite{Ejiri_When_is_the_Albanese_morphism_an_algebraic_fiber_space_in_positive___characteristic?}.  For this purpose, we first apply the standard result of Mehta and Ramanathan \cite[Proposition 9]{Mehta_Ramanathan_Frobenius_splitting_and_cohomology_vanishing_for_Schubert_varieties} to see that a weakly ordinary variety satisfying the condition $K_X \sim 0$ is globally $F$-split.  Then, using \cite[Theorem 1.1 and Theorem 1.2]{Ejiri_When_is_the_Albanese_morphism_an_algebraic_fiber_space_in_positive___characteristic?} we see the Albanese $X \to \Alb_X$ is in fact a relatively normal and $F$-split surjective algebraic fiber space.  Consequently, we know that the general fibers of $X \to \Alb_X$ are normal, and all the fibers are $F$-pure and hence reduced.

Moreover, since the condition of $F$-splitting is preserved under \'etale covers, this also implies that the augmented irregularity as defined in \eqref{eq:augmented_irregularity} is finite.  We may therefore take a Galois \'etale cover $Z \to X$ such that $\dim \Alb_Z = \wh{q}(X)$. We note that as $K_Z\sim 0$ and $Z$ is globally $F$-split, for $Z \to \Alb_Z$ also hold  the above features of $X \to \Alb_X$.  We shall prove that the Albanese morphism $Z \to \Alb_Z$ becomes a product after taking an \'etale cover and a further diagonalizable torsor over the base.

\subsubsection{General fibers are strongly $F$-regular}
 
The standard tools to control the behaviour of fibrations such as $Z \to \Alb_Z$ are the semi-positivity results for the relative canonical sheaves provided in characteristic $p$ in the paper of the first author \cite{Patakfalvi_Semi_positivity_in_positive_characteristics}.  The ubiquitous requirement for the application of the aforementioned results is the strong $F$-regularity of the general fibers of the investigated morphism.  Using the arguments above, until now we only managed to prove that the fibers are $F$-pure.  In order to improve the situation and get that the general fibers of $Z \to \Alb_Z$ are strongly $F$-regular it suffices to show that the Frobenius pullbacks $Z^e = Z \times_{F^e_{\Alb_Z}} \Alb_Z$ are strongly $F$-regular along the generic fibers of natural projections $Z^e \to \Alb_Z$, for every $e > 0$ \cite{Patakfalvi_Schwede_Zhang_F_singularities_in_families}.  The natural tool now is the theory of test ideals $\tau(Y) \subseteq \cO_Y$, associated to normal varieties $Y$, developed by Hochster and Huneke (see \cite{Hochster_Huneke_Tight_closure_and_strong_F-regularity} for the original work and \cite{Schwede_Tucker_A_survey_of_test_ideals} for a comprehensible survey).  The ideals control strong $F$-regularity of $Y$ in the sense that the condition $\tau(Y) = \cO_Y$ is satisfied exactly along the locus where $Y$ is strongly $F$-regular.  In our situation, to prove that $\tau(Z^e) = \cO_{Z^e}$ along the generic fiber of the projection we apply the transformation rule for test ideals under finite maps provided in \cite{Schwede_Tucker_On_the_behavior_of_test-ideals_under_finite_morphisms} to the relative Frobenius $F^e_Z \colon Z \to Z^e$.  We emphasize that, as required by the results of \emph{loc.cit}, $Z^e$ is normal along the generic fiber of the projection (see the middle of \autoref{ss:intro_albanese_morphism_irregularity}).  For the precise argument, we refer to \autoref{sec:singularities_general_fibres}.

\subsubsection{Flatness}

As a next step towards the proof, in \autoref{sec:flatness} we show that the Albanese morphism of $Z$ is in fact flat.  For this purpose, we mimic the characteristic zero arguments given in \cite[Theorem]{Lu_Tu_Zhang_Zheng_On_semistability_of_Albanese_maps}.  Our contribution here is mainly the realization that in characteristic $p$ there are appropriate semi-positivity results (see \autoref{sec:semi-positivity_engine}) to execute the strategy.  To sum up, our current arguments state that the morphism $Z \to \Alb_Z$, which we intend to prove becomes a product, is a flat algebraic fiber space with strongly $F$-regular general fiber.

\subsubsection{Restriction to curves}

We consequently proceed to the proof of the next approximation of the desired result, that is, we show that the family $Z \to \Alb_Z$ is isotrivial over a general complete intersection curve $T \subset \Alb_Z$ that goes through an arbitrarily fixed closed point $t_{\specialpt} \in \Alb_Z$ and a general fixed closed point $0 \in \Alb_Z$.  We set $V = T \times_{\Alb_Z} Z$.  Since the all fibers of $Z \to \Alb_Z$ are reduced, and the general fiber is normal and strongly $F$-regular, we see that $V$ is normal and the morphism $V \to T$ is a flat fibration with strongly $F$-regular general fiber.  Moreover, using the base change formulas for the relative canonical divisor we see that $K_{V/T} \sim (K_{Z/\Alb_Z})_{ | V} \sim 0$.  We note that in general this part of the argument requires a little care.  We provide the relevant base change results for the relative canonical in \autoref{sec:canonical_bundle_formulas}.  The details of the argument are provided as the part of the main proof given in \autoref{sec:isotriviality_of_the_Albanese}.

\subsubsection{Numerical flatness}
\label{ss:intro_numerical_flatness}

In order to prove that the morphism $f \colon V \to T$ is isotrivial we first show that there exists an appropriate $f$-ample divisor on $L$ on $V$ such that the relative section sheaves $f_*\cO_V(mL)$, for $m>0$, satisfy certain notion of triviality called \emph{numerical flatness}.  A vector bundle $\cE$ on a projective smooth scheme $X$ is numerically flat if both $\cE$ and $\cE^\vee$ are nef.  We refer to \autoref{ss:numerically_flat_bundles} for a more detailed description of the notion.  

In our context, we prove that the sheaves $f_*\cO_V(mL)$ are numerically flat if
 $L$ is an $f$-ample divisor such that $L^{d+1} = 0$, where $d+1 = \dim V$.  For
 the detailed proof we refer to \autoref{sec:numerical_flatness}.  In this paragraph we give a brief description of the arguments.  

First of all, we show that an $f$-ample divisor $L$ satisfying $L^{d+1} = 0$ is in fact nef. In fact, it is  enough to show that $L + f^* \varepsilon H$ is nef for fixed ample divisor $H$ on $T$ and for every $\varepsilon>0$. The point is that  a Riemann-Roch computation shows that in this case there is an effective $\Gamma \sim_{\bQ} L + f^* 
\varepsilon H$, see \autoref{lem:asymptotic_RR_relative_ample}. Then, semi-positivity theory applied to $\varepsilon' \Gamma = (K_V + \Delta_V + \varepsilon' \Gamma) + (-K_V - \Delta_V)$ yields the above nefness using that $-K_V - \Delta_V \sim_{\bQ} 0$, see \autoref{thm:nef} and \autoref{lem:semipositivity}. Here $\Delta_V$ is a natural boundary on $V$ that makes the linear equivalence $K_{Z/\Alb Z}|_V \sim_{\bQ} K_V + \Delta_V$, see \autoref{sec:canonical_bundle_formulas}. Having shown the nefness of $L$, the nefness of $f_* \sO_V(mL)$ follows from standard semi-positivity theory again, see \autoref{thm:K_trivial_numerically_flat}. 

Second, we show that $f_* \sO_V(mL)$ is anti-nef. So, assume the contrary.  According to \autoref{lem:nef_L_min_max}, this is equivalent to the lowest piece $\sE$ of the Harder-Narasimhan filtration  of $F^{l,*} f_* \sO_V(mL)$ having positive degree for every $l \gg 0$. Now, by choosing $l \gg0$  we may assume that $\sE$ is strongly semi-stable, and then $\sE^{\otimes r}$ is semi-stable for every $r>0$ \cite[Thm 6.1]{Langer_Semistable_sheaves_in_positive_characteristic}. The main idea is that by going to $r \gg 0$, via the multiplication map $\sE^{\otimes r} (-t)  \to (f_*\sO_X(rmL))(-t)$ this yields a section in $H^0(T, f_*\sO_X(rmL)(-t)) \cong H^0(X, rmL - X_t)$ that should not exist as $L^{d+1} = 0$, see \autoref{thm:dual_nef}.

\subsubsection{Isotriviality for finite fields}
\label{ss:intro_isotriviality_over_finite_fields}

In characteristic $p>0$, the numerical flatness turns out to be a particularly strong notion if the underlying variety $X$ is defined over the algebraic closure of the finite field.  More precisely, using the results of Langer and the classical theorem of Lange and St\"uhler one can prove that a numerically flat bundle defined on a variety $X/\overline{\FF}_q$ is trivializable on the cover $Y \to X$ which is a composition of a finite \'etale morphism and a power of the Frobenius (see \autoref{lem:pullback_trivializes}).   Assuming that the variety $V$ is defined over the algebraic closure of a finite field and taking such a cover $\tau \colon S \to T$ for the bundle $f_*\cO_V(mL)$, where $m$ is chosen so that the natural multiplication maps $\Sym^d f_*\cO_V(mL) \to f_*\cO_V(dmL)$ are surjective, we see that the relative section ring 
\[
R_{V_S/S}(mL_S) = \bigoplus_{d \in \NN} f_{S*}\cO_{V_S}(dmL_S) \expl{\isom}{flat base change} \bigoplus_{d \in \NN} \tau^*f_*\cO_V(dmL)
\] 
consists of numerically flat bundles which are quotients of trivial bundles.  Such bundles are in fact trivial themselves and hence, since $S$ is a projective curve, $R_{V_S/S}(mL_S)$ comes as pullback of a ring defined over the base field.  This  implies that $V_S \to S$ is a product family, and hence gives isotriviality over curves $T$ if the initial variety $X$ is defined over the algebraic closure of a finite field.  The precise statements are presented in \autoref{sec:isotriviality_finite_fields}.  Similar arguments actually lead to the proof of \autoref{thm:isotriviality}.

\subsubsection{Reduction to finite fields}

In order to reduce to the case where the base field is the algebraic closure of a finite field, we use the spreading out technique and the base change properties for the suitable relative polarized isomorphism scheme.  We first observe that above isotriviality result implies that the natural map 
\[
\Isom_T\left((V,mL),(V_t \times_k T,mL_t \times_k T)\right) \to T,
\]
where $t \in T$ is a $k$-rational base point, is surjective if the base field $k$ is the algebraic closure of a finite field.  In order to get a similar statement in general, for arbitrary perfect base field $k$, we take a spreading out $(\cV \to \cT,\cL,\sigma \colon \Spec(R) \to \cT)$ over a finitely generated $\FF_q$-algebra $R$ of the morphism $V \to T$ along with the divisor $L$ and a choice of a base point $t \in T$.  We consider the relative isomorphism scheme 
\[
\cI = \Isom_{\cT}\left((\cV,m\cL),(\cV_{\sigma} \times_{\Spec(R)} \cT, m\cL_{\sigma} \times_{\Spec(R)} \cT)\right) \to \cT.
\]
By base change property of isomorphism schemes, we see that the result of the previous section implies that the  morphism $\cI \to \cT$ defined over $R$ is surjective when restricted to every closed point of $\Spec(R)$.  By a standard scheme theoretic argument, this yields surjectivity at the geometric generic point, and hence the required isotriviality even in the polarized setting.  We note that the number $m$ showing up in \autoref{ss:intro_isotriviality_over_finite_fields} can be chosen in advance before the spreading out so that it yields surjectivity of the multiplication map for every finite field reduction.

\subsubsection{From isotriviality over curves to the isotriviality over $\Alb_Z$ and to the splitting over a flat cover}

First, we observe that the choice of a line bundle $L$  in \autoref{ss:intro_numerical_flatness} and the number $m$ in \autoref{ss:intro_isotriviality_over_finite_fields} can in fact be performed uniformly, see \autoref{thm:decomposition_theorem_alternative} and \autoref{thm:decomposition_theorem}.
 Then, we consider the isomorphism scheme
\[
I = \Isom_{\Alb_Z}\left((Z,mL_Z),(\Alb_Z \times_k Z_0,\Alb_Z \times_k mL_0)\right)
\]
where $0 \in \Alb_Z$ is a fixed $k$-rational base point and $L_0 = L_{Z|Z_0}$. 
By the base-change properties of the isomorphism scheme, and the above explained isotriviality over the curves $T \subseteq \Alb_Z$, the image of $I\to \Alb_Z$ contains every such $T$. However, as $t_{\specialpt} \in T$ was arbitrarily fixed, this means that $\pi : I \to \Alb_Z$ is surjective. 
Additionally,  by using again the base-change properties of the isomorphism scheme, we see that the morphism $Z \to \Alb_Z$ becomes a product (even in the polarized way) after the base change along $\pi$. The formal version of the argument is provided in the proof of \autoref{thm:decomposition_theorem_alternative}.

We remark that in the logarithmic setting the actual argument is quite delicate.  We provide the details of the necessary base change results for the logarithmic version of the isomorphism scheme in \autoref{sec:Isom_scheme}.

\subsubsection{Finiteness of automorphism groups}

The above argument does not say anything about the structure of the trivializing cover $\pi \colon I \to \Alb_Z$.  In order to rectify this situation, we observe that $\pi$ is in fact a torsor under the polarized automorphism scheme $G = \Aut(Z_0,L_0)$ of the fiber $Z_0$.  By the previous considerations we see that $Z_0$ is a strongly $F$-regular, Gorenstein variety satisfying $K_{Z_0} \sim 0$.  In \autoref{sec:finite_automorphisms} we show that polarized automorphisms schemes of such varieties are in fact finite group schemes.  The proof is by contradiction.  Assuming otherwise, we infer that $G$ admits a subgroup isomorphic to the additive or multiplicative group acting on $Z_0$ with generically trivial stabilizers.  Using the classical observation of Rosenlicht \cite[Theorem 2 and 10]{Rosenlicht_Some_basic_theorems_on_algebraic_groups} this easily implies that $Z_0$ is ruled which gives a contradiction using a simple argument based on Kodaira dimension.  We emphasize that the last part of argument requires our strong bounds on singularities of $Z_0$.

\subsubsection{Nori fundamental group scheme} 

We are now ready to finish the proof.  In the previous section, we showed that the trivializing morphism $I \to \Alb_Z$ is in fact a torsor under a finite flat groups schemes over $k$.  Such objects were extensively studied by Nori in his works concerning generalizations of \'etale fundamental groups.  In particular, in \cite[Proposition]{Nori_Fundamental_Group_Scheme_Of_An_Abelian_Variety} it is proven that the reduced part of every torsor under a finite flat group scheme over an abelian variety $A$ is dominated by the $A[n]$-torsor given by the multiplication map $[n] \colon A \to A$, for some $n \in \NN$.  Applying this to our situation, this means the morphism $I_{\red} \to \Alb_Z$ is covered by the multiplication map $[n] \colon \Alb_Z \to \Alb_Z$, which in turn implies that $Z \to \Alb_Z$ becomes a product after taking a base change under $[n]$.  We conclude the proof by observing that $[n]$ is in fact a composition of an \'etale morphism and a diagonalizable torsor because $\Alb_Z$ is $F$-split.  The details of the above argument are provided in \autoref{sec:bogomolov_beauville_actual_proofs}.


\subsection{Acknowledgements}

The authors would like to thank Piotr Achinger, Fabio Bernasconi, Javier Carvajal-Rojas, J\'anos Koll\'ar, Adrian Langer, Max Lieblich, Yuya Matsumoto, Karl Schwede and Burt Totaro for useful conversations and remarks.  We are also grateful to the referee for many insightful comments.

During the work on the article the authors were supported by grant \#200021/169639 of the Swiss National Science Foundation.

This material is partially based upon the work of the first author supported by the National Science Foundation under Grant No. DMS-1440140 while the author was in residence at the Mathematical Sciences Research Institute in Berkeley, California, during the Spring of 2019 semester.

\section{General preliminaries}

In this section we gather some preliminary results required in the following considerations.


\subsection{Base field}

The field $k$ is perfect, and apart from \autoref{sec:Isom_scheme} and \autoref{sec:appendix}, of characteristic $p>0$.

\subsection{Basic notation and definitions}

A \emph{variety} is a separated, integral scheme of finite type over our fixed base field $k$. 
A \emph{fibration} means a surjective proper morphism $f \colon X \to T$ between varieties, such that the natural map $\sO_T \to f_* \sO_X $ is an isomorphism.  An open set $U$ of a Noetherian scheme $X$ is \emph{big} if $\codim_X (X \setminus U ) \geq 2$.

For a variety $X$ defined over a field $k$, by $H^i(X_\et,\QQ_\ell)$ (resp. $H^i_c(X_\et,\QQ_\ell)$) we denote the $\ell$-adic cohomology (resp. $\ell$-adic cohomology with compact support) given by the formula $H^i(X_\et,\QQ_\ell) = \left(\underset{n \to \infty}{\lim} H^i_\et(X_{\overline{k}},\ZZ/\ell^n\ZZ)\right) \otimes \QQ_\ell$.

\subsubsection{$F$-singularities}

The following definition is a summary of the classical definitions concerning $F$-singularities given for example in \cite[Definition 4.1 and Definition 4.2]{Patakfalvi_Schwede_Tucker_Positive_characteristic_algebraic_geometry}.

\begin{definition}
\label{def:F_singularities}
Let $X$ be a normal variety, and let $\Delta$ be an effective $\QQ$-Weil divisor.
\begin{enumerate}
\item The pair $(X,\Delta)$ is globally $F$-split if there exists an integer $e>0$ for which the natural map $\cO_X \to F^e_*\cO_X(\lceil (p^e - 1)\Delta \rceil)$ splits in the category of $\cO_X$-modules.
\item The pair $(X,\Delta)$ is sharply $F$-pure if it can be covered by globally $F$-split open subvarieties. 
\item The pair $(X,\Delta)$ is globally  $F$-regular if for every divisor $D$ and for every integer $e \gg 0$ the natural map $\cO_X \to F^e_*\cO_X(\lceil (p^e - 1)\Delta \rceil + D)$ splits in the category of $\cO_X$-modules.
\item The pair $(X,\Delta)$ is strongly $F$-regular if it can be covered by globally  $F$-regular open subvarieties.
\end{enumerate}
\end{definition}

\begin{remark}
It turns out that strong $F$-regularity can be verified using the concept of test ideal sheaves (see \cite{Hochster_Huneke_Tight_closure_and_strong_F-regularity,Schwede_Tucker_A_survey_of_test_ideals}).  For brevity, we usually refer to them as simply test ideals.  We recall that for a pair $(X,\Delta)$ the corresponding test ideal sheaf $\tau(X,\Delta)$ is defined (see \cite[Definition 6.22]{Patakfalvi_Schwede_Tucker_Positive_characteristic_algebraic_geometry}) as the minimal ideal subsheaf $\cJ \subset \cO_X$ such that for all $e \geq 1$ and all local (defined on an affine open subset) homomorphisms $\phi \colon F^e_*\cO_X(\lceil (p^e - 1)\Delta \rceil \to \cO_X$ the condition $\phi(F^e_*\cJ) \subset \cJ$ is satisfied.  In our consideration, we use the fact that a pair $(X,\Delta)$ is strongly $F$-regular if and only if the test ideal $\tau(X,\Delta)$ is isomorphic to $\cO_X$ (see \cite[Theorem 3.19]{Schwede_Tucker_A_survey_of_test_ideals} for a divisor free statement).  We will use the transformation rules for test ideals under finite maps, proven in \cite{Schwede_Tucker_On_the_behavior_of_test-ideals_under_finite_morphisms}, in order to analyze the singularities of the Albanese morphism (see \autoref{sec:singularities_general_fibres}).
\end{remark}

\subsection{Reflexive sheaves and reflexive operations}

If $\sF$ is a  coherent sheaf on a normal variety $X$, then $\Sym^{[a]}(\sF):= (\Sym^a(\sF))^{**}$ is the reflexive symmetric power. The most important thing to remember is that all types of reflexive operations are determined in codimension $1$, that is, in this particular case if  $\sT \subseteq \sF$ is the maximal torsion subsheaf, and $\iota : U \hookrightarrow X$ is the locus where $\factor{\sF}{\sT}$ is locally free, then $\Sym^{[a]}(\sF)= \iota_* \left(\Sym^a\left(\left.\factor{\sF}{\sT}\right|_U\right)\right)$ \cite{Hartshorne_Stable_reflexive_sheaves}. We define similarly reflexive pullbacks $f^{[*]} \sF:= (f^* \sF)^{**}$, where $f : Y \to X$ is a morphism of finite type, or reflexive tensor products $\sF [ \otimes ] \sG:= (\sF \otimes \sG )^{**}$, where $\sG$ is also a coherent sheaf on $X$. For the latter reflexive operations similar extension property holds in terms of $\iota_*$ as for the reflexive symmetric product.

In general we use the basic statements about reflexive coherent sheaves on normal varieties \cite[page 124-129]{Hartshorne_Stable_reflexive_sheaves} without giving a precise reference each time. 

\subsection{General, generic and geometric generic fibers}

If $f : X \to T$ is a fibration, then the generic  fiber is the fiber $X_\eta$ over the generic point $\eta \in T$.  Similarly, the geometric generic fiber is $X_{\overline{\eta}}$, where $\overline{\eta} = \Spec \overline{k(\eta)}$. On the other hand, we say that a property holds for a general fiber, if there is a dense open set $U \subseteq T$ such that the given property holds for all $X_t$ where $t \in U$ is a closed point.  We note that we only work over a perfect not necessarily algebraically closed field, and therefore there are points which are closed but not rational (i.e., the residue field of a point is a finite extension of the base field).  The general pattern is that, as $k$ is assumed to be perfect, singularity properties holds for general fibers if and only if they hold for the geometric generic fiber. See \cite[Prop 2.1]{Patakfalvi_Waldron_Singularities_of_General_Fibers_and_the_LMMP} for the  incarnation of this pattern for normality, reducedness and regularity. 

\subsection{Numerically flat bundles}
\label{ss:numerically_flat_bundles}

In the main part of the paper, we will use the theory of numerically flat bundles in characteristic $p$ provided by Langer in \cite{Langer_On_the_S-fundamental_group_scheme_I,Langer_On_the_S-fundamental_group_scheme_II}.

\begin{definition}
\label{def:num_flat_def}
Let $T$ be a smooth projective variety over a perfect field $k$.  We say that a vector bundle  $\sE$ on $T$ is \emph{numerically flat}, if both $\sE$ and $\sE^*$ are nef. 
\end{definition}

\begin{proposition}[{\cite[Proposition 5.1]{Langer_On_the_S-fundamental_group_scheme_I}}]
\label{prop:num_flat_equivalent_defs}
Let $T$ be a smooth variety over $k$, let $H$ be an ample divisor and let $\cE$ be a vector bundle.  The following conditions are equivalent:
\begin{enumerate}
	\item $\cE$ is numerically flat,
	\item $\cE$ is nef and of degree zero.
    \item $\cE$ is strongly $H$-semistable with $c_1(\cE).H^{n-1} = c_2(\cE).H^{n-2} = 0$.
\end{enumerate}
\end{proposition}

We recall that by definition strong $H$-semistability is equivalent to $H$-semistability of all Frobenius pullbacks.  We also note that by the above result numerically flat bundles are in fact strongly semistable with respect to any polarization.

\begin{lemma}
\label{lem:pullback_trivializes}
Suppose $\sE$ is a numerically flat vector bundle on a normal projective variety $T$ over $\bF_q$.  Then there is a finite cover $\tau \colon S \to T$ by a normal projective variety such that $\tau^* \sE \cong \sO_S^{\oplus \rk \sE}$.  The cover $\tau$ might be chosen to be a composition of a finite \'etale covering and an iteration of the  Frobenius morphism.
\end{lemma}

\begin{proof}
According to \cite[Thm 1.1]{Langer_On_the_S-fundamental_group_scheme_II}, there is a scheme $M$ of finite type over  $\bF_q$ parameterizing numerically flat vector bundles of rank $r$ over $T$. So, for each numerically vector bundle $\sF$ over $T$ correspond a few $\bF_q$-rational points of $M$. As $|M(\bF_q) | < \infty$, we obtain that there are integers $e\geq 0$ and $e'>0$ such that $\left(F^{e + e'}\right)^* \sE \cong \left(F^{e }\right)^* \sE$. Hence by the main result of \cite{Lange_Stuhler_Vektorbundel_auf_curven}, there is a finite \'etale cover $\rho \colon S \to T$ such that $\rho^* \left(F^{e }\right)^* \sE$ is trivial.  Set then $\tau:= F^e \circ \rho $.  
\end{proof}

\subsection{Albanese morphism}

For the arithmetic applications, we shall work over perfect non-necessarily algebraically closed field.  In this context the \emph{Albanese morphism} of a normal complete variety $X$ is an initial object in the category of morphism $X \to T$, where $T$ are torsors over abelian varieties.  It exists by the detailed discussion provided in the Mathoverflow answer \cite{MO_Albanese}.  More precisely, in the fifth paragraph the author claims that the Albanese is a torsor under the maximal abelian subvariety of the reduced part of the Picard scheme.  In particular, this behaves well with respect to base change over extensions of perfect fields (such extensions preserve reducedness). All these arguments are already present in Serre's original work \cite{Serre_Morphismes_universels_variete_d_Albanese}. More modern treatments are in \cite{Wittenberg_On_Albanese_torsors_and_the_elemntary_obstruction,Achter_Casalaina-Martin_Vial_A_complete_answer_to_Albanese_base_change_for_incomplete_varieties}.
\begin{proposition}
\label{prop:albanese_variety_basic_properties}
Let $X$ be a normal projective variety defined over a perfect field $k$.  Then the dimension of the Albanese variety $\dim \Alb_X$ is equal to $\frac{1}{2}H^1(X_{\et},\QQ_\ell)$.  
\end{proposition}

\begin{proof}
After base extension, we may assume that $X$ has a rational point.  Then the Albanese variety is isomorphic to the dual $\Pic^0(X)_{\rm red}^{\vee}$ (see \cite[Theorem 5.3]{Badescu_Algebraic_surfaces}), and consequently the statement follows from \autoref{prop:picard_variety_basic_properties}.
\end{proof}

\begin{proposition}
\label{prop:picard_variety_basic_properties}
Let $X$ be a normal projective variety defined over a perfect field $k$.  Then the identity component $\Pic^0(X)$ of the Picard scheme is projective.  Its reduced part $\Pic^0(X)_{\rm red}$ is an abelian variety of dimension $\frac{1}{2}  \dim_{\QQ_\ell} H^1(X_{\et},\QQ_\ell)$.
\end{proposition}

\begin{proof}
We need to prove $\Pic^0(X)$ is of dimension $\frac{1}{2}\dim_{\QQ_\ell} H^1(X_{\et},\QQ_l)$.  For this purpose, we consider an exact sequence of \'etale sheaves on $X_{\overline{k}}$:
\[
0 \ra \ZZ/\ell^N \ZZ \ra \GG_m \ra \GG_m \ra 1.
\] 
The long exact sequence of cohomology yields an isomorphism $H^1_{\et}(X_{\overline{k}},\ZZ/\ell^n \ZZ) \isom \Pic(X_{\overline{k}})[\ell^N]$.  As $\left|\Pic(X_{\overline{k}})[\ell^N] \setminus \Pic^0_{\red} (X_{\overline{k}})[\ell^N]\right|$ is bounded in $N$, and $\Pic^0_{\red} (X_{\overline{k}})[\ell^N] \cong \left( \bZ/\ell^N \bZ \right)^{2 \dim \Pic^0(X)}$, taking a limit with $N \to \infty$ we get the desired claim.
\end{proof}

\begin{corollary}
\label{cor:univ_homeo_alb_dim}
Let $X \to Y$ be a universal homeomorphism between normal schemes defined over a perfect field $k$.  Then the Albanese dimension of $X$ and $Y$ are equal.
\end{corollary}


\subsection{$F$-splittings of varieties and morphisms}
\label{sec:relatively_fsplit}

We shall need the following relative version of \autoref{def:F_singularities} (1).  Let $f \colon X \to S$ be a morphism of schemes over a perfect field $k$.  We recall that the $e$-th relative Frobenius morphism $F^e_{X/S}$ is defined by the following diagram:
\[
\xymatrix{
  X \ar[rd]^{F^e_{X/S}}\ar@/_1.2pc/[rdd]_{\pi} \ar@/^1.2pc/[rrd]^{F^e_X}& & \\
   & X^e_S \ar[r]\ar[d]_{\pi'}\ar@{}[dr]|-{\square} & X \ar[d]^{\pi} \\
   &  S \ar[r]_{F^e_S} & S.
}
\]

\begin{definition}[{\cite[Definition 5.1]{Ejiri_When_is_the_Albanese_morphism_an_algebraic_fiber_space_in_positive___characteristic?}}]\label{def:relative_F_split}
We say that a pair $(f,\Delta)$ is $F$-split if the map:
\[
\cO_{X^e_S} \to F^e_{X/S,*}\cO_X \to F^e_{X/S,*}\cO_X(\lceil (p^e - 1)\Delta \rceil).
\]  
admits a splitting in the category of $\cO_X$-modules.  In particular, a scheme $X$ defined in characteristic $p>0$ is $F$-split if the natural map $X \to \Spec(\FF_p)$ is $F$-split (equivalently, the natural map $\cO_X \to F_*\cO_X$ is split as a homomorphism of $\cO_X$-modules, cf. \autoref{def:F_singularities}).
\end{definition}

\noindent Suppose $f \colon X \to S$ is an $F$-split morphism, and let $T \to S$ be a morphism.  Taking a base change of a splitting we see that $f_T \colon X_T \to T$ is also $F$-split.  In particular, fibers of an $F$-split morphism are $F$-split \cite[Proposition 5.6]{Ejiri_When_is_the_Albanese_morphism_an_algebraic_fiber_space_in_positive___characteristic?}.

\begin{proposition}
\label{prop:technical_F_splittings_and_section}
Let $X$ be a normal scheme defined over a perfect field $k$.  Then for every Weil divisor $D$ there exists an isomorphism 
\[
\sHom_{\cO_X}(F^e_*\cO_X(D),\cO_X) \isom F^e_*\cO_X\left((1 - p^e)\left(K_X + \frac{D}{p^e - 1}\right)\right).
\]  
In particular, there is a bijection between the morphisms $F^e_*\cO_X(D) \to \cO_X$ and global sections of the reflexive sheaf $F^e_*\cO_X\left((1-p^e)\left(K_X + \frac{D}{p^e - 1}\right)\right)$. 
\end{proposition}

\begin{proof}
We first remark that the divisor $(1-p^e)\left(K_X + \frac{D}{p^e - 1}\right)$ in the statement is a Weil divisor with integer coefficients, and so the statement does make sense.  We may restrict our considerations to the smooth locus $X_{\rm sm}$ because $X$ is normal and hence all the sheaves in question are reflexive.   Using the Grothendieck duality we thus obtain a bijection:
\begin{align*}
\sHom(F^e_*\cO_{X_{\rm sm}}(D),\cO_{X_{\rm sm}}) \isom \sHom(F^e_*\cO_{X_{\rm sm}}(D),\omega_{X_{\rm sm}}) \otimes \omega_{X_{\rm sm}}^{-1} \expl{\isom}{Grothendieck duality and projection formula} F^e_*\cO_{X_{\rm sm}}((1 - p^e)K_{X_{\rm sm}} - D) 
\\ \isom  F^e_*\cO_{X_{\rm sm}}\left((1-p^e)\left(K_{X_{\rm sm}} + \frac{D}{p^e - 1}\right)\right).
\end{align*}
The final statement follows by taking global sections.
\end{proof}

We are now ready to recall the following results which we will use to understand the singularities of the fibers of the Albanese morphism.

\begin{proposition}[{\cite[Theorem 1.1 and Theorem 1.2]{Ejiri_When_is_the_Albanese_morphism_an_algebraic_fiber_space_in_positive___characteristic?}}]
\label{prop:albanese_of_fsplit}
Let $(X,\Delta)$ be a pair of a normal variety $X$ and an effective $\QQ$-Weil divisor $\Delta$.  Then $(X,\Delta)$ is $F$-split if and only if the Albanese morphism $X \to A$ is $F$-split with respect to $\Delta$ (see \autoref{def:relative_F_split}), and $A$ is ordinary.  Moreover, if Cartier index of $K_X + \Delta$ is coprime to $p$, the Albanese morphism is a surjective fibration.
\end{proposition}

\begin{remark}
One can address the Stein factorization property of the Albanese morphism using the following argument based on \cite{Hacon_Patakfalvi_Zhang_Birational_characterization_of_abelian_varieties_and_ordinary_abelian_varieties_in_characteristic_p_greater_than_0}.  Let $\pi \colon X \to Y$ be the Stein factorization of $X \to A$.  Taking pushforward of an $F$-splitting of $X$, we see that $Y$ is $F$-split and hence of non-positive Kodaira dimension.  Using \cite[Proposition 1.4]{Hacon_Patakfalvi_Zhang_Birational_characterization_of_abelian_varieties_and_ordinary_abelian_varieties_in_characteristic_p_greater_than_0} we consequently see that $\pi$ is separable and consequently $X$ is of Kodaira dimension zero.  Granted its surjectivity, the map $Y \to A$ is now a separable morphism of varieties of $F$-split varieties of Kodaira dimension zero and hence \'etale.  This implies that $Y$ is an abelian variety itself.
\end{remark}

\begin{proposition}\label{prop:fibres_are_regular_codimension_one}
Let $X$ be a normal variety $X$.  Moreover, let $S$ be a normal variety and $f \colon X \to S$ an $F$-split morphism.  Then the geometric generic fiber of $f$ is regular in codimension one.
\end{proposition}

\begin{proof}
The proof is the same as that of \cite[Proposition 5.6]{Ejiri_When_is_the_Albanese_morphism_an_algebraic_fiber_space_in_positive___characteristic?}. 
\end{proof}

\begin{proposition}[{\cite[Proposition 9]{Mehta_Ramanathan_Frobenius_splitting_and_cohomology_vanishing_for_Schubert_varieties}}]\label{prop:ordinary_is_fsplit}
Let $X$ be a normal Cohen--Macaulay variety with trivial canonical class defined over a perfect field of characteristic $p>0$.  Then $X$ is weakly ordinary if and only if it is $F$-split.
\end{proposition}

\begin{proof}
For the convenience of the reader we recall the proof.  First, we observe that the dualizing sheaf $\omega_X$ is $S_2$ and agrees with $\cO_X(K_X)$ on the regular locus, and hence $X$ is Gorenstein and $\omega_X \isom \cO_X$.  Consequently, the natural map $\cO_X \to F_*\cO_X$ is split if and only if the trace map given by the Grothendieck dual
\[
{\rm Tr}_X \colon F_*\omega_X \isom \cHom(F_*\cO_X,\omega_X) \to \omega_X
\]
is split.  Since $\omega_X$ is trivial, this happens if and only if the corresponding map on global sections is surjective.  Using Serre duality this is equivalent to $H^d(X,\cO_X) \to H^d(X,\cO_X)$ being injective, and hence bijective since $k$ is perfect.
\end{proof}


\subsection{The relative canonical bundle}
\label{sec:canonical_bundle_formulas}

For a surjective morphism $f : X \to T$ of varieties we define the relative canonical sheaf $\omega_{X/T}$ as $\omega_{X/T} := h^{-r} ( f^! \sO_T)$, where $r$ is the relative dimension $\dim X - \dim T$,  $f^!$ is the upper shriek functor of Grothendieck duality as defined in \cite{Hartshorne_Residues_and_duality}, and $h^{-r}(\_)$ means taking the $(-r)$-th cohomology sheaf. 

In general the philosophy is that whenever we may associate a divisor to $\omega_{X/T}$, we call it $K_{X/T}$. This philosophy is made precise below in two cases: when $f$ is equidimensional and $X$ and $T$ are normal, and when $T$ is smooth. 

\subsubsection{The case of a projective equidimensional morphism.} To be precise, here equidimensional means  that there is an integer $r>0$, such that all the fibers of $f$ have pure dimension $r$. In this case, note that if $U \subseteq T$ is any big open set, and $\iota : f^{-1} U \to X$ is the natural embedding (note that $f^{-1} U \subseteq X$ is also big by the equidimensional assumption), then $\omega_{X/T} \cong \iota_* \omega_{f^{-1}U/U}$. Indeed, this statement is local, so we may shrink arbitrarily the base, and then we may take a Noether normalization of $f$:
\begin{equation*}
\xymatrix{
X \ar@/^1.5pc/[rr]^f \ar[r]^{\tau} & \bP^n_T \ar[r]^g & T
}
\end{equation*}
Then
{\small
\begin{equation*}
\tau_* \omega_{X/T} 		
\explshift{5pt}{\cong}{Grothendieck duality}
\sHom_{\bP^n}\left(\tau_* \sO_X,  \omega_{\bP^n_T/T}\right) 
\explparshift{160pt}{60pt}{\cong}{$\omega_{\bP^n_T/T} \cong \pr_T^* \omega_{\bP^n} \cong j_* \omega_{\bP^n_U/U}$, where \\$j : g^{-1}U \to \bP^n_T$ is the natural embedding and $\pr_T: X \times T \to T$ is the projection}
j_* \sHom_{g^{-1}U}\left(\tau_* \sO_X,  \omega_{\bP^n_T/T}\right)
\explshift{75pt}{\cong}{Grothendieck duality}
j_* \left( \tau|_{f^{-1}U} \right)_* \omega_{f^{-1}U/U} \cong \tau_* \iota_* \omega_{f^{-1}U/U}.
\end{equation*}
}

\subsubsection{The case of a projective equidimensional morphism, normal base and normal total space.} 
\label{sec:relative_canonical_bundle_normal_spaces}
Assume now additionally that $X$ and $T$ are normal. Let $U \subseteq T$ be the big open set which is the intersection of the regular locus of $T$ and the flat locus of $f$. Then, 	  the above isomorphism $\omega_{X/T} \cong \iota_* \omega_{f^{-1}U/U}$ yields that 
\begin{equation*}
\omega_{X/T} \cong \iota_* \omega_{f^{-1}U/U} 
\expl{\cong}{by \cite[Lem 2.4]{Patakfalvi_Semi_negativity_of_Hodge_bundles_associated_to_Du_Bois_families} using that $U$ is Gorenstein}
 \iota_* \left( \omega_{f^{-1} U} \otimes f^* \omega_U^{-1} \right)  \cong \iota_* \sO_X( K_{f^{-1}U} - f^{-1} K_U ) \cong \sO_X(K_X - f^* K_T),
\end{equation*} 
where the pullback means the usual pullback for equidimensional surjective morphism between normal varieties (i.e., doing it over the flat locus which is big in the target, and then extending uniquely). 

\subsubsection{The case of smooth base and normal total space.}
\label{sec:relative_canonical_bundle_smooth_base}

If $T$ is smooth, $X$ is normal, but $f$ possibly non-equidimensional, $\omega_{X/T}$ is still a divisorial sheaf. Indeed, in this case $\omega_{X/T} \cong \omega_X \otimes f^* \omega_T^{-1}$ by \cite[Lem 2.4]{Patakfalvi_Semi_negativity_of_Hodge_bundles_associated_to_Du_Bois_families}. Hence, as $\omega_X$ is reflexive \cite[Cor 5.69]{Kollar_Mori_Birational_geometry_of_algebraic_varieties} we obtain that so is $\omega_{X/T}$ and hence it is divisorial. By the above isomorphism we also have $K_{X/T} = K_X - f^* K_T$. 

\subsubsection{Adjunction when the base is smooth and the total space is normal.}
\label{sec:relative_canonical_bundle_adjunction_smooth_base}

In the case of \autoref{sec:relative_canonical_bundle_smooth_base},  let $S \subseteq T$ be a smooth one codimensional subvariety, and assume that the scheme theoretical preimage $f^{-1} S$ is integral at the codimension $1$ points of $X$. That is, we assume that the divisorial pullback $Y:=f^*S$ is integral. Let $\tau \colon Z \to Y$ be the normalization, and set $g:= f|_Y$. Let $\Delta$ be an effective $\bQ$-divisor on $X$ such that $Y \not\subseteq \Supp \Delta$ and $K_X + \Delta$ is $\bQ$-Cartier. Finally, let $\Delta_Z$ be the different of $\Delta + Y$, that is, in particular $K_Z + \Delta_Z \sim_{\bQ} K_X + Y + \Delta|_Z$.
Then,  we have
\begin{multline}
\label{eq:relative_canonical_bundle_adjunction_smooth_base}
K_{Z/S} + \Delta_Z 
\expl{=}{\autoref{sec:relative_canonical_bundle_smooth_base}}
 K_Z  + \Delta_Z - (g \circ \tau)^* K_S
\sim_{\bQ}
 K_X + \Delta +  Y|_Z - (g \circ \tau)^* (K_T +S|_S)) 
\\ =
 K_X + \Delta +  Y - f^* (K_T +S)|_Z 
 \expl{=}{$Y = f^* S$ as divisors}
K_X  + \Delta - f^* K_T|_Y = K_{X/T}+ \Delta|_Z
\end{multline}


\section{Semi-positivity engine}
\label{sec:semi-positivity_engine}

\emph{The base field in this section and in \autoref{sec:flatness} and \autoref{sec:numerical_flatness} is perfect and of characteristic $p>0$.}   In some of the proofs we claim that we may assume that the base field is algebraically closed.  This is justified since the extension from a perfect field to one of its algebraic closures preserves all positivity and singularity properties. 

We will use the following technical statement multiple times throughout the article. 

Recall that pseudo-effectivity of non-$\bQ$-Cartier divisors is defined as follows: a divisor $D$ on a normal projective variety $X$ is pseudo-effective if for any (resp. all) ample divisor $H$ and any $\varepsilon>0$, $D + \varepsilon H$ is $\bQ$-effective. In particular for a divisor $D$ on $X$:
\begin{itemize}
\item if $C$ is a general element of a (generically irreducible) covering family of curves on $X$, such that $C$ is contained in the $\bQ$-Cartier locus of $D$, then $C \cdot D  \geq 0$, and
\item if there is a $\bQ$-effective Weil divisor $E$ such that $D + \varepsilon E$ is $\bQ$-effective for every rational number $\varepsilon>0$, then $D$ is pseudo-effective. 
\end{itemize}

\begin{theorem}
\label{lem:semipositivity}
Assume we are in the following situation:
\begin{enumerate}
\item  $f: X \to T$ is an equidimensional fibration between normal projective varieties with normal general fibers,
\item   $U \subseteq T$ is a non-empty open set, 
\item $\Gamma \geq 0$ be a $\bQ$-divisor such that $K_{X/T} + \Gamma|_{f^{-1}U}$ is $\bQ$-Cartier,
\item  $L$ is a nef  $\bQ$-Cartier $\bQ$-divisor such that $K_{X/T} + \Gamma + L$ is $f$-nef over $U$, and
\item     either
\begin{enumerate}
\item \label{itm:semipositivity:SFR}  $(X_t, \Gamma_t)$ is strongly $F$-regular for every closed point $t \in U$, or
\item \label{itm:semipositivity:SFP}  $(X_t, \Gamma_t)$ is sharply $F$-pure  for every closed point $t \in U$ with the Cartier index of $K_{X/T} + \Gamma|_U$ prime-to-$p$. 
\end{enumerate}
\end{enumerate}
 Then $K_{X/T} + \Gamma+L$ is pseudo-effective.

{\scshape Case $(*)$:}
If $T$ is a curve, and $K_{X/T} +  \Gamma+L$ is $f$-nef $\bQ$-Cartier (so globally, not only over $U$), then $K_{X/T} + \Gamma + L$ is not only pseudo-effective but also nef. 
\end{theorem}

\begin{proof}
By base extension we may assume that $k$ is algebraically closed.
We do this because we  reduce the statement to \cite[Thm 5.1]{Ejiri_Weak_positivity_theorem_and_Frobenius_stable_canonical_rings_of_geometric_generic_fibers}, which is stated over an algebraically closed base field. There are two main issues during this reduction. These issues and our strategy to deal with them are the following:
\begin{enumerate}[label=(\alph*)]
\item \label{itm:semipositivity:issue_nef} One issue is that there is no $L$ in \cite[Thm 5.1]{Ejiri_Weak_positivity_theorem_and_Frobenius_stable_canonical_rings_of_geometric_generic_fibers}.  We solve this by the following trick: we perturb $L$ with a little ample $\varepsilon H$ so that it becomes ample, and then by choosing a general effective $\bQ$-divisor $\bQ$-linearly equivalent to $L + \varepsilon H$ we merge the perturbed $L + \varepsilon H$ into the boundary. The key here is that this process does not make the singularities on the general fibers worse by the results of \cite{Schwede_Zhang_Bertini_theorems_for_F-singularities}. 
\item \label{itm:semipositivity:issue_index} The other main issue is that one has to guarantee two  conditions of \cite[Thm 5.1]{Ejiri_Weak_positivity_theorem_and_Frobenius_stable_canonical_rings_of_geometric_generic_fibers}. First, \cite[Thm 5.1]{Ejiri_Weak_positivity_theorem_and_Frobenius_stable_canonical_rings_of_geometric_generic_fibers} requires the Weil-index to be prime-to-$p$, and second  $(ii)$  of \cite[Thm 5.1]{Ejiri_Weak_positivity_theorem_and_Frobenius_stable_canonical_rings_of_geometric_generic_fibers} requires that $S^0\left(X_{\overline{\eta}}, \Gamma_{\overline{\eta}}; m \left( K_{X_{\overline{\eta}}}+ \Gamma_{\overline{\eta}} \right) \right) = H^0 \left(  m \left( K_{X_{\overline{\eta}}}+ \Gamma_{\overline{\eta}} \right) \right) $ for every divisible enough integer $m > 0$, where $\eta$ is the generic point of $T$.
These two conditions are closely related because according to  \cite[Prop 2.23]{Patakfalvi_Semi_positivity_in_positive_characteristics}, the second condition is satisfied if the index is not divisible by $p$, assuming the perturbation described in point \ref{itm:semipositivity:issue_nef}. 

We deal with these two requirements differently in the two cases of singularity assumptions. In the case of assumption \autoref{itm:semipositivity:SFP}, the above requirements are almost satisfied except there can be non prime-to-$p$ Weil index away from the general fibers. This is easy to solve  by removing the vertical components of $\Gamma$ from the boundary.

In the case of assumption \autoref{itm:semipositivity:SFR} however we have to work more.  We use  \cite[Lem 3.15]{Patakfalvi_Semi_positivity_in_positive_characteristics} in this case to perturb  to the (even Cartier-)index prime-to-$p$ situation. 
\end{enumerate}
Below we work out the details of the argument. 

First, we need to separate the vertical and the horizontal parts of $\Gamma$ in the case of assumption \autoref{itm:semipositivity:SFP}, because there is no way we can guarantee  in that case that the Weil index of $\Gamma$ in the vertical part is prime-to-$p$, which is an assumption of \cite[Thm 5.1]{Ejiri_Weak_positivity_theorem_and_Frobenius_stable_canonical_rings_of_geometric_generic_fibers}. 
So, in the respective cases of singularity assumptions we set
\begin{enumerate}[wide=10pt]
\item[\autoref{itm:semipositivity:SFR}:] $\Gamma_h:=\Gamma$ and $\Gamma_v:= 0$, and
\item[\autoref{itm:semipositivity:SFP}:] $\Gamma_h$ is the horizontal part, and $\Gamma_v$ is the vertical part.
\end{enumerate}
The upshot is that in either cases $\Gamma = \Gamma_h + \Gamma_v$ and furthermore in case \autoref{itm:semipositivity:SFP} we already have that the Weil index of $\Gamma_h$ is prime to $p$, by the assumption that $K_{X/T} + \Gamma|_U$ has prime-to-$p$ Cartier index.

Let $H$ be an ample divisor on $X$. Then $L + \varepsilon  H$ is ample for every $\varepsilon>0$. So, for every fixed  $0<\varepsilon\in \bQ$ the denominator of which is prime-to-$p$,  there is a prime-to-$p$ integer $m>0$ such that $|m(L+ \varepsilon  H)|$ is very ample. Choose a general $D \in |m(L+ \varepsilon  H)|$. Then for a general $t \in T$ we have that $D_t \in |mL_t|$ is also general. In particular, in the case of the respective singularity assumptions we have that:
\begin{equation}
\label{eq:semipositivity:SFP_SFR_after_L}
\begin{minipage}{340pt}
\begin{enumerate}
\item $\left(X_t, \left(\Gamma_h\right)_t + \frac{1}{m}D_t \right)$ is also strongly $F$-regular according to \cite[Corollary 6.10 (iii)]{Schwede_Zhang_Bertini_theorems_for_F-singularities}, and
\item $\left(X_t, \left(\Gamma_h\right)_t + \frac{1}{m}D_t \right)$ is also sharply $F$-pure with Gorenstein index prime-to-$p$ according to \cite[Corollary 6.10 (i)]{Schwede_Zhang_Bertini_theorems_for_F-singularities}.
\end{enumerate}
\end{minipage}
\end{equation}
The above perturbation is adequate for solving issue \ref{itm:semipositivity:issue_nef}. So, set $\Gamma_{\varepsilon}= \frac{1}{m}D$. Besides the singularity properties stated in \autoref{eq:semipositivity:SFP_SFR_after_L}, the only feature of  $\Gamma_{\varepsilon}$ important to remember  is that 
\begin{equation}
\label{eq:semipositivity:Gamma_epsilon_Q_equiv}
\Gamma_\varepsilon \sim_{\bQ} L + \varepsilon H.
\end{equation}
To also deal with issue \ref{itm:semipositivity:issue_index} in the case of assumption \autoref{itm:semipositivity:SFR}, we add a further divisor $\Gamma_{e, \varepsilon}$ to the boundary. This is defined as follows in the respective cases of singularity assumptions:
\begin{equation}
\label{eq:semipositivity:Gamma_e_epsilon}
\begin{minipage}{340pt}
\begin{enumerate}
\item  Choose $E \in |K_{X/T} + a H|$ for some fixed $a>0$ independent of $\varepsilon>0$. Set then  $\Gamma_{e, \varepsilon}:= \frac{1}{p^e -1}(E + \Gamma_h + \Gamma_\varepsilon)$ for some $e \gg0$ (warning: $K_{X/T} + a H$ is not Cartier, so no general choice of $E$ is possible, it must go through all the points where $K_{X/T} + a H$ is not Cartier).
\item  Set $\Gamma_{e, \varepsilon}:=0$.
\end{enumerate}
\end{minipage}
\end{equation}
Finally,   in either cases set $\Lambda:= \Gamma_h + \Gamma_{\varepsilon} + \Gamma_{e, \varepsilon}$. Summing up in the case of assumption \autoref{itm:semipositivity:SFP} we obtain: 
\begin{equation}
\label{eq:semipositivity:full_perturbation_SFP}
K_{X/T} + \Lambda + \Gamma_v 
\expl{\sim_{\bQ}}{definition of $\Lambda$}
 K_{X/T} + \Gamma_h + \Gamma_{\varepsilon} + \Gamma_{e, \varepsilon} + \Gamma_v
\expl{\sim_{\bQ}}{$ \Gamma= \Gamma_h + \Gamma_v$, \autoref{eq:semipositivity:Gamma_epsilon_Q_equiv} and \autoref{eq:semipositivity:Gamma_e_epsilon}}
 K_{X/T}  + \Gamma + L + \varepsilon H,
\end{equation}
 and in the case of assumption \autoref{itm:semipositivity:SFR} we obtain
\begin{multline}
\label{eq:semipositivity:full_perturbation_SFR}
K_{X/T} + \Lambda + \Gamma_v 
\expl{\sim_{\bQ}}{definition of $\Lambda$}
 K_{X/T} + \Gamma_h + \Gamma_{\varepsilon} + \Gamma_{e, \varepsilon} + \Gamma_v
\expl{\sim_{\bQ}}{$ \Gamma= \Gamma_h + \Gamma_v$, \autoref{eq:semipositivity:Gamma_epsilon_Q_equiv} and \autoref{eq:semipositivity:Gamma_e_epsilon}}
\\ \Gamma_v + \frac{p^e}{p^e-1} (K_{X/T} + \Gamma_h + L)  + \left( \varepsilon \frac{p^e}{p^e-1} + a \frac{1}{p^e -1}\right) H
\expl{\longrightarrow}{$\varepsilon\to0, e \to \infty$}
K_{X/T} + \Gamma + L 
\end{multline}
In particular, in either case, it is enough to show that $K_{X/T} + \Lambda + \Gamma_v$ is pseudo-effective (resp. nef in {\scshape case $(*)$}) for fixed $\varepsilon =\frac{1}{r}$, where $r>0$ is divisible enough, and $e \gg 0$, where $e$ is chosen to be very big after choosing a very divisible $r$. So, from now we change our goal to prove the pseudo-effectivity (resp. nefness) of $K_{X/T} + \Lambda + \Gamma_v$ with the above fixed choices of $\varepsilon$ and $e$.

The next step is to apply \cite[Thm 5.1]{Ejiri_Weak_positivity_theorem_and_Frobenius_stable_canonical_rings_of_geometric_generic_fibers}
 to $(X, \Lambda)$. So, we have to verify the requirements of \cite[Thm 5.1]{Ejiri_Weak_positivity_theorem_and_Frobenius_stable_canonical_rings_of_geometric_generic_fibers}
 one by one. Let us start with the Weil index being prime-to-$p$. In the case of assumption \autoref{itm:semipositivity:SFP}, $\Lambda = \Gamma_h + \frac{1}{m} D$, where $m$ is prime-to-$p$ and $\Gamma_h$ has prime-to-$p$ Weil index. So,  the Weil index of $\Lambda$ being prime-to-$p$ is satisfied in this case. 

In the case of assumption \autoref{itm:semipositivity:SFR}, we achieved this requirement by adding $\Gamma_{e, \varepsilon}$ into $\Lambda$. In fact, in this case even the Cartier index of $K_{X/T} + \Lambda$ is prime-to-$p$ and also $(X_t, \Lambda_t)$ has strongly $F$-regular singularities for $t \in T$ general closed point, by \autoref{eq:semipositivity:SFP_SFR_after_L} and \cite[Lem 3.15]{Patakfalvi_Semi_positivity_in_positive_characteristics}. To be precise, to go from a single $t \in T$ here one has to invoke \cite[Thm B]{Patakfalvi_Schwede_Zhang_F_singularities_in_families} as well. 

Additionally, in both cases we have  $\left(X_{\overline{\eta}},\Lambda_{\overline{\eta}} \right)$ is sharply $F$-pure with  Cartier index not divisible by $p$ (by the previous paragraph in the case of assumption \autoref{itm:semipositivity:SFR} and by \autoref{eq:semipositivity:SFP_SFR_after_L} in the other case).
In particular,  by \cite[Prop 2.23]{Patakfalvi_Semi_positivity_in_positive_characteristics}, in both cases, if $\eta$ is the generic point of $T$, for every integer $m>0$ divisible enough we have
\begin{equation*}
S^0\left( X_{\overline{\eta}},\Lambda_{\overline{\eta}}; m\left( K_{X_{\overline{\eta}}} +\Lambda_{\overline{\eta}} \right) \right) =H^0\left( X_{\overline{\eta}}, m\left( K_{X_{\overline{\eta}}} +\Lambda_{\overline{\eta}} \right) \right) .
\end{equation*}
Furthermore,  $f$ is separable, because $f_* \sO_X \cong \sO_T$ and the geometric generic fiber of $h$ is  normal. So, all condition of \cite[Thm 5.1]{Ejiri_Weak_positivity_theorem_and_Frobenius_stable_canonical_rings_of_geometric_generic_fibers}
are satisfied (the finite generation comes out of being ample on the general fiber, by \autoref{eq:semipositivity:full_perturbation_SFP} and \autoref{eq:semipositivity:full_perturbation_SFR}). 
In particular, we obtain that $f_* \sO_X (m(K_{X/T} + \Lambda))$ is weakly positive for every $m$ divisible enough. To  be precise, \cite[Thm 5.1]{Ejiri_Weak_positivity_theorem_and_Frobenius_stable_canonical_rings_of_geometric_generic_fibers} states that $f_* \sO_X (m(K_{X} + \Lambda)) \otimes \omega_T^{-m}$ is weakly positive. However, this is the same statement as the weak positivity of $f_* \sO_X (m(K_{X/T} + \Lambda))$, as weak positivity is decided on any big open sets and over the regular locus of $T$, the two sheaves are the same (see \autoref{sec:relative_canonical_bundle_normal_spaces}).

By definition, the above weak positivity means that for a fixed ample divisor $A$ on $T$ there is an integer $b>0$ such that  $\Sym^{[b]}(f_* \sO_X (m (K_{X/T} + \Lambda))) \otimes \sO_T(bA)$ is generically globally generated \cite[Def 4.2]{Ejiri_Weak_positivity_theorem_and_Frobenius_stable_canonical_rings_of_geometric_generic_fibers}. 
  However, then the following natural composition homomorphism (which is non-zero for $m>0$ divisible enough by the relative ampleness along the general fibers)   shows that
 $K_{X/T} + \Lambda  + \frac{1}{m}f^*A$ is $\bQ$-effective:
 \begin{multline}
 \label{eq:semipostivitiy:homomorphism_from_weakly_positive}
 f^* \left(\Sym^{[b]}( f_* \sO_X (m (K_{X/T} + \Lambda))) \otimes \sO_T(bA)\right) 
\explparshift{370pt}{-100pt}{\to}{Obtained by applying $f^* ( \_ \otimes \sO_T(bA))$ to $\Sym^{[b]}( f_* \sO_X (m (K_{X/T} + \Lambda)))\to f_* \sO_X (bm (K_{X/T} + \Lambda))$. The latter is constructed by noting that $f_* \sO_X (bm (K_{X/T} + \Lambda))$ is reflexive (explained in the next paragraph). Hence the multiplication map $\Sym^{b}( f_* \sO_X (m (K_{X/T} + \Lambda))) \to f_* \sO_X (bm (K_{X/T} + \Lambda))$ factors through $\Sym^{[b]}( f_* \sO_X (m (K_{X/T} + \Lambda)))$.}
  f^* \left( f_* \sO_X (bm (K_{X/T} + \Lambda)) \otimes \sO_T(bA)\right)
\\ \expl{\cong}{projection formula, as $A$ is Cartier} f^* f_* \sO_X(bm(K_{X/T} + \Lambda)+ f^* bA)  
\explparshift{200pt}{50pt}{\to}{evaluation homomorphism $f^* f_*(\_) \to (\_)$, or can be viewed also as one of the natural transformations associated to the adjoint functors $f_*$ and $f^*$}   \sO_X(bm(K_{X/T} + \Lambda)+ f^* bA).
 \end{multline}
For the reflexivity of $f_* \sO_X (bm (K_{X/T} + \Lambda))$ note that in general a coherent sheaf $\sF$ on a normal variety $W$ is reflexive if and only if for every open set $U \subseteq W$ and every big open subset $V \subseteq U$, the restriction homomorphism $\sF(U) \to \sF(V)$ is an isomorphism \cite[Prop 1.6]{Hartshorne_Stable_reflexive_sheaves}. This is a property which is preserved under $f_*$ as $f$ is equidimensional. Indeed, if $V \subseteq U$ is a big open subset of $T$ then so is $f^{-1}V \subseteq f^{-1}U$ for $f$ equidimensional.
 
  As $\Gamma_v$ is effective,  $K_{X/T} + \Lambda + \Gamma_v + \frac{1}{m}f^*A$ is also $\bQ$-effective
 As this is true for every $m>0$ divisible enough, $K_{X/T} + \Lambda+ \Gamma_v$ is pseudo-effective. This concludes the main case of the lemma.
 
For {\scshape case $(*)$} note that we only have to prove that on every horizontal curve $C$ on $X$ we have $(K_{X/T} + \Lambda+ \Gamma_v) \cdot C \geq 0$. However, by the $f$-ample assumption and by the generic global generation of $\Sym^{[b]}( f_* \sO_X (m (K_{X/T} + \Lambda))) \otimes \sO_T(bA)$, the homomorphism of \autoref{eq:semipostivitiy:homomorphism_from_weakly_positive} is generically surjective over $C$ for every $q>0$ divisible enough. Hence,  
\begin{equation*}
\left(K_{X/T} + \Lambda + \Gamma_v + \frac{1}{m}f^*A \right) \cdot C
\geq \left(K_{X/T} + \Lambda  + \frac{1}{m}f^*A \right) \cdot C
 \geq 0
 \end{equation*}
  for every $m>0 $ divisible enough. In particular, $(K_{X/T} + \Lambda + \Gamma_v ) \cdot C \geq 0$, which concludes also {\scshape Case $(*)$}. 
\end{proof}

\section{Flatness}
\label{sec:flatness}

The proof of the next theorem closely follows the argument of  \cite[Theorem]{Lu_Tu_Zhang_Zheng_On_semistability_of_Albanese_maps}, with little modifications. Our main contribution is the realization that there exists a semi-positivity statement in positive characteristic (\autoref{lem:semipositivity}) that is needed to implement the argument of \cite[Theorem]{Lu_Tu_Zhang_Zheng_On_semistability_of_Albanese_maps} in positive characteristic.

\begin{theorem}
\label{thm:nef_anti_rel_canonical_flat}
Let $f : (X, \Delta) \to T$  be a surjective fibration  from a normal projective pair to a smooth variety such that  the general fibers   are normal, $-(K_{X/T} + \Delta)$ is a nef $\bQ$-Cartier divisor and  either
\begin{enumerate}
\item \label{itm:nef_anti_rel_canonical_flat:SFR}  $(X_t, \Delta_t)$ is strongly $F$-regular for $t \in T$ general, or
\item \label{itm:nef_anti_rel_canonical_flat:SFP} $(X_t, \Delta_t)$ is sharply $F$-pure  for $t \in T$ general and the the Cartier index of $K_{X/T} + \Delta$ is  prime to $p$ over some non-empty open set of $T$. 
\end{enumerate}
Then $f$ is equidimensional. 
\end{theorem}

\begin{proof}

By base extension we may assume that $k$ is algebraically closed. We do this because at some point we will take general hyperplane sections. 

If $f$ is equidimensional, then there is nothing to prove. Hence, we may assume that $f$ is not equidimensional. 

\emph{We claim that we may also assume that there is a closed point $t_0 \in T$ such that $\codim_X f^{-1}(t_0)  =1$.} Suppose that there is no $t_0$ as in the claim. Choose then $t_0$ with  $\codim_X f^{-1}(t_0)$ minimal. Let  $T'$ be a general hypersurface of $T$ through $t_0$, and let $X'$ be the normalization of the scheme theoretic preimage $f^{-1}\left(T'\right)$. The above claim then follows from our \emph{second claim that $f' : X' \to T'$ satisfies the assumptions of the theorem, with adequately chosen  boundary $\Delta_{X'}$ (made explicit below) and with  $\codim_{X'} \left(f'\right)^{-1}(t_0) = \codim_X f^{-1}(t_0) -1$.} First we note that $Z:=f^{-1}\left(T'\right)$ is irreducible with normal general fibers: let $\Exc(\_)$ denote the locus where the fibers of the given morphism do note have the expected dimension, which is a closed set \cite[Exc. II.3.22.d]{Hartshorne_Algebraic_geometry}, and let $g \colon Z \to T'$ be the induced morphism. 
Then, by Krull's Hauptidealsatz and the genericity of $T'$, $Z$ is equidimensional with $\dim Z = \dim X -1$. Additionally, 
\begin{equation}
\label{eq:nef_anti_rel_canonical_flat:dim_estimate_Exc}
 \dim \left(\Exc(g) \setminus g^{-1}(t_0)  \right) 
 \expl{\leq}{Krull's hauptidealsatz and the genericity of $T'$} 
 \dim \Exc(f) -1< \dim X -1 = \dim Z,
\end{equation}
and
\begin{equation}
\label{eq:nef_anti_rel_canonical_flat:dim_estimate_t_0_inverse}
\dim g^{-1}(t_0) = \dim f^{-1} (t_0) 
\expl{<}{the assumption that $\codim_X f^{-1}(t_0) >1$ }
 \dim X -1 = \dim Z.
\end{equation}
So, \autoref{eq:nef_anti_rel_canonical_flat:dim_estimate_Exc} and \autoref{eq:nef_anti_rel_canonical_flat:dim_estimate_t_0_inverse} yields that no component of $Z$ can be contained in $\Exc(g)$. However, then all components of $Z$ dominate $T'$. On the other hand,  the normality of the general fibers of $X$ implies the same for $Z$.  Combining the previous two sentences, we obtain that $Z$ is irreducible with normal general fibers, and the morphism $X' \to Z$ is isomorphic at the general fibers. In particular, we obtain for the divisorial pullback that $f^* T'=Z$, and  as cohomology and base-change always holds at general points, $f'$ is also a fibration.

In particular, by \autoref{eq:relative_canonical_bundle_adjunction_smooth_base},  there is a canonical bundle formula $K_{X/T} + \Delta|_{X'} \sim_{\bQ} K_{X'/T'} + \Delta_{T'}$, where $\Delta_{T'}$ is simple base-change at general fibers, but it can be more subtle over special fibers. This together with  $X \to Z$ being  an isomorphism at general fibers implies that all assumptions on the general fibers of $f$ still holds for $f'$. It also implies that $K_{X'/T'} + \Delta_{T'}$ is an antinef $\bQ$-Cartier divisor. This concludes our claim.  That is, we may assume that 
\begin{equation}
\label{eq:nef_anti_rel_canonical_flat:codim_assumption}
\codim_X f^{-1}(t_0)=1.
\end{equation}
Let $T' \to T$ be a flatification \cite[(5.2.2)]{Raynaud_Gruson_Criteres_de_platitude_et_de_projectivite_Techniques_de_platification_d_un_module}.  Let  $S$ be the normalization of $T'$. Then, there is a commutative diagram:
\begin{equation*}
\xymatrix{
Y \ar[r]^{\rho} \ar[d]^{h} & X \ar[d]^f \\
S \ar[r]^{\tau} & T
}
\end{equation*}
where $h$ is equidimensional and $Y$ and $S$  are normal. Note that as we took normalization for $S$ and $Y$, $h$ is not necessarily flat. Nevertheless, we know that there is a big regular open set $U \subseteq S$ over which $h$ is flat. As $\rho$ is birational, we may write 
\begin{equation*}
K_Y + \Theta = \rho^*( K_X + \Delta).
\end{equation*}
Similarly, $\tau$ is birational, so we may write:
\begin{equation}
\label{eq:nef_anti_rel_canonical_flat:base}
K_S + \Sigma = \tau^* K_T,
\end{equation}
where $\Sigma \leq 0$. Additionally, as $T$ is smooth, $\Supp (- \Sigma)$ contains all $\tau$-exceptional divisors.  
Applying $h^*$ to \autoref{eq:nef_anti_rel_canonical_flat:base} over $U$ and then taking unique extension ($h^{-1} U$ is big in $Y$ as $h$ is equidimensional), we obtain 
a crepant equation: 
\begin{equation}
\label{eq:nef_anti_rel_canonical_flat:crepant}
K_{Y/S} + \underbrace{\Delta_Y}_{=\Theta - h^* \Sigma} =  \rho^* (K_{X/T} + \Delta).
\end{equation}
It is important to note that $\Delta_Y$ might as well have negative coefficients.
Let $E \subseteq Y$ be the strict transform of a divisorial component of the non-flat locus $V \subseteq X$ of $f$. Because $\codim_{X \setminus \Exc(f)} V \setminus \Exc(f)  \geq 2$, this is the same as requiring $E$ to be a component of $\Exc(f \circ \rho)$. By the equidimensionality of $h$, $E$ maps to a divisor $D$ of $S$, which must be contracted by $\tau$ ($D$ cannot have divisorial image in $T$ because then $E$ would intersect $X \setminus \Exc(f)$). By the smoothness of $T$,  $\coeff_D \Sigma<0$, and as $E$ maps to a divisor of $X$, $\coeff_E \Theta \geq 0$. Hence, 
\begin{equation}
\label{eq:nef_anti_rel_canonical_flat:coeff}
\coeff_E \Delta_Y = \coeff_E \Theta - \coeff_E h^* \Sigma
\expl{>}{$\coeff_D \Sigma<0$ and $\coeff_E \Theta \geq 0$} 
0.
\end{equation}
It follows that, we may write 
\begin{equation}
\label{eq:nef_anti_rel_canonical_flat:Delta_Y}
\Delta_Y = \Gamma + H - G,
\end{equation}
 where  $\Gamma, G$ and $H$ are effective, and  $\Gamma$  agrees with $\Delta_Y$ over general fibers, $G$ is $\rho$-exceptional and $H$ is the part of $\Delta_Y$ that is supported on the strict transforms of the divisorial components of $\Exc(f)$. Furthermore, by \autoref{eq:nef_anti_rel_canonical_flat:codim_assumption} and \autoref{eq:nef_anti_rel_canonical_flat:coeff}, $H \neq 0$. \autoref{lem:semipositivity} then yields that the following divisor is pseudo-effective:
 \begin{equation}
 \label{eq:nef_anti_rel_canonical_flat:initial_pair}
\underbrace{ \underbrace{K_{Y/S} + \Gamma}_{\parbox{80pt}{\scriptsize $(Y, \Gamma)$ satisfies the same singularity assumptions as $(X,\Delta)$ on general fibers}} + 
 \underbrace{(- \rho^* (K_{X/T} + \Delta) }_{\textrm{nef}}}_{\textrm{$f$-nef around the generic fiber}}
 \expllower{60pt}{\sim_{\bQ}}
 {\autoref{eq:nef_anti_rel_canonical_flat:crepant} and \autoref{eq:nef_anti_rel_canonical_flat:Delta_Y}} -H +G 
 \end{equation}
  By intersecting with the preimage of a general complete intersection curve on $X$,  we see that $-H +G$ is not pseudo-effective. This is a contradiction.
\end{proof}

We also provide a general proposition concerning reducedness of fibers avoiding the statements on relative $F$-splitting of the Albanese morphism.

\begin{proposition}
\label{prop:reduced_fibres}
Let $f \colon (X,\Delta) \to T$ be a flat fibration from a projective normal pair to a smooth variety such that the divisor $-K_{X/T} - \Delta$ is $\QQ$-Cartier and nef, and the general fiber $(X_t,\Delta_t)$ is strongly $F$-regular.  Then the following statements hold true:
\begin{enumerate}
\item every geometric fiber of $f$ is reduced,
\item the components of $\Delta$ do not contain any irreducible components of the fibers of $f$.
\end{enumerate}
\end{proposition}

\begin{proof}
We first prove (1).  For the sake of contradiction, we assume that $t \in T$ is such that $f^{-1}(t)$ is non-reduced.  We choose a general curve $C$ passing through $t$ and intersecting the strongly $F$-regular locus of $f$.  There is a finite map $\tau \colon S \to C$ from a smooth projective curve such that if $Z$ is the normalization of $X_S$, then $g \colon Z \to S$ has reduced fibers \cite[Prop 2.1]{de_Jong_Starr_Every_rationally_connected_variety_over_the_function_field_of_a_curve_has_a_rational_point}. Let $\Delta_Z$ be the induced crepant boundary on $Z$ provided by \cite[Prop 2.1]{Codogni_Patakfalvi_Positivity_of_the_CM_line_bundle_for_families_of_K-stable_klt_Fanos} (where the general characteristic zero assumption of \cite{Codogni_Patakfalvi_Positivity_of_the_CM_line_bundle_for_families_of_K-stable_klt_Fanos} is not used).  We recall that $\Delta_Z$ satisfies the equation
\[
K_{Z/S} + \Delta_Z = \rho^*(K_{X/T} + \Delta),
\] 
where $\rho \colon Z \to X$ is the natural morphism. The construction of \cite[Prop 2.1]{Codogni_Patakfalvi_Positivity_of_the_CM_line_bundle_for_families_of_K-stable_klt_Fanos}, is based on \cite[Lem 9.13]{Kovacs_Patakfalvi_Projectivity_of_the_moduli_space_of_stable_log_varieties_and_subadditvity_of_log_Kodaira_dimension},
and in particular $\Delta_Z$ contains the divisor of the map $\omega_{Z/S} \hookrightarrow  \beta^* \omega_{X_S/S}$ induced by the Grothendieck trace $\beta_* \omega_{Z/S} \to \omega_{X_S/S}$, where $\beta \colon Z \to X_S$ is the natural morphism. Using that $S$ is Gorenstein this is the same as the map $\omega_{Z} \hookrightarrow \beta^* \omega_{X_S}$ induced by the Grothendieck trace $\beta_* \omega_{Z} \to \omega_{X_S}$. However, the divisor of $\omega_{Z} \hookrightarrow \beta^* \omega_{X_S}$ is exactly the divisorial part of the conductor for $G_1$ schemes \cite[Lem 2.14]{Patakfalvi_Waldron_Singularities_of_General_Fibers_and_the_LMMP} (and $X_S$ is $G_1$, because $X$ is $G_1$, hence by the smoothness of $T$, $X$ is also relatively Gorenstein over $T$ in codimension $1$, and hence so is $X_S$ over $S$, which finally yields that $X_S$ is $G_1$ by the smoothness of $T$).

So, $\Supp \Delta_Z$ includes the conductor of $Z \to X_S$, and as $X_S$ has at least one non-reduced fiber, this conductor has to contain at least one vertical divisor. Hence, $\Delta_Z \geq E$, where $E$ is a non-zero effective vertical $\bQ$-divisor. Then,  the following equation shows that \autoref{lem:semipositivity} implies that the divisor $-E$ is pseudo-effective.
\begin{equation}
\label{eq:semipositivity_reduced}
\underbrace{\underbrace{K_{Z/S} + (\Delta_Z -E)}_{\parbox{130pt}{\footnotesize $\left(Z_s,\left(\Delta_Z -E \right)_s \right)=\left(Z_s,\left(\Delta_Z  \right)_s \right)$ is strongly $F$-regular for  general closed points $ s \in S$}}  + 
\underbrace{( - (K_{Z/S} + \Delta_Z) )}_{\textrm{nef}}}_{\textrm{$f$-nef generically}}
 \equiv - E 
\end{equation}
This is a contradiction. Hence, the fibers of $f$ are indeed reduced. 

In order to prove (2) we reason as above to reduce the situation to the case where the base schemes $S$ is a curve, and then use the same argument using an equivalent of \eqref{eq:semipositivity_reduced}.
\end{proof}

\section{Numerical flatness}
\label{sec:numerical_flatness}

In the following chapter we prove the main technical result of this paper concerning numerical flatness of relative section rings.

\begin{setting}
\label{notation:triviality_K_trivial}
Let $ f : (X, \Delta) \to T$ be a fibration from a normal pair of dimension $d+1$ to a smooth projective curve 
with normal general fiber and with $-K_{X/T} - \Delta$ a nef $\bQ$-Cartier divisor.
Let us assume additionally that  the general fibers $(X_t, \Delta_t)$ are
 strongly $F$-regular.
(Note that this requirement for general fibers is equivalent to the analogous requirement for the geometric generic fiber by \cite[Theorem B]{Patakfalvi_Schwede_Zhang_F_singularities_in_families})

Let $\wt{L}$ be a very ample divisor on $X$ such that 
\begin{itemize}
\item $R^i f_* \sO_X\left(\wt{L}\right)=0$ for every $i>0$, 
\item $\Sym^m f_*\cO_X\left(\wt{L}\right) \to f_*\cO_X\left(m \wt{L}\right)$ is surjective for every $m>0$.
\end{itemize}
Let $G$ be a general fiber of $f$. 
By \autoref{lem:relative_ample_top_self_intersection}, we have $L^{d+1} = 0$ for $L:= n\wt{L}-mG$, where $n=(d+1)\widetilde{L}_t^d$ and $m= \widetilde{L}^{d+1}$.  This ends the description of our notation.
\end{setting}

\begin{lemma}\label{lem:relative_ample_top_self_intersection}
Let $X$ be a projective variety of dimension $d+1$, let $f \colon X \to T$ be a fibration onto a curve $T$,  let $G$ be the general fiber of the fibration, and let $N$ be an ample Cartier divisor on $X$.   Then  $(nN - mG)^{d+1} = 0$ for $n= (d+1) N_t^d$ and $m= N^{d+1}$.
\end{lemma}

\begin{proof}
We have
\begin{multline}
(nN - mG)^{d+1} 
\\ = n^{d+1}N^{d+1} - (d+1)n^{d}mN^{d} \cdot G + \sum_{2 \leq i \leq d+1} {d+1 \choose i}n^{d+1-i}m^i N^{d+1-i} \cdot G^i
\\ 
\expl{=}{$G^i = 0$, for $i \geq 2$} n^{d}\left(nN^{d+1} - (d+1)mN^{d} \cdot G \right) 
\explshift{-10pt}{=}{$N^{d} \cdot G = N_t^d$}
0
\end{multline}
\end{proof}

The next lemma is shown in \cite[Lem A.2 \& Rem A.3]{Codogni_Patakfalvi_Positivity_of_the_CM_line_bundle_for_families_of_K-stable_klt_Fanos} in characteristic zero. Unfortunately, the argument there is based on a Grothendieck--Riemann--Roch argument, and hence cannot be performed without passing to a resolution of singularities of the total space. Hence, here we have to make an alternative argument using the weaker statement of asymptotic Riemann--Roch \cite[Thm IV.2.15]{Kollar_Rational_curves_on_algebraic_varieties}, which works for singular schemes of any characteristic. An unfortunate side-product is that the argument is a bit more cumbersome than that of \cite[Lemma A.2]{Codogni_Patakfalvi_Positivity_of_the_CM_line_bundle_for_families_of_K-stable_klt_Fanos}. 

 \begin{lemma}
 \label{lem:asymptotic_RR_relative_ample}
 If $f : X \to T$ is a projective morphism from a variety of dimension $d+1$ to a smooth projective curve and $L$ is  an $f$-ample $\bQ$-Cartier divisor on $X$, then for integers $m>0$ divisible enough, 
 \begin{equation*}
 h^0(X,mL)   \geq \frac{m^{d+1}}{(d+1)!} L^{d+1} +  O\left( m^d \right).
 \end{equation*}
 \end{lemma}
 
 \begin{proof}
 The main issue solved during the argument is that \cite[Thm IV.2.15]{Kollar_Rational_curves_on_algebraic_varieties} does not apply directly to $L$ as it is not necessarily nef. Nevertheless, it is $f$-ample, and hence there is $\bQ \ni a >0$ such that $L + a f^* H$ is nef, where $H$ is an ample divisor on $T$. So, we have to apply \cite[Thm IV.2.15]{Kollar_Rational_curves_on_algebraic_varieties} to $L + af^* H$, and then we have to connect on the both ends of the following computation the invariants of $L$ to those of $L+ a f^* H$. First, we do the part of the computation pertaining to $L + a f^*H$. That is, for every integer $m>0$ divisible enough we have:
\begin{multline}
\label{eq:asymptotic_RR_relative_ample:without_ordo_part}
\chi\left(m\left(L+af^*H\right)\right) 
 \explshift{-30pt}{=}{relative Serre vanishing using that $m$ is divisible enogh}
\chi\left(f_* \sO_X\left(m\left(L+af^*H\right)\right)\right)%
 \explshift{30pt}{=}{projection formula}
\chi(\sO_T(maH) \otimes f_* \sO_X(mL)))
\\ \expl{=}{Riemann-Roch for vector bundles on $T$}
\chi( f_* \sO_X(mL))) + ma (\deg H) \left( \rk f_* \sO_X(mL) \right)
\end{multline}
This then yields the statement by the following computation: 
\begin{multline*}
\frac{m^{d+1}}{(d+1)!} L^{d+1} 
\explshift{40pt}{=}{expanding the parentheses in $\left( L + af^*H\right)^{d+1}$}
\frac{m^{d+1}}{(d+1)!}\left( L + af^*H\right)^{d+1}  - \frac{m^{d+1}}{d!} L_t^d a \deg H + O\left( m^d \right)
\\ \explshift{200pt}{=}{\cite[Thm IV.2.15]{Kollar_Rational_curves_on_algebraic_varieties} applied to $m_0(L+ af^*H)$, where $m_0(L+ a f^*H)$ is Cartier}
\chi\left(m\left(L+af^*H\right)\right)  - \frac{m^{d+1}}{d!} L_t^d a \deg H + O\left( m^d \right)
\\ \expl{=}{\autoref{eq:asymptotic_RR_relative_ample:without_ordo_part}}
\chi( f_* \sO_X(mL))) + ma (\deg H) \left( \rk f_* \sO_X(mL) \right)  - \frac{m^{d+1}}{d!} L_t^d a \deg H + O\left( m^d \right)
\\ \explshift{300pt}{=}{\cite[Thm IV.2.15]{Kollar_Rational_curves_on_algebraic_varieties} applied to $m_0L_t$, where $m_0L_t$ is Cartier and $t \in T$ is arbitrary}
\chi( f_* \sO_X(mL)) + O\left( m^d \right)
\leq h^0(T,f_* \sO_X(mL)) + O\left( m^d \right)
=
h^0(X, mL) +  O\left( m^d \right)
\end{multline*}
 
 \end{proof}

\begin{theorem}
\label{thm:nef}
In the situation of \autoref{notation:triviality_K_trivial}, $L$ is nef. 
\end{theorem}

\begin{proof}
First, note that for any $0< \varepsilon \ll 1$, the following is true: there is a $\varepsilon'>0$, depending on $\varepsilon$, such that for all $t\in T$ general (where generality also depends on $\varepsilon$) and for  all $\Gamma' \in |L_t|$, $\left(X_t, \Delta_t + \varepsilon' \Gamma' \right)$ is strongly $F$-regular. Indeed, the above statement is clear for $t$ and $\Gamma'$ fixed. Then, the openness of strong $F$-regularity \cite[Theorem B]{Patakfalvi_Schwede_Zhang_F_singularities_in_families} and Noetherian induction yields the above version with $t$ and $\Gamma'$ not fixed. 

So, fix $\varepsilon>0$. Consider $L + \varepsilon f^* H$. It is enough to prove, by limiting $\varepsilon \to 0$, that $L+ \varepsilon f^* H$ is nef. Set $N:= q(L + \varepsilon f^* H)$, for $q \gg 0$. By \autoref{lem:asymptotic_RR_relative_ample}, 
\begin{equation*}
h^0(X, N)  \geq  q^{d+1} \frac{(L + \varepsilon f^* H)^{d+1}}{(d+1)!} + O(q^d).
\end{equation*}
As
\begin{equation*}
(L + \varepsilon f^* H)^{d+1} 
\expl{=}{$L^{n+1} =0$}
\varepsilon (d+1) (\deg H)  \left(L_t^n \right)    >0,
\end{equation*}
we obtain that $h^0(X,N) \neq 0$ (using $q \gg 0$). So, choose  $\Gamma  \in |N|$. Then, $\left(X_t, \Delta_t+ \varepsilon' \Gamma_t \right)$ strongly $F $-regular for $t \in T$ general. Hence, according to \autoref{lem:semipositivity}, the following divisor is nef  
\begin{equation*}
\underbrace{\underbrace{K_{X/T} +\Delta + \varepsilon' \Gamma}_{\parbox{92pt}{\scriptsize $\left(X_t, \Delta_t+ \varepsilon' \Gamma_t \right)$ is  strongly $F $-regular for $t \in T$ general} } + (\underbrace{ - (K_{X/T} + \Delta)}_{\textrm{nef}})}_{\textrm{$f$-ample}}
 \equiv 
 \varepsilon' q(L + \varepsilon f^* H) 
\end{equation*}
Hence $L + \varepsilon f^* H$  is nef, which concludes our statement.  
\end{proof}

Recall that a vector bundle $\sE$ on a smooth, projective curve $T$ is strongly semi-stable if $\left(F^*\right)^n \sE$ is semi-stable for all integers $n \geq 0$.

\begin{notation}
\label{notation:HN_Frobenius}
Let $\sE$ be a vector bundle on a projective curve $T$. Let 
\begin{equation*}
 0=\sE_{l,0} \subseteq 
 \sE_{l,1} \subseteq 
\dots
 \sE_{l,s_l}  = F^{l,*} \sE.
\end{equation*} 
be the Harder--Narasimhan filtration of $F^{l,*} \sE$. 
Then, it is shown \cite[Claim 2.7.1]{Langer_Semistable_sheaves_in_positive_characteristic} that the sequences $ \frac{\mu_{\max}\left(F^{l,*} \sE \right) }{p^l}:= \frac{\mu\left( \sE_{l,1} \right)}{p^l}$ and $ \frac{\mu_{\min}\left(F^{l,*} \sE \right) }{p^l} := \frac{\mu\left( \factor{\sE_{l,s_l}}{\sE_{l, s_l -1}} \right)}{p^l}$ are stabilizing and for $l \gg 0$, the quotients $\factor{\sE_{l,i}}{\sE_{l,i-1}}$ are strongly semi-stable for all $i$. So one defines $L_{\max}(\sE)$ and $L_{\min}(\sE)$ as the stable values of the above sequences.
\end{notation}

\begin{lemma}
\label{lem:nef_L_min_max}
Let $\sE$ be a vector bundle on a smooth, projective curve $X$. Using \autoref{notation:HN_Frobenius},  $\sE$ is nef if and only if $L_{\min}(\sE) \geq 0$. Similarly, $\sE^*$ is nef if and only if
$L_{\max}(\sE) \leq 0$.
\end{lemma}

\begin{proof}
We use \autoref{notation:HN_Frobenius} throughout the proof. Furthermore, 
we are using the following definition of nefness: $\sE$ is nef, if for every finite morphism $\tau : Y \to X$ from a smooth projective curve, and every quotient line bundle $\tau^* \sE \to \sL$, $\deg \sL \geq 0$. Also, we  show only the first equivalence, as the second equivalence follows from it by applying $\left(F^{l,*} \sE\right)^* \cong F^{l,*} \sE^*$. 

Now, assume that $\sE$ is nef. Then, so is $F^{l,*} \sE$ for every $l \gg 0$, and so are  all the quotient bundles of $F^{l,*} \sE$  (for example using the above definition). However, $F^{l,*} \sE/ \sE_{l,s_l -1}$ is a quotient of $F^{l,*} \sE$ with slope equal to $p^l L _{\min}(\sE)$ for $l \gg 0$. This shows that $L_{\min} \geq 0$.

In the other direction, assume that $L_{\min} \sE \geq 0$, and let $\tau \colon Y \to X$ be a finite morphism from a smooth projective curve, and $\tau^* \sE \to \sL$ a quotient line bundle. For proving that $\deg \sL \geq 0$ we may freely pull-back further this quotient. Hence, we may assume that $\tau$ factors through $F^l$ for $l \gg 0$ as:
\begin{equation*}
\xymatrix{
Y \ar@/^2pc/[rr]^\tau \ar[r]^\rho &  X \ar[r]_{F^l} & X
}
\end{equation*}
Hence, $\sL$ will have a non-zero map from $\rho^* \sE_{l,i} / \rho^* \sE_{l,i-1} \cong \rho^* (\sE_{l,i}/\sE_{l,i-1})$ for some $i$. As $\sE_{l,i}/\sE_{l,i-1}$ is strongly semi-stable with slope at least $p^l L_{\min}(\sE)$, this implies that $\deg \sL \geq 0$ ,because by \cite[Lemma 1.1]{Gieseker_On_a_theorem_of_Bogomolov_on_C} taking pullbacks by separable maps preserves semistability.
\end{proof}

\begin{notation}
\label{notation:Frobenius_pullback}
We recall our notation for the Frobenius pullback 
\[
\xymatrix{
  X^l \ar[rd] \ar@/_1.4pc/[rdd] \ar@/^1.2pc/[rrd]^{F^l_X}& & \\
   & W := X_{T^l} \ar[r]\ar[d]_{g = f_{T^l}}\ar@{}[dr]|-{\square} & X \ar[d]^{f} \\
   & S := T^l \ar[r]_{F^l_S} & T.
}
\]
\end{notation}

\begin{theorem}
\label{thm:dual_nef}
In the situation of \autoref{notation:triviality_K_trivial}, for any integer $a> 0$, $(f_* \sO_X(aL))^*$ is nef.
\end{theorem}

\begin{proof}
Assume $(f_* \sO_X(aL))^*$ is not nef. 
According to \autoref{lem:nef_L_min_max}, and using \autoref{notation:Frobenius_pullback}, this means that for every $l \gg 0$,
\begin{equation*}
0< \mu_{\max}  \left( F^{l,*} f_* \sO_X(aL) \right)  
\expl{=}{$F^{l,*} f_* \sO_X(aL)\cong g_* \sO_W(aL_{T^l}) $ by flat base-change \cite[Prop III.9.6]{Hartshorne_Algebraic_geometry} }
 \mu_{\max}  \left(  g_* \sO_W(aL_{T^l}) \right) 
\end{equation*}  
Furthermore, as $l \gg 0$, we also have that 
\begin{equation*}
\mu_{\max}  \left(  g_* \sO_W(aL_{T^l}) \right) = \mu (\sE),
\end{equation*}
where $\sE$ is a strongly semi-stable subbundle of $g_* \sO_W(aL_{T^l})$. According to \cite[Thm 6.1]{Langer_Semistable_sheaves_in_positive_characteristic}, $\sE^{\otimes r}$ is semi-stable with slope $r\mu (\sE)$. In particular, we may fix $r >0$ such that $\mu\left(\sE^{\otimes r}\right) > 2g\left({T^l}\right)+1= 2 g(T) +1$. Consider then the composition
\begin{equation*}
\sE^{\otimes r} \hookrightarrow \left(  g_* \sO_W(aL_{T^l}) \right)^{\otimes r} \to g_*\sO_W(raL_{T^l}) .
\end{equation*}
This is non-zero (diagonal tensors do not go to zero). So, we obtain that 
\begin{equation*}
\mu_{\max} \left(  g_*\sO_W(raL_{T^l}) \right) > 2g\left(T^l\right) +1 = 2 g(T) +1.
\end{equation*}
Hence for any closed point $t \in {T^l}$:
\begin{equation*}
0 
\expl{\neq}{\cite[Prop 5.7]{Codogni_Patakfalvi_Positivity_of_the_CM_line_bundle_for_families_of_K-stable_klt_Fanos}}
 H^0\left({T^l}, \big( g_*\sO_W(raL_{T^l}) \big)(-t) \right) 
\explpar{160pt}{  =}{$\left( X_{T^l} \right)_t = X \times_T T^l \times_{T^l} t = X \times_T t = X_t$, so $X_t$ is regarded as a divisor on $X_{T^l}$ } 
 H^0(W, raL_{T^l} - X_t).
\end{equation*}
Choose $0\neq\Gamma \in \left| raL_{T^l} - X_t\right| $. Let $Z$ be the normalization of $X_{T^l}$, and let $\Delta_Z$ be the induced crepant boundary given by \cite[Prop 2.1]{Codogni_Patakfalvi_Positivity_of_the_CM_line_bundle_for_families_of_K-stable_klt_Fanos} (again, we recall that the general characteristic zero assumption of \cite{Codogni_Patakfalvi_Positivity_of_the_CM_line_bundle_for_families_of_K-stable_klt_Fanos} is not used). That is, if $\rho: Z \to X$ is the induced morphism, then we have
\begin{equation*}
K_{Z/T^l} + \Delta_Z \sim_{\bQ} \rho^* (K_{X/T} + \Delta).
\end{equation*}
Additionally, by \cite[Prop 2.1.(iii)]{Codogni_Patakfalvi_Positivity_of_the_CM_line_bundle_for_families_of_K-stable_klt_Fanos}, over the locus where $f$ has normal fibers, $X_{T^l}$ is already normal, and $\Delta_Z$ just the usual pullback of $\Delta$. Hence, over this locus the fibers $\left(Z_t, \left( \Delta_Z\right)_t \right)$ are strongly $F$-regular. 
However, then for any $0 < \varepsilon \ll 1$, $\left(Z_t, \left(\Delta_Z\right)_t + \varepsilon \Gamma_t\right)$ is  strongly $F$-regular pair for $t \in T$ general, where the notion of generality depends on $\varepsilon$ \cite[Theorem B]{Patakfalvi_Schwede_Zhang_F_singularities_in_families}.
In particular, according to \autoref{lem:semipositivity}, the following divisor is nef:
\begin{equation*}
\varepsilon \Gamma \sim_{\bQ} 
\underbrace{
\underbrace{K_{Z/T^l} + \Delta_Z +   \varepsilon \Gamma}_{\parbox{118pt}{\scriptsize $\left(Z_t, \left( \Delta_Z \right)_t + \varepsilon \Gamma_t\right)$  is strongly $F$-regular for $t \in T$ general  }} 
+(\underbrace{-(K_{Z/T^l} + \Delta_Z)}_{\textrm{nef}}) 
}_{\textrm{$f$-ample}}
\end{equation*}
 However,  then
\begin{equation*}
\Gamma^{d+1} =\left( raL_{T^l} - Z_t\right)^{d+1} 
\expl{=}{$L^{d+1}=0$} 
- (d+1) ra L_t^d <0,
\end{equation*}
which contradicts nefness of $\Gamma$. Hence, our initial assumption was false, which concludes our proof.  
\end{proof}

We leave the proof of the following lemma to the reader.

\begin{lemma}
\label{lem:nef_via_Frobenius_pullback}
Let $\sE$ be a vector bundle on a smooth projective curve $T$. Then, $\sE$ is nef if and only if for every integer $e>0$, $\left( \left(F^e\right)^* \sE \right) (t) $ is nef for some closed point $t \in T$. 
\end{lemma}

Recall that a vector bundle $\sE$ on a projective scheme $X$ is called \emph{numerically flat}, if both $\sE$ and $\sE^*$ are nef (see \autoref{def:num_flat_def}). 

\begin{theorem}
\label{thm:K_trivial_numerically_flat}

In the situation of \autoref{notation:triviality_K_trivial}, 
$\sF:=f_* \sO_X( mL)$ is numerically flat for every integer $m > 0$. 
\end{theorem}

\begin{proof}
By base extension we may assume that $k$ is algebraically closed. We do this, because at some point we need to use \cite[Prop 3.6]{Patakfalvi_Semi_positivity_in_positive_characteristics}, which is stated over algebraically closed fields. 

By \autoref{thm:dual_nef}, $\sF^*$ is nef. 
Furthermore,  using \autoref{thm:nef}, $L$ is nef. Hence, if the Cartier index of $K_{X/T} + \Delta$  is prime to $p$, we may apply \cite[Prop 3.6]{Patakfalvi_Semi_positivity_in_positive_characteristics} to obtain that $\sF$ is nef, because
\begin{equation*}
K_{X/T}  + \Delta + \underbrace{(-( K_{X/T}+\Delta)) + m L)}_{\textrm{$f$-ample and nef}} \equiv mL .
\end{equation*}
However, if the above index is divisible by $p$, we need a subtle perturbation argument. This is allowed by the fact that $(-( K_{X/T}+\Delta)) + m L)$ is $f$-ample and nef, and hence a perturbation coming from the base makes it ample. The precise argument is as follows.

First, we see that by \autoref{prop:reduced_fibres} all the fibers of $f$ are reduced.  We now approach the nefness of $\sF$. According to \autoref{lem:nef_via_Frobenius_pullback} it is enough to show that for every $e>0$, $\left(\left(F^e\right)^* \sF\right)(t)$ is nef for some closed point $t \in T$. So, fix an integer $e>0$. By our claim the fibers of  $g : V:=X_{T^e} \to T^e$ are reduced. As we also assume that the general fibers are normal, it follows that $V$ is normal. Then, for $\rho : V \to X$ we have
\begin{equation*}
\left(F^e\right)^* \sF  \cong g_* \sO_V(m  \rho^* L).
\end{equation*}
Choose $t \in T^e$. Then $V_t + m \rho^* L $ is ample. According to \cite[Lem 3.15]{Patakfalvi_Semi_positivity_in_positive_characteristics}, we may choose an ample divisor $A$ on $Z$, a sequence $0< \varepsilon_n \to 0$, and $\Gamma \in |K_Z +  A|$,  such that  $(V,\Delta_V + \varepsilon_n (\Gamma+ \Delta_V))$ has index prime-to-$p$ and $\left(V_t,\left(\Delta_V + \varepsilon_n (\Gamma+ \Delta_V)\right)_t \right)$ is strongly $F$-regular for general closed point $t \in T^n$. Then, for $n>0$, \cite[Prop 3.6]{Patakfalvi_Semi_positivity_in_positive_characteristics} implies that the following sheaf is nef
\begin{multline*}
g_* \sO_V(
\underbrace{K_{V/T^e} + \Delta_V + \varepsilon_n (\Gamma + \Delta_V)}_{\parbox{150pt}{\scriptsize $\left(V_t,\left(\Delta_V + \varepsilon_n (\Gamma+ \Delta_V)\right)_t \right)$ is strongly $F$-regular for general closed point $t \in T^n$, and Cartier index is prime to $p$}}  + 
\underbrace{(\underbrace{- (K_{V/T^e} + \Delta_V)}_{\textrm{nef}} + \underbrace{\underbrace{m \rho^* L + V_t}_{\textrm{ample}} - \varepsilon_n (\Gamma + \Delta_V)}_{\textrm{ample for $n \gg0$}}))}_{\textrm{ample for $n \gg 0$}} 
\\ \cong g_* \sO_V( m \rho^*L + V_t ) 
\cong
 \left( \left(F^e\right)^* \sF \right)(t)
\end{multline*}
This concludes our proof. 
\end{proof}

\section{Singularities of the general fibers of sharply $F$-pure fibrations}
\label{sec:singularities_general_fibres}

In this section, we prove a general statement roughly saying that if a variety $X$ that has strongly $F$-regular singularities is fibered over $T$ such that the general fiber $X_{t_{\gen}}$ is normal and sharply $F$-pure, then $X_{t_{\gen}}$ is in fact strongly $F$-regular. This statement is local on $X$. Nevertheless, we would like to use it primarily in the following global situation: $f \colon X \to T$ is the Albanese morphism and $X$ is globally $F$-split.  To apply in this situation the results of \autoref{sec:flatness} and \autoref{sec:numerical_flatness}, we do need to prove that the general fibers are strongly $F$-regular (see the assumptions of the statements in the above sections). This is the main motivation for the results of the present section. Note that even if we assume that $X$ is smooth we do have to deal with the above issues because the results of \cite{Ejiri_When_is_the_Albanese_morphism_an_algebraic_fiber_space_in_positive___characteristic?}
only tell us even in this case that the general fibers are normal and $F$-pure.

We begin by recalling the necessary statements from \cite{Schwede_Tucker_On_the_behavior_of_test-ideals_under_finite_morphisms} where the authors analyze the behaviour of test ideals under finite maps.  First, mimicking the approach from \emph{loc.cit} (see Setting 6.19 therein), we describe the setting in which we work.

\begin{setting}\label{setting:schwede_tucker}
Let $\pi \colon Y \to X$ be a surjective finite map between normal varieties over a perfect field $k$ of characteristic $p>0$.  Let $\cL = \cO_Y(L)$ be a line bundle on $Y$ and let 
\[
\mathfrak{T} \colon \pi_* \cL \to \cO_X
\]
be a non-zero homomorphism of $\cO_X$-modules (we call it the \emph{trace}).  By the Grothendieck duality (on the normal locus of $X$) for finite maps we obtain:
\begin{multline}
\label{eq:schwede_tucker_correspondence}
\cHom(\pi_*\cL,\cO_X) 
\explparshift{205pt}{30pt}{\isom}{$\cHom(\pi_*\cL,\cO_X)$ is reflexive, hence we may use projection formula on the regular locus of $X$ and then we can extend the isomorphism using reflexivity } 
\cHom(\pi_*\cL,\omega_{X/k}) [\otimes] \omega_{X/k}^{*} 
\explshift{90pt}{ \isom }{Grothendieck duality for $\pi$}
\pi_*\cHom(\cL,\omega_{Y/k}) [\otimes] \omega_{X/k}^{*} 
\\ 
\explshift{200pt}{\isom}{again projection formula over the regular locus of $X$ and then reflexive extension} \pi_*\left(\cL^{*} \otimes \omega_{Y/k} [\otimes] \pi^{[*]} \omega_{X/k}^{*}\right) = \pi_*\cO_Y(-L + K_Y - \pi^*K_X)).
\end{multline}
Consequently, the homomorphism $\mathfrak{T} \in \Hom(\pi_*\cL,\cO_X)$ corresponds to an effective divisor 
\[
R_{\mathfrak{T}} \equiv -L + K_Y - \pi^*K_X.
\]
In fact, if we quotient out by multiplication by a unit on the homomorphism side, this correspondence becomes 1-to-1 by \autoref{eq:schwede_tucker_correspondence}.

Choose now an affine open set $U \subseteq X$ and a trivialization $\xi : \sL|_{f^{-1}U} \xrightarrow{\isom} \sO_{f^{-1}U}$. This yields a homomorphism
\begin{equation*}
\xymatrix@C=80pt{
\mathfrak{T}_{\xi} : \pi_* \sO_{f^{-1}U} \ar[r]^{\mathfrak{T}|_U \circ \pi_* \left(\xi^{-1}\right)} & \sO_U
}
\end{equation*}
It follows from the above definition that $R_{\mathfrak{T}_{\xi}}= R_{\mathfrak{T}}|_U$. 
\end{setting}

Next, we cite an immediate consequence of 
\cite[Theorem 6.25]{Schwede_Tucker_On_the_behavior_of_test-ideals_under_finite_morphisms}. This consequence can be obtained  by applying \cite[Theorem 6.25]{Schwede_Tucker_On_the_behavior_of_test-ideals_under_finite_morphisms} to $\mathfrak{T}_{\xi_U}$ where $\{U\}$ is an open cover  of $X$ with trivializations $\xi_U$ of $\sL$ over $f^{-1}U$. More precisely, $\mathfrak{T}_{\xi}$ is a homomorphism $\pi_* \sO_{f^{-1}U} \to \sO_U$, which extends uniquely to a $K(X)$-module homomorphism $\pi_* K(Y) \to K(X)$. We apply \cite[Theorem 6.25]{Schwede_Tucker_On_the_behavior_of_test-ideals_under_finite_morphisms} to the latter over $U$. Doing this over each element $U$ of our chosen open cover yields the following global statement:   

\begin{theorem}[{\cite[Theorem 6.25]{Schwede_Tucker_On_the_behavior_of_test-ideals_under_finite_morphisms}}]
\label{thm:schwede_tucker}
Let $\pi \colon Y \to X$ and $\mathfrak{T} \colon \pi_*\cL \to \cO_X$ be as in \autoref{setting:schwede_tucker}, and set $\Delta_{Y/\mathfrak{T}} = \pi^*\Delta_X - R_{\mathfrak{T}}$.  Then we have the following relation between the respective test ideals:
\[
\mathfrak{T}\left(\pi_*(\tau(Y,\Delta_{Y/\mathfrak{T}}) \cdot \cL)\right) = \tau(X,\Delta_X). 
\] 
In particular, if the pair $(Y,\Delta_{Y/\mathfrak{T}})$ is strongly $F$-regular and $\mathfrak{T}$ is surjective then $(X,\Delta_X)$ is strongly $F$-regular too.
\end{theorem}

\begin{theorem}
\label{thm:general_Fiber_SFR2}
Let $f \colon (X,\Delta) \to T$ be a fibration of normal varieties such that $K_X + \Delta$ has Cartier index prime-to-$p$, $(X,\Delta)$ is strongly $F$-regular and the geometric generic fiber $\left(X_{\overline{\eta}},\Delta_{\overline{\eta}}\right)$ of $f$ is normal and $F$-pure, then $\left(X_{\overline{\eta}},\Delta_{\overline{\eta}}\right)$ is strongly $F$-regular. 
\end{theorem}

\begin{proof}
For the entire proof we use the notation associated to iterated Frobenius pullbacks introduced in \autoref{notation:Frobenius_pullback}, with the only difference being that we use $e$ instead of $l$ for the  index of the iteration. 

First, we observe that by \cite[last 3 lines of Thm A]{Patakfalvi_Schwede_Zhang_F_singularities_in_families}, it suffices to prove that for every divisible enough integer $e>0$  the Frobenius pullback $\left(X_{T^e},\Delta_{T^e}\right)$ is strongly $F$-regular in a neighborhood of  the generic fiber. 
So, from now, $e>0$ is always any divisible enough integer. In particular, $\sL:=\sO_X\left( \left(1-p^e\right)(K_{X/T} + \Delta)\right)$ is a line bundle. We consider now the relative logarithmic Grothendieck trace, one of the definitions (\cite[Def 2.8-Rem 2.11]{Patakfalvi_Schwede_Zhang_F_singularities_in_families}) of which  is that it is the $\sO_{X^e_T}$-linear  homomorphism 
\[
\phi^e_{X/T, \Delta} \colon F^e_{X/T,*} \sL \to 
\cO_{X^e_T},
\] 
such that, using the language of \autoref{setting:schwede_tucker}, we have
\begin{equation}
\label{eq:general_Fiber_SFR:divisor_of_log_trace}
R_{\phi^e_{X/T, \Delta}} = (p^e-1)\Delta. 
\end{equation}
According to \cite[Lemma 2.16 \& 1st-3rd paragraph of Sec 2.3]{Patakfalvi_Schwede_Zhang_F_singularities_in_families}, using that $\overline{k(\eta)}^{1/p^e} =\overline{k(\eta)}$, we may identify the geometric generic fiber of $\phi^e_{X/T,\Delta}$ with the trace for the geometric generic fiber, that is, we have:
\begin{equation}
\label{eq:general_Fiber_SFR:trace_compatbility}
\left(\phi^e_{X/T, \Delta}\right)_{\overline{\eta}}= \phi^e_{X_{\overline{\eta}}, \Delta_{\overline{\eta}}}.
\end{equation} 
Since $\left(X_{\overline{\eta}},\Delta_{\overline{\eta}}\right)$ is sharply $F$-pure, we see that $\phi^e_{X_{\overline{\eta}}}$ is surjective and hence, after potentially shrinking $T$, we may assume that $\phi^e_{X/T,\Delta}$ is surjective too. Also, by possibly shrinking $T$, as $X_{\overline{\eta}}$ is normal, we may assume that $X_{T^e}$ is normal.    Consequently, we may apply \autoref{thm:schwede_tucker} for $X$, $Y$, $\pi$ and $\Delta_X$ of that theorem  set to be $X_{T^e}$, $X$,  $ F^e_{X/T}$ and $\Delta_{T^e}$ of the present proof, and by setting additionally $\cL := \cO_X\left(\left(1 - p^e\right)(K_{X/T} + \Delta)\right)$ and $\mathfrak{T} := \phi^e_{X/T, \Delta}$.  By a direct computation in the context of \autoref{setting:schwede_tucker} we verify that
\begin{equation*}
\Delta_{X/\phi^e_{X/T,\Delta}}  = \left(F^{e}_{X/T}\right)^*\left(\Delta_{T^e} \right) - R_{ \phi^e_{X/T, \Delta}} 
 \expl{=}{\autoref{eq:general_Fiber_SFR:divisor_of_log_trace}}
 \left( F^e \right)^* \Delta  - \left(p^e-1\right) \Delta = p^e \Delta - \left( p^e -1 \right) \Delta = 
  \Delta.
\end{equation*}
This means that the test ideal $\tau\left(X_{T^e},\Delta_{T^e}\right) = \cO_{X_{T^e}}$ and hence $\left(X_{T^e},\Delta_{T^e}\right)$ is strongly $F$-regular by \autoref{thm:schwede_tucker}.  This concludes our proof. 
\end{proof}

\section{The log-Isom scheme}
\label{sec:Isom_scheme}

In this section we define precisely the variants of $\Isom$ schemes that we use in this article. The main issue is defining the corresponding logarithmic $\Isom$ scheme, i.e., how one interprets in a modular way the restriction that in the logarithmic setting the boundary divisor should  be invariant under automorphisms. Modular here means that the construction should be compatible with base-change. This compatibility for closed points and flat covers of the base is essential for our applications. Nevertheless, the present section is admittedly technical, and it is reasonable to skip it for the first read. 

One main issue is to deal with the changing of the coefficients after restrictions. For example one could have coefficients $1/2$ and $1$ in the boundary which become all $1$ after either restricting to a fiber or an inseparable base-change. In particular, it can happen that some components of the boundary were not allowed to be swapped by an automorphism originally, but they are allowed after the base-change. The only way we are able to remedy this issue is to assume that the log-canonical divisors are $\bQ$-Cartier. That is, given a pair $(X, \Delta)$ such that $m(K_X + \Delta)$ is Cartier, we encode $(X, \Delta)$ by the induced homomorphism $\omega_X^{\otimes m} \to \sM:= \sO_X(m(K_X + \Delta))$. This is similar to the approach presented in \cite[Sec.~6]{Kovacs_Patakfalvi_Projectivity_of_the_moduli_space_of_stable_log_varieties_and_subadditvity_of_log_Kodaira_dimension},  the only significant difference being that the polarization there was log canonical, but here it is given by a line bundle $\sL$ independent of the log canonical divisor. We note that there is also another approach based on the work of Koll\'ar \cite{Kollar_Families_of_divisors}.

\begin{setting}
\label{setting:isom_1}
Consider the following situation:
\begin{enumerate}
\item  \label{itm:isom_1:geom_normal} $f^r \colon X^r \to T \ (r=1,2)$ are two flat families of geometrically normal projective varieties of dimension $n$ over a normal, Noetherian base, 
\item  $\sL_r$ are $f^r$-ample line bundles on $X^r$,
\item  $\Delta_r$ are  effective divisors on $X^r$, such that 
\begin{enumerate}
\item no irreducible component of any fiber $X^r_t$ is contained in $\Supp \Delta_r$, and
\item $K_{X^r/T} + \Delta_r$ is $\bQ$-Cartier
\end{enumerate}
\item we fix an integer $m>0$, such that $m(K_{X^r/T} + \Delta_r)$ is Cartier for both $r=1$ and $2$, we set $\sM_r:= \sO_{X^r}( m(K_{X^r/T} + \Delta_r))$, and we fix  induced homomorphisms $\iota_{\Delta_r,m}: \omega_{X^r/T}^{\otimes m} \to \sM_r$ for $r=1,2$ (in general this homomorphism is not uniquely determined by $\Delta_r$,  only up to multiplying by an element of $\Gamma(X, \sO_X)^\times$), and
\item let $U^r \subseteq X^r$ be the open set where $f^r$ is  smooth. By  assumption \autoref{itm:isom_1:geom_normal} $U^r$ is a relatively big open set over $T$. 
\end{enumerate}
\end{setting}

\begin{definition}
\label{def:isom}
In the situation of  \autoref{setting:isom_1}, we consider the following moduli functors $\Sch_T \to \Sets$ from the category of schemes over $T$ to the category of sets. Here, $S$ is any scheme over $T$, and  $p_r : X^r \times_T S \to X^r$ is the projection onto the first factor. Also, we define the functors only on the objects of $\Sch_T$. On arrows the functors go to the usual pullback maps. 
\begin{equation*}
\Isom_T\left( \left(X^1, \sL_1\right),\left(X^2, \sL_2 \right)\right)(S) := \left\{ \   \alpha : X^1 \times_T S \to X^2 \times_T S \  \left|\    \parbox{85pt}{$\alpha$ is an isomorphism over $S$, and $p_{1}^* \sL_1 \cong_S \alpha^* p_2^* \sL_2$} \right. \  \right\},
\end{equation*}
\begin{multline}
\label{eq:isom:log_isom}
\Isom_T\left(\left(X^1, \Delta_1;\sL_1 \right), \left(X^2,\Delta_2;\sL_2 \right)\right)(S) 
\\ := 
\left\{ \  \begin{pmatrix} \alpha : X^1 \times_T S \to X^2 \times_T S, \\[10pt] \xi: \alpha^* \left(\sM_2 \times_T S \right) \to  \sM_1 \times_T S  \end{pmatrix} \ 
 \left|\    \parbox{225pt}{$\alpha$ and $\xi$ are isomorphisms such that
 \begin{enumerate}
 \item $ p_1^* \sL_1\cong_S \alpha^* p_2^* \sL_2$, 
 \item the following diagram commutes:
 \end{enumerate}
 $
 \xymatrix@C=60pt{
  \alpha^* \omega_{U^2 \times_T S/S}^{\otimes m}  \ar[d]_{d \left( \alpha|_{U^1\times_T S} \right)} \ar[r]^-{\alpha^* \left( \left. \iota_{\Delta_2,m}\right|_{U^2 \times_T S} \right)} & \alpha^* \left( \sM_2|_{U^2 \times_T S} \right) 
 \ar[d]^{\xi|_{U^1 \times_T S}} \\
  \omega_{U^1\times_T S/S}^{\otimes m} \ar[r]^{\left.\iota_{\Delta_1,m}\right|_{U^1 \times_T S}} & \sM_1|_{U^1 \times_T S}
 }
 $
 }  \right. \  \right\}.
%
\end{multline}
We also set 
\begin{equation*}
 \Isom_T\left(\left(X^1, \Delta_1;\sL_1 \right)\right):= \Isom_T\left(\left(X^1, \Delta_1;\sL_1 \right), \left(X^1, \Delta_1;\sL_1 \right)\right).
\end{equation*}
\end{definition}

\begin{remark}
\label{rem:isom}
We note the following:
\begin{enumerate}
\item $\sL \cong_S \sL'$ for two line bundles $\sL$ and $\sL'$ on $X^r \times_T S$ means that $\sL' \otimes \sL^{-1} \cong p_S^* \sM$ for some line bundle $\sM$ on $S$.
\item For ease of notation we omitted the restrictions to $U^1 \times_T S$ on $\alpha$ in the commutative diagram of \autoref{eq:isom:log_isom}.
\item As $U^r$ is the smooth locus of $f^r$, $U^r \times_T S$ is the smooth locus of $f^r \times_S T$. In particular, $\alpha(U^1 \times_T S) = U^2 \times_T S$. Hence, we may regard $\alpha$ also as a morphism $U^1 \times_T S \to  U^2 \times_T S$ over $S$. 

\item As we only deal with the relative canonical bundle of smooth morphisms in the commutative diagram of \autoref{eq:isom:log_isom}, we take them to be the particular model $\wedge^{\dim T/S} \Omega_{U^r\times_T S/S}$. In particular, the homomorphism $d \left( \alpha|_{U^1\times_T S} \right)$ is uniquely defined, not only up to an isomorphism. In particular, if $S' \to S$ is a morphism, then we obtain that $d \left( \left( \alpha|_{U^1\times_T S} \right) \times_S S' \right)= d \left( \alpha|_{U^1\times_T S}  \right) \times_S S'$. 

\item \label{itm:isom:rel_big} as $U^r$ is relatively big, the commutativity of the diagram of \autoref{eq:isom:log_isom} is equivalent to the commutativity of 
\begin{equation*}
 \xymatrix@C=60pt{
  \alpha^{[*]} \omega_{X^2 \times_T S/S}^{[m]}  \ar[d] \ar[r] & \alpha^* \left( \sM_2|_{U^2 \times_T S} \right) \ar[d]^\xi \\
  \omega_{U^1\times_T S/S}^{[m]} \ar[r]  & \sM_1|_{U^1 \times_T S}
 },
\end{equation*}
where the unlabeled arrows are the unique extensions of the maps of \autoref{eq:isom:log_isom}. We note that similarly to the absolute setting, in the relative setting reflexive sheaves are also equivalent to relatively $S_2$ flat sheaves (\cite[Prop 3.1 - Cor 3.8]{Hassett_Kovacs_Reflexive_pull_backs} and \cite[Appendix]{Patakfalvi_Schwede_Zhang_F_singularities_in_families}).

\item $\Isom_T\left(\left(X^1, \Delta_1;\sL_1 \right), \left(X^2,\Delta_2;\sL_2 \right)\right)$ is a slight abuse of notation as it does not only depend on $\Delta_r$ but on the actual choice of $\iota_{\Delta_r,m}$. Nevertheless, the different choices of $\iota_{\Delta_r,m}$ differ by a multiplication by a unit, and hence they induce isomorphic $\Isom$ functors via unique isomorphisms (one composes $\xi$ with the above mentioned units). 
\end{enumerate}
\end{remark}

\begin{remark}\label{remark:isom_of_two_is_a_torsor}
In the situation of \autoref{def:isom}, $\Isom_T\left(\left(X^1, \Delta_1;\sL_1 \right)\right)$ is a group scheme over $T$, and both $\Isom_T\left(\left(X^1, \Delta_1;\sL_1 \right)\right)$ and $\Isom_T\left(\left(X^2, \Delta_2;\sL_2 \right)\right)$ acts on $\Isom_T\left(\left(X^1, \Delta_1;\sL_1 \right), \left(X^2, \Delta_2;\sL_2 \right)\right)$.
\end{remark}

\begin{construction}
\label{constr:isom}
Here we describe the construction of the fine moduli space for the functor $\Isom_T\left( \left( X^1, \sL_1 \right), \left( X^2, \sL_2 \right) \right)$ of \autoref{def:isom}. This moduli space exists as a quasi-projective scheme over $T$. We start the construction with the functor
\begin{multline*}
\Isom_T^{\pre}\left( \left( X^1, \sL_1 \right), \left( X^2, \sL_2 \right) \right)(S)
\\
 := 
\left\{ \   \alpha : X^1 \times_T S \to X^2 \times_T S \  \left|\    \parbox{202pt}{$\alpha$ is an isomorphism over $S$, and for every  \\ $m \in \bZ$ and every $s \in S:$\\ $\chi\left( (p_1^*\sL_1 \otimes \alpha^* p_2^* \sL_2)^{\otimes m}_s\right) = \chi\left( \left(p_1^*\sL^{\otimes 2m}_1 \right)_s\right)$} \right. \  \right\}.
\end{multline*}
According to \cite[page 221-20]{Grothendieck_Fondements_de_la_geometrie_algebrique}, $\Isom_T^{\pre}\left( \left( X^1, \sL_1 \right) \left( X^2, \sL_2 \right) \right)$ exists as a quasi-projective scheme over $T$. More precisely, one has to take the $\Hom_S(X,Y)^P$ from \cite[page 221-20]{Grothendieck_Fondements_de_la_geometrie_algebrique} constructed using the line bundles $\sL$ and $\sM$ from the first line of the middle paragraph of \cite[page 221-20]{Grothendieck_Fondements_de_la_geometrie_algebrique} by setting  $S=T$, $X=X^1$, $Y=X^2$ and $P(n) =\chi\left( \left(p_1^*\sL^{\otimes 2m}_1 \right)_s\right)$.  This is the moduli space representing the functor obtained by letting $\alpha$ in the definition of $\Isom_T^{\pre}$ being arbitrary morphism  over $S$ as opposed to an isomorphism. Then, the scheme representing $\Isom_T^{\pre}$ is the subscheme corresponding to isomorphisms. This is an open subscheme by  line 6 of \cite[page 221-20]{Grothendieck_Fondements_de_la_geometrie_algebrique}.

Then, one can define a morphism
\begin{equation*}
\pol: \Isom_T^{\pre}\left( \left( X^1, \sL_1 \right), \left( X^2, \sL_2 \right) \right) \to \Pic(X^{1}/T)
\end{equation*}
 by the following assignment for every $S \in \Sch_T$: 
\begin{multline*}
\pol(S):  \Isom_T^{\pre}\left( \left( X^1, \sL_1 \right), \left( X^2, \sL_2 \right) \right)(S) \ni \alpha 
\\ \mapsto p_1^*\sL_1 \otimes \alpha^* p_2^* \sL_2^{-1} \in \Pic(X^{1}/T)(S) \Big( 
\expl{=}{by definition, \cite[Prop V.2.1, p 232-04]{Grothendieck_Fondements_de_la_geometrie_algebrique}}
 \factor{\Pic(X^{1} \times_T S)}{\Pic(S)} \Big).
\end{multline*}
As both $\Isom_T^{\pre}(X, \sL)$ and $\Pic(X/T)$ are representable \cite[Thm V.3.1]{Grothendieck_Fondements_de_la_geometrie_algebrique} \cite[Thm 6.3]{Altman_Kleiman_Compactifying_the_Picard_scheme} by locally quasi-projective schemes, $\pol$ is in fact a morphism of schemes. Additionally as functors 
\begin{equation}
\label{eq:Isom_pre_Isom}
\Isom_T\left( \left( X^1, \sL_1 \right), \left( X^2, \sL_2 \right) \right) \cong \Isom_T^{\pre}\left( \left( X^1, \sL_1 \right), \left( X^2, \sL_2 \right) \right) 
\explshift{-320pt}{\times_{\Pic(X^{1}/T)}}{using the morphism $\pol$ on the left side, and on the right side  the morphsim $T \to \Pic(X^{1}/T)$ induced by $\sO_{X^{1}}$} 
T .
\end{equation}
Hence, $\Isom_T^{\pre}\left( \left( X^1, \sL_1 \right), \left( X^2, \sL_2 \right) \right)$ is a closed subfunctor of $ \Isom_T^{\pre}(X, \sL)$ and its fine moduli space exists as a scheme and it is given by the same formula as \autoref{eq:Isom_pre_Isom}. 
\end{construction}

\begin{construction}
\label{constr:log_isom}
The fine moduli space for the functor $\Isom_T\left( \left( X^1, \Delta_1; \sL_1 \right), \left( X^2,\Delta_2; \sL_2 \right) \right)$ of \autoref{def:isom} also exists as a quasi-projective scheme over $T$. To see this, we have to go through a longer construction, a big part of the notation of which is shown on the following diagram.
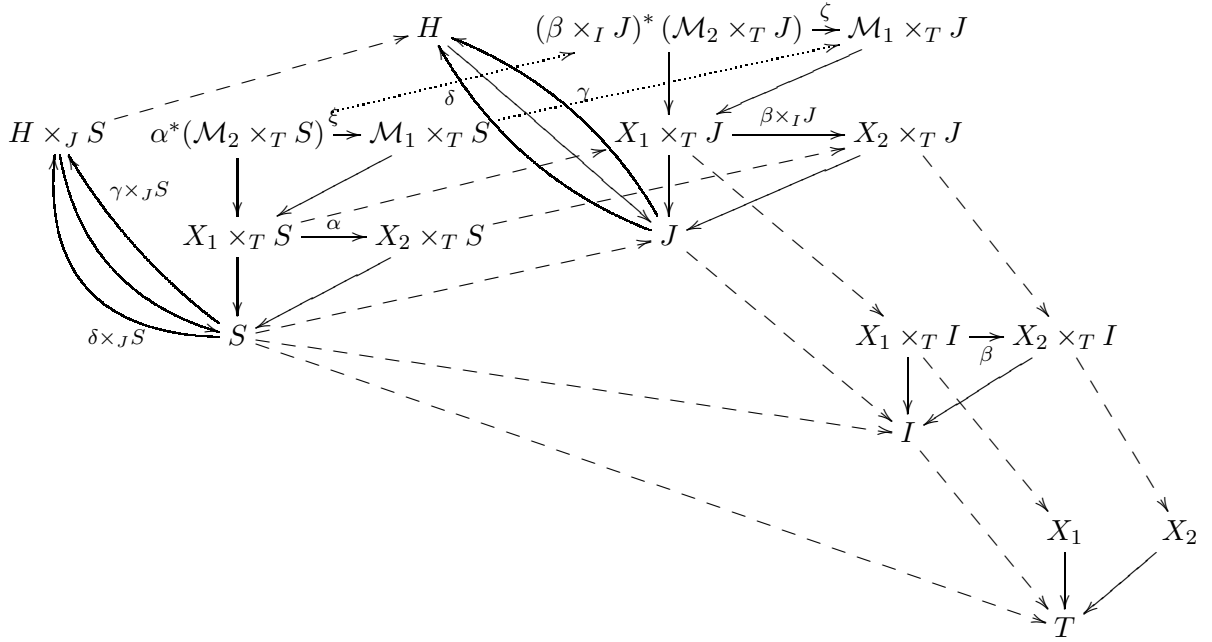
\begin{equation*}
\xymatrix@R=22pt@C=10pt{
& & H \ar[ddr] & \left(\beta \times_I J\right)^* \left( \sM_2 \times_T J\right) \ar[d] \ar[r]^-\zeta &  \sM_1 \times_T J \ar[dl]  \\
%
%
H \times_J S \ar@/_1.5pc/[ddr] \ar@{-->}[urr]
&
\alpha^* ( \sM_2 \times_T S) \ar[d] \ar[r]^{\xi} \ar@{.>}[rru] & \sM_1 \times_T S \ar[ld] \ar@{.>}[rru]  & X^{1} \times_T J \ar[d] \ar@{-->}[rdd]  \ar[r]^{\beta \times_I J} & X^{2} \times_T J \ar[ld] \ar@{-->}[rdd] \\
%
%
& X^{1} \times_T S 
\ar@{-->}[rru] \ar[d] \ar[r]^{\alpha} & 
X^{2} \times_T S \ar[dl]  \ar@{-->}[rru] 
& J \ar@{-->}[rdd] \ar@/^1pc/[uul]^(0.8){\delta} \ar@/_1pc/[uul]_{\gamma} &  & & \\
%
%
& S \ar@/^2.5pc/[uul]^(0.3){\delta \times_J S} \ar@/^0.5pc/[uul]_(0.7){\gamma \times_J S} \ar@{-->}[drrr] \ar@{-->}[dddrrrr] \ar@{-->}[rru] &
 &  & X^{1} \times_T I \ar@{-->}[ddr] \ar[d] \ar[r]_{\beta} 
&
 X^{2} \times_T I \ar@{-->}[ddr] \ar[dl]  \\
%
%
& & & & I \ar@{-->}[ddr]  \\
& & & & & X^{1} \ar[d] & X^{2} \ar[dl] \\
& & & & & T
}
\end{equation*}
First, set $I:=\Isom_T\left( \left( X^1, \sL_1 \right), \left( X^2, \sL_2 \right) \right)$, and let $\beta : X^1 \times_T I  \to X^2 \times_T I$ be the universal isomorphism. In particular, we have $\beta^* \left( \sL_2 \times_T I \right) \cong_I \left( \sL_1 \times_T I \right)$.  
Set $J:=\Isom_I(\beta^*  (\sM_2 \times_T I), (\sM_1\times_T I))$, where $\Isom_I(\_,\_)$ is the open set of $\Hom_I(\_,\_)$ \cite[Def-Lem 33]{Kollar_Hulls_and_Husks} parametrizing isomorphisms only. Let $\zeta : (\beta \times_I J)^*  (\sM_2 \times_T J) \to  \sM_1\times_T J$ be the universal family. Let 
\begin{equation*}
\eta : (\beta \times_I J)^{*} \omega_{X^2 \times_T J /J}^{[m]} 
\expl{\cong}{isomorphisms preserve reflexivity}
 (\beta \times_I J)^{[*]} \omega_{X^2 \times_T J /J}^{[m]} \to \omega_{X^1 \times_T J/J}^{[m]}
\end{equation*}
the induced homomorphism. Set $H:=\Hom_J \left((\beta \times_I J)^{[*]} \omega_{X^2 \times_T J /J}^{[m]}, \sM_1 \times_T J \right)$, and the following sections of $H$ over $J$:
\begin{equation*}
\gamma:=
\explpar{150pt}{\left(\underbrace{\iota_{\Delta_1,m} \times_T J}\right)}{$\in \left\{ \begin{matrix} \Hom_J\left(\omega_{X^{1}/T}^{\otimes m} \times_T J, \sM_1 \times_T J \right) \\ \cong \Hom_J\left(\omega_{X^{1}\times_T J/J}^{[m]} , \sM_1 \times_T J \right) \end{matrix} \right. $}
\circ \eta \in H (J)
\textrm{, and }
\delta:=\zeta \circ (\beta \times_I J)^{[*]}  
\explpar{150pt}{\left( \underbrace{\iota_{\Delta_2,m} \times_T J}\right)}{$\in \left\{\begin{matrix} \Hom_J\left(\omega_{X^{2}/T}^{\otimes m} \times_T J, \sM_2 \times_T J \right) \\  \cong \Hom_J\left(\omega_{X^{2} \times_T J/J}^{[m]} , \sM_2 \times_T J \right) \end{matrix} \right.$}
  \in H(J).
\end{equation*}
Set then the scheme $\Isom_T\left( \left( X^1, \Delta_1; \sL_1 \right), \left( X^2,\Delta_2; \sL_2 \right) \right) \subseteq J$ be the closed subscheme where the sections $\gamma$ and $\delta$ of $H \to J$ agree (or with other words $\gamma^{-1} (\delta)$ or $\delta^{-1} (\gamma)$, where the $\gamma$ and $\delta$ in the parentheses are regarded as a closed subscheme). 

To see that the above construction in fact yields a fine moduli space for the functor $\Isom_T\left( \left( X^1, \Delta_1; \sL_1 \right), \left( X^2,\Delta_2; \sL_2 \right) \right)$, fix
\begin{multline*}
\left(\alpha : X^1 \times_T S \to X^2 \times_T S, \  \xi: \alpha^* \left( \sM_2 \times_T S \right) \to \sM_1 \times_T S \right) 
\\ \in  \Isom_T\left( \left( X^1, \Delta_1; \sL_1 \right), \left( X^2,\Delta_2; \sL_2 \right) \right)(S).
\end{multline*}
Then, there is an induced morphism $S \to I$ over $T$ given by $\alpha$ such that $\beta \times_I S = \alpha$. But then 
\begin{equation*}
\alpha^* (\sM_2 \times_T S ) =  (\beta \times_I S)^*  (\sM_2 \times_T S ) = \left( \beta^* (\sM_2 \times_T I) \right) \times_I S.
\end{equation*}
So, $\xi$ is an isomorphism $\left( \beta^* (\sM_2 \times_T I) \right) \times_I S \to (\sM_1 \times_T I ) \times_I S$. Hence, by the universal property of $J$,  $S \to I$ lifts to $S \to J$ such that $\xi = \zeta \times_J S$. Note that, as $I$ and $J$ are fine moduli spaces themselves as well, 
To conclude showing that $\Isom_T\left( \left( X^1, \Delta_1; \sL_1 \right), \left( X^2,\Delta_2; \sL_2 \right) \right)$ is a fine moduli space for its functor, we have to show that $S \to J$ factors through $\gamma^{-1}(\delta)$ if and only if it satisfies the commutative diagram of \autoref{eq:isom:log_isom}. By \autoref{lem:section_factoring_through}, the condition of $S \to J$ factoring through $\gamma^{-1}(\delta)$ is equivalent to $\gamma_S = \delta_S$. However, as the target of $\gamma_S$ and $\delta_S$ are line bundles and $X^{1} \times_T S \to S$ is a relatively $S_2$ morphism, the $\gamma_S= \delta_S$ holds if and only if $\gamma_S|_{U^{1} \times_T S}= \delta_S|_{U^{1} \times_T S}$ (see also point \autoref{itm:isom:rel_big} of  \autoref{rem:isom}). However, the latter is exactly the commutativity of the diagram of \autoref{eq:isom:log_isom}.
\end{construction}

\begin{lemma}
\label{lem:section_factoring_through}
Let $f : X \to T$ be a morphism with two sections $\sigma_i : T \to X \  ( i=1,2)$. Let $V := \sigma_1^{-1} ( \sigma_2(T)) (= \sigma_2^{-1} ( \sigma_1(T)))$. Then any morphism $S \to T$ factors through $V$ if and only if $\left( \sigma_1 \right)_S = \left( \sigma_2 \right)_S$.
\end{lemma}

\begin{proof}
By considering affine charts, we may assume that all spaces and morphisms are affine.  Let $B$, $A$ and $C$ be the rings corresponding to $X$, $T$ and $C$, respectively, and let $\lambda_i : B \to A$ be the homomorphisms corresponding to $\sigma_i$. Set $I_i:= \ker (\lambda_i)$. In the affine language, $S \to T$ factoring through $V$ means that $\lambda_1(I_2) \subseteq \ker (A \to C)$. On the other hand, $\left( \sigma_1 \right)_S = \left( \sigma_2 \right)_S$ corresponds to $\lambda_1 \otimes_A C(I_2 \otimes_A C)=0$  or equivalently to $\lambda_1 \otimes_A C(I_2 \otimes_A 1)=0$. However, $\lambda_1(I_2) \subseteq \ker (A \to C)$ and $\lambda_1 \otimes_A C(I_2 \otimes_A 1)=0$ are equivalent because they phrase the conditions that the image of $I_2$ via the two routes from the top left corner to the bottom right corner of the following commutative diagram is $0$:
\begin{equation*}
\xymatrix@R=10pt{
B \ar@/_1pc/[d]_{\lambda_1} \ar[r] & B \otimes_A C  \ar@/^1pc/[d]^{\lambda_1 \otimes_A C} \\
A \ar[u] \ar[r] & C \ar[u]
}
\end{equation*}
\end{proof}

\begin{proposition}
\label{prop:base_change_for_isom_scheme}
In the situation of \autoref{setting:isom_1}:
\begin{enumerate}

\item \label{itm:base_change_for_isom_scheme:base_change}
$
\Isom_T\left( \left(X^1, \Delta_1; \sL_1\right),\left(X^2,\Delta_2; \sL_2 \right)\right) $ is compatible with base change to a regular base. That is if  $S \to T$ is a morphism from a regular variety, then
\begin{multline*}
\qquad \Isom_T\left( \left(X^1, \Delta_1; \sL_1\right),\left(X^2,\Delta_2; \sL_2 \right)\right) \times_T S 
\\ \cong \Isom_S\left( \left(X^1_S, \left(\Delta_1\right)_S;\left( \sL_1\right)_S\right),\left(X^2_S,\left(\Delta_2\right)_S; \left(\sL_2\right)_S \right)\right).
\end{multline*}
\item \label{itm:base_change_for_isom_scheme:sections} Over regular bases, the sections of $
\Isom_T\left( \left(X^1, \Delta_1; \sL_1\right),\left(X^2,\Delta_2; \sL_2 \right)\right) $ correspond to the intuitive definition. That is,  for every morphism $S \to T$ from a regular variety we have 
\begin{multline*}
\qquad \Isom_T\left( \left(X^1, \Delta_1; \sL_1\right),\left(X^2,\Delta_2; \sL_2 \right)\right) (S) 
\\= \left\{\alpha \in \Isom_S \left(  \left(X^{1}\right)_S,\left(X^{2}\right)_S \right) \left| \alpha^* \left(\sL_2\right)_S \cong_S \left(\sL_1\right)_S, \ \alpha^* \left(\Delta_2\right)_S = \left(\Delta_1 \right)_S \right. \right\}
\end{multline*}

\end{enumerate}
\end{proposition}

\begin{proof}
All these points follow immediately from the definition. For the last one, use \autoref{rem:isom}.
\end{proof}

\begin{remark}
We would like to emphasize that the naive description stated in point \autoref{itm:base_change_for_isom_scheme:sections} of \autoref{prop:base_change_for_isom_scheme} does not even make much sense over non-reduced bases. There the functor can be regarded directly via the definition in \autoref{def:isom}.
\end{remark}

\section{Isotriviality over curves in special cases}
\label{sec:isotriviality_finite_fields}

In this section we show that in the situation of \autoref{notation:triviality_K_trivial}, the morphism $f$ has isomorphic fibers, in the following two special cases: if the base field $k$ is finite, or if $-(K_{X/T} + \Delta)$ is semi-ample. Unfortunately, because numerical flatness of vector bundles is not an open condition, we are not able to show similar statement without some additional assumptions as the above two. 
In particular, this works if the family is obtained by pulling back over a  general curve in the Albanese variety of a pair $(X,\Delta)$ with $-(K_X + \Delta)$ semi-ample, leading to  one of our main theorems, \autoref{thm:decomposition_theorem}. 

We precede the actual arguments with a few lemmas concerning vector bundles of degree zero.

\begin{lemma}
\label{lem:degree_0_is_trivial}
If $\sE \subseteq \sO_T^{\oplus m}=:\sF$ is a vector bundle of degree zero on a smooth projective curve $T$ over any field $k$, then $\sE = \sO_T \otimes_k V$ for some $V \subseteq H^0(T, \sF)$. 
\end{lemma}

\begin{proof}
We prove it by induction on $m - \rk \sE$. If $m - \rk \sE=0$, then $\sE$ has full rank, and so the only way $\deg \sE=0$ can happen is that $\sE = \sF$. Hence we may assume that $\rk \sE < m$, and that we know the statement for the  higher co-rank case. Then we may find $\sO_T \hookrightarrow \sF $ such that $\sO_T \cap \sE = 0$. In particular, the quotient $\factor{\sF}{\sO_T} =: \sF' \cong \sO_T^{m-1}$ induces an embedding $\sE \hookrightarrow \sF'$. So, we obtain by the induction hypothesis that $\sE \cong \sO_T^{\oplus \rk \sE}$. Then the image of $H^0(T, \sE) \to H^0(T, \sF)$ gives the required $V$. 
\end{proof}

\begin{lemma}
\label{lem:isotrivial_finite_field_special_version}
Let $f \colon X \to T$ be a flat morphism of a projective scheme onto a smooth projective curve over $k$, and let $L$ be an $f$-ample line bundle.  Assume that the following conditions are satisfied
\begin{enumerate}
    \item \label{itm:isotrivial_finite_field_special_version:surjective} the multiplication map $\Sym^m f_*\cO_X(L) \to f_* \cO_X(mL)$ is surjective, for every $m > 0$,
    \item the sheaves $f_*\cO_X(mL)$ are  degree zero vector bundles, for every $m  \geq 1$,
    \item the bundle $f_*\cO_X(L)$ is numerically flat.
\end{enumerate}
Then if $k= \bF_q$ and $t \in T(\bF_q)$, then there is a finite cover $\tau \colon S \to T$ such that $X_S \cong S \times_{\Spec k} X_t$ over $S$.  In fact, the natural map $\Isom_T((X,L),(X_t \times T,L_t \times T)) \to T$ is surjective.  

In particular, the conclusion holds if  \autoref{itm:isotrivial_finite_field_special_version:surjective} is assumed  and if $f_*\cO_X(mL)$ is numerically flat for every $m>0$.
\end{lemma}

\begin{proof}
As numerically flat vector bundles have degree 0, we need to prove only the main statement. By replacing $L$ with some positive multiple we may also assume that the formation of $f_* \sO_X(jL)$ commutes with base-change for every integer $j>0$.

 By assumptions, the bundle $f_* \sO_X(L)$ is numerically flat.  Consequently \autoref{lem:pullback_trivializes} yields a cover $\tau \colon S \to T$ such that $\tau^* f_* \sO_X(L)$ is a trivial vector bundle, that is, $\tau^* f_* \sO_X(L) \cong \sO_S^{\oplus r}$ for some integer $r>0$.  Set $k_S:=H^0(S, \sO_S)$ and $s \in S$ an arbitrary pre-image of $t\in T$. By replacing $S$ with $S \otimes_{k_S} k(s)$ we may assume that $s \in S(k_S)$. 
 
  Let $\rho \colon X_S \to X$ and $g \colon X_S \to S$ be the induced morphisms. Introduce the notation $\xi_m$ for the surjection 
\begin{equation*}
\xi_m : \Sym^m \tau^* f_* \sO_X(L) \cong\Sym^m g_* \sO_{X_S}(\rho^* L) \to g_* \sO_{X_S}(m\rho^*L).
\end{equation*}
Then $\ker \xi_m$ is a degree zero vector bundle in the trivial vector bundle $\Sym^m g_* \sO_{X_S}(\rho^* L)$. Hence, \autoref{lem:degree_0_is_trivial} yields that it is induced by some $V \subseteq H^0(S, \Sym^m g_* \sO_{X_S}(\rho^* L))$. Then, it follows that the multiplication maps $\Sym^m g_* \sO_{X_S}(\rho^* L) \to g_* \sO_{X_S}(m\rho^*L)$ are $\sO_X$-linear homomorphisms between trivial vector bundles, induced by the global section maps 
\[
\Sym^m H^0(S, g_* \sO_{X_S}(\rho^* L)) \to H^0(S, g_* \sO_{X_S}(m\rho^*L)).
\]
This implies that  the relative canonical ring $R_S(X_S,\rho^* L)$ we have 
\begin{multline*}
\hspace{-10pt}
R_S(X_S, \rho^* L) \cong 
\left( \bigoplus_{m \in \bZ} H^0\big(S, g_* \sO_{X_S}(\rho^* mL)\big) \right) \otimes_{k_S} \sO_S 
\explparshift{300pt}{-200pt}{\cong}{as $g_* \sO_{X_S}(\rho^* mL)$ is trivial and $s \in S(k_S)$, the restriction homomorphism $H^0\big(S, g_* \sO_{X_S}(\rho^* mL)\big) \to \big(g_* \sO_{X_S}(\rho^* mL)\big) \otimes k(s) $ is an isomophism} 
\left( \bigoplus_{m \in \bZ} \big( g_* \sO_{X_S}(\rho^* mL)\big) \otimes k(s) \right) \otimes_{k_S} \sO_S 
\\ \explparshift{160pt}{10pt}{\hookrightarrow}{isomorphism outside of degree $0$, as the formation of $f_* \sO_X(jL)$ commutes with base-change for every integer $j>0$}
R\big(X_s, (\rho^*L)_s\big) \otimes_{k_s} \sO_S
\explshift{120pt}{\cong}{flat base-change, as $k \hookrightarrow k_s$ is flat}
(R(X_t, L_t) \otimes_k k_S ) \otimes_{k_S} \sO_S 
\cong 
R(X_t, L_t) \otimes_{\bF_q} \sO_S .
\end{multline*}
As  this embedding is an isomorphism outside of degree zero, it becomes an isomorphism when $\Proj(\_ )$ applied to it. 
  This yields the desired claim.
\end{proof}
\begin{lemma}
\label{lem:isomorphism_over_finite_field_generic_fiber}
Consider a commutative diagram as follows:
\begin{equation*}
\xymatrix{
X \ar[r]^f \ar@/^1.5pc/[rr]^h & T \ar[r]^g & S \ar@/^1.5pc/[l]^{\sigma},
}
\end{equation*}
where
\begin{enumerate}
\item  $S$ is the spectrum of a finitely generated algebra over $\bZ$, 
\item $\sigma$ is a section of $g$,
\item  $f$ and $g$ are flat and projective,
\item $T$ is regular, 
\item \label{itm:isomorphism_over_finite_field_generic_fiber:finite_field} we fix an $f$-ample line bundle $\sL$ such that for every $s \in S$ with $|k(s)|< \infty$, 
\begin{equation*}
%
%
\Isom_{T_s} \left(\left(X_s, \sL_s\right), \left(X_{\sigma(s)} \times_{\Spec k(s)} T_s, \sL_{\sigma(s)} \times_{\Spec k(s)} T_s\right)\right)
\end{equation*}
is surjective and flat on $T_s$. 
\end{enumerate}
Then whenever $\eta \to S$ is a spectrum of a field mapping to the generic point of $S$, the homomorphism  $\Isom_{T_\eta} \left(\left(X_\eta, \sL_\eta\right), \left(X_{\sigma(\eta)} \times_{\Spec k(\eta)} T_\eta, \sL_{\sigma(\eta)} \times_{\Spec k(\eta)} T_\eta\right)\right) \to T_\eta
$ is surjective.

If we assume additionally that 
\begin{enumerate}[resume]
\item all geometric fibers of $f$ are normal,
\item there is an effective $\bQ$-divisor $\Delta$ on $X$ not containing any fiber in its support such that $K_{X/T} + \Delta$ is $\bQ$-Cartier,
\end{enumerate}
and instead of \label{itm:isomorphism_over_finite_field_generic_fiber:finite_field}
we  assume that
\begin{itemize}
\item[(5)'] 
for every $s \in S$ with $|k(s)|< \infty$, 
\begin{equation*}
%
%
\Isom_{T_s} \left(\left(X_{s}, \Delta_{s}; \sL_{s}\right), \left(X_{\sigma(s)} \times_{\Spec k(s)} T_s, \psi_s^* \Delta ;\sL_{\sigma(s)} \times_{\Spec k(s)} T_s\right)\right)
\end{equation*}
is surjective and flat on $T_s$, where $\psi : X \times_{T, \sigma \circ g}  T \to X$ is the projection morphism onto the first factor.  
\end{itemize}
Then, for $\eta \to S$ as above, we obtain that 
\begin{equation*}
\Isom_{T_\eta} \left(\left(X_{\eta}, \Delta_{\eta}; \sL_{\eta}\right), \left(X_{\sigma(\eta)} \times_{\Spec k(\eta)} T_s, \psi_\eta^* \Delta  ;\sL_{\sigma(\eta)} \times_{\Spec k(\eta)} T_{\eta}\right)\right) \to T_{\eta}
\end{equation*}
 is surjective and flat.

\end{lemma}

\begin{proof}
Set in the respective cases
\begin{equation*}
V:=\Isom_{T} \left(\left(X, \sL \right), \left( X \times_{T, \sigma \circ g}  T,  \sL  \times_{T, \sigma \circ g}  T\right)\right) .
\end{equation*}
or
\begin{equation*}
V:=\Isom_{T} \left(\left(X, \Delta ;  \sL \right), \left( X \times_{T, \sigma \circ g}  T, \psi^* \Delta  ;  \sL  \times_{T, \sigma \circ g}  T\right)\right).
\end{equation*}
Then, for all $s \in S$, we have the respective canonical isomorphisms
\begin{equation*}
V_s \cong 
\Isom_{T_s} \left(\left(X_s, \sL_s\right), \left(X_{\sigma(s)} \times_{\Spec k(s)} T_s, \sL_{\sigma(s)} \times_{\Spec k(s)} T_s\right)\right),
\end{equation*}
or
\begin{equation*}
V_s \cong 
\Isom_{T_s} \left(\left(X_s,\Delta_s; \sL_s\right), \left(X_{\sigma(s)} \times_{\Spec k(s)} T_s, \psi_s^* \Delta, \sL_{\sigma(s)} \times_{\Spec k(s)} T_s\right)\right).
\end{equation*}
Hence, our assumption tells us that $V_s \to T_s$ is surjective, whenever $k(s)$ is finite, and we have to deduce that $V_\eta \to T_\eta$ is surjective.  As $V \to T$ is finite type, the image $W$ of $V\to T$ is constructible.  So, we know that $W_s=T_s$ whenever $k(s)$ is a finite field, that is, for all closed points $s \in S$, and for our surjectivity statements we have to show that then $W_\eta = T_\eta$ as well. Assume the contrary. Then $W^c:= T \setminus W$ is a constructible set, which avoids all fibers of $g$ over closed points, but it intersects the generic fiber non-trivially. Hence, the constructible set $g(W^c)$ avoids all closed points of $S$, but contains the generic point. This is a contradiction, concluding our surjectivity statements.

For the flatness statements we use that the locus $\tW:=\{t \in T| V \to T  \textrm{ is flat at } t \in T\}$ is constructible (combination of openness of the flat locus in $V$ \cite[Th\'eor\`eme 11.3.1]{Grothendieck_Elements_de_geometrie_algebrique_IV_III} and Chevalley's theorem about the image of constructible set is constructible). Just as above for $W$ we know that $\tW_s = T_s$ whenever $|k(s)| < \infty$, and then just as above for $W$ we obtain that $\tW_\eta= T_\eta$. 
\end{proof}

\begin{lemma}
\label{lem:from_total_space_trivial_to_pair_trivial}
Let $Y$ and $S$ be normal projective varieties, and set $k_S:= H^0(S, \sO_S)$. 

If $\Gamma$ is an effective $\bQ$-divisor on $Y \times_{\Spec k} S$ such that $K_{Y \times S/S} + \Gamma$ is an anti-nef $\bQ$-Cartier divisor, then $(Y \times_{\Spec k} S, \Gamma) \cong (Y \otimes_k k_S, \Sigma) \times_{\Spec k_S} S$ for some effective $\bQ$-divisor $\Sigma$ on $Y \otimes_k k_S$.
\end{lemma}

\begin{proof}
In this proof, all products are taken over $\Spec k$, which we do not notate in subscript under the $\times$ signs for better readability. Let $\gamma \colon Y \times S \to Z$ be the Stein factorization  of the projection $\alpha \colon Y\times S \to Y$. 
Note that $Z=Y \otimes_k k_S$, because if we apply flat base-change  to the morphism $S \to \Spec k$ via the base-change $Y \to \Spec k$, then we obtain that
$h_* \sO_{Y \times S} \cong H^0(S, \sO_S) \otimes_k \sO_Y$.

As divisors are determined in codimension $1$, we may replace $Y$ by its regular locus, and hence we may assume that $Y$ is regular. This way we do lose the projectivity of $Y$, but this actually does not matter as we only consider $h$-nefness in the present proof.  On the other hand, since $Y$ is regular, we gain that $K_Y$ is Cartier. This is used essentially in the following computation, which implies that every $h$-vertical curve that intersects $\Supp \Gamma$ is necessarily contained in $\Supp \Gamma$: 
\begin{equation*}
 \Gamma =  (K_{Y \times S/S} + \Gamma) - K_{Y \times S/S} 
= 
\underbrace{\underbrace{(K_{Y \times S/S} + \Gamma)}_{\textrm{$h$-anti-nef}} - \underbrace{h^* K_Y}_{\textrm{$h$-numerically trivial}}}_{\textrm{$\bQ$-Cartier and $h$-anti-nef}}.
\end{equation*}
This implies that every component of $\Gamma$ is $h$-vertical. That is, $\Gamma = \gamma^* \Sigma$  for an adequate divisor $\Sigma$ on $Z= Y \otimes_k k_S$.  This concludes our proof.
\end{proof}

\begin{remark}
In the following proof we do not have to chack flatness of $X \to T$ when applying the $\Isom$ argument, as $X$ is a variety mapping surjectively onto a smooth curve $T$ and so $X$ is automatically flat over $T$. The situation is the same in the proof of \autoref{cor:semi_ample}. On the other hand, in other proofs where nefness of the anti log-canonical divisor is assumed but the base is of arbitrary dimension, e.g., in the proof of \autoref{thm:decomposition_theorem_alternative}, we use \autoref{thm:nef_anti_rel_canonical_flat} to obtain the necessary flatness. 
\end{remark}

\begin{corollary}
\label{cor:finite_field}
In the situation of \autoref{notation:triviality_K_trivial}, if  $k= \bF_q$ and $t {}\in T(\bF_q)$, 
then there is a finite cover $\tau \colon S \to T$ by a smooth projective curve such that $(X_S,\rho^* \Delta; \rho^* \sL) \cong_S S \times_k (X_t,\Delta_t, \sL_t)$, where  $\rho: X_S \to X$ is the projection morphism induced by \autoref{eq:finite_field:splitting}. (Here isomorphism over $S$ for the line bundle means that we allow twist by pull-backs of line bundles on $S$, in accordance with \autoref{def:isom}.)

Addendum: in particular, $\Isom_T((X,\Delta,L),(X_t \times T,\Delta_t \times T, L_t \times T)) \to T$ is  surjective and flat.
\end{corollary}

\begin{proof}
It
follows directly from \autoref{lem:isotrivial_finite_field_special_version} using \autoref{thm:K_trivial_numerically_flat} that there is $\tau \colon S \to T$ such that 
\begin{equation}
\label{eq:finite_field:splitting}
(X_S, \rho^*\sL) \cong_S S \times_k (X_t, \sL_t).
\end{equation}
Additionally, by passing to a higher $S$, we may assume that there is an $s \in S(k_S)$ such that $\tau(s)=t$, where $k_S = H^0(S, \sO_S)$. 

We claim that for such a $\tau$ we have $(X_S, \rho^* \Delta; \rho^* \sL) \cong S \times_k (X_t,\Delta_t; \sL_t)$. To show this claim, note first that by \autoref{eq:finite_field:splitting}, we obtain that the fibers of $X \to T$ are all normal. Hence,  $K_{X_S/S} + \rho^* \Delta = \rho^* (K_{X/T} + \Delta)$ by \cite[Prop 2.1]{Codogni_Patakfalvi_Positivity_of_the_CM_line_bundle_for_families_of_K-stable_klt_Fanos}. In particular, $K_{X_S/S} + \rho^* \Delta$ is anti-nef.  Then, \autoref{lem:from_total_space_trivial_to_pair_trivial} implies that $\rho^* \Delta=  \gamma^* \Sigma$ for the composition morphism 
\begin{equation*}
\gamma : X_S 
\explparshift{300pt}{150pt}{\cong}{as the isomorphism \autoref{eq:finite_field:splitting} is over $S$, this isomorphism is automatically over $k_S= H^0(S, \sO_S)$ $\implies$ $s \in S$ yields a section of $\gamma$}
S \times_{\Spec k} X_t \to X_t \otimes_k k_S
\end{equation*} 
and for a divisor $\Sigma$ on $ X_t \otimes_k k_S$. Furthermore, we can determine $\Sigma$, by restricting over $s$:
\begin{equation*}
\Sigma 
\expl{=}{$s \in S$ yields a section of $\gamma$}
 (\gamma^* \Sigma)_s = (\rho^* \Delta)_s= \Delta_t \otimes_k k_S
\end{equation*}
Therefore, $\rho^* \Delta= \Delta_t \times_k S$. 

For the  addendum, use \autoref{prop:base_change_for_isom_scheme} and in the case of flatness combine it with a corollary of the existence of a flattening stratification stating that if a morphism is flat after finite base-change, it is flat already without applying the base-change \cite[\href{https://stacks.math.columbia.edu/tag/0533}{Tag 0533}]{stacks-project}.
\end{proof}

\begin{lemma}
\label{lemma:semi_ample_lemma}
In the situation of \autoref{notation:triviality_K_trivial}, if $-(K_{X/T} + \Delta)$ is semi-ample (e.g., $K_{X/T} + \Delta \sim_{\bQ} 0$) and $t \in T(k)$ is a rational point, then natural morphism $\Isom_T((X,\Delta,L),(X_t \times T,\Delta_t \times T, L_t \times T)) \to T$ is  surjective and flat.
\end{lemma}

\begin{proof}
Choose a model $\xymatrix{ (X_S, \Delta_S) \ar[r]_{f_S} &  T_S \ar[r]_{g_S} & S \ar@/_1pc/[l]_{\sigma_S}  } $  of $\xymatrix{(X, \Delta) \ar[r]_f & T \ar[r]_g & \Spec k \ar@/_1pc/[l]_{\sigma_t} }$ over the spectrum $S$ of a regular finitely generated $\bF_p$-algebra, such that for every closed point $s \in S$, $f_s  : (X_s, \Delta_s) \to T_s$ satisfies the assumptions of \autoref{notation:triviality_K_trivial}.  This is doable by the semi-ample assumption, and by the openness of strongly $F$-regular singularities \cite[Thm B]{Patakfalvi_Schwede_Zhang_F_singularities_in_families}. Moreover, take a spreading out $\cL$ of the line bundle $L$ such that the bundles $\cL_s$ satisfy $\cL_s^{d+1} = 0$.    Then \autoref{cor:finite_field} and \autoref{prop:base_change_for_isom_scheme} tell us that for every closed point $s \in S$, the structure morphism $\Isom_{T_s} \left(\left(X_s,\Delta_s; \sL_s\right), \left(X_{\sigma(s)} \times_{\Spec k(s)} T_s, \psi_s^* \Delta, \psi_s^* \sL_{\sigma(s)} \right)\right) \to T_s$ is surjective, where $\phi$ is the projection $X \times_{T, \sigma_S \circ g_S} T \to X$.  We emphasize that the assumption on normality of the fibers is satisfied because we may first prove that the family is trivial in the boundary free setting.  Additionally, \autoref{cor:finite_field} and \autoref{prop:base_change_for_isom_scheme}, also tells us that the same structure morphism is flat after a finite base, change. Hence flattening stratification, or the more elementary algebra lemma \cite[\href{https://stacks.math.columbia.edu/tag/0533}{Tag 0533}]{stacks-project} tells us  that it is also flat without applying the base-change. Then \autoref{lem:isomorphism_over_finite_field_generic_fiber} concludes our proof. 
\end{proof}

\begin{corollary}
\label{cor:semi_ample}Let $f \colon (X, \Delta) \to T$ be a fibration from a normal pair of dimension $d+1$ to a smooth projective curve with normal general fiber and with $-K_{X/T} - \Delta$ a semi-ample $\bQ$-Cartier divisor.  Assume additionally that  the general fibers of $f$ are strongly $F$-regular. Then $\left(X_{\ot}, \Delta_{\ot}\right) \cong \left(X_{\ot'}, \Delta_{\ot'}\right)$ for any geometric points $\ot, \ot' \in T$ with the same residue fields.  
\end{corollary}

\begin{proof}
According to \autoref{lem:relative_ample_top_self_intersection}, we may choose a line bundle $L$ on $X$ such that $L^{n+1} = 0$, and hence the morphism $f \colon (X,\Delta) \to T$ together with a line bundle $L$ satisfies the conditions of \autoref{notation:triviality_K_trivial}.  Consequently, we may apply \autoref{lemma:semi_ample_lemma} to conclude. 
\end{proof}


\section{Isotriviality of the Albanese morphism}
\label{sec:isotriviality_of_the_Albanese}

In this section we combine the result of the previous sections into the main theorem of the present paper.  

\begin{lemma}
\label{lem:Bertini}
If $k$ is algebraically closed, $W$ is a  smooth variety over $k$,   $v \neq w \in W(k)$, and $L$ is a very ample divisor on $W$, then for every $n\gg 0$ a general element of 
\begin{equation*}
V_n=\big\{\ H \in |nL| \ \big| \ v, w \in \Supp H \ \big\}\end{equation*}
 is smooth. 
\end{lemma}

\begin{proof}
Consider the linear system $V_n$ for any $n \gg 0$. By \cite[\href{https://stacks.math.columbia.edu/tag/0FD5}{Tag 0FD5}]{stacks-project},  $V_n$ globally generates and it defines an immersion over $U:=W \setminus \{v,w\}$. Hence, by \cite[\href{https://stacks.math.columbia.edu/tag/0FD6}{Tag 0FD6}]{stacks-project}, for a general $H \in V_n$, $H|_U$ is smooth. Additionally, by jet separation and by the big enough choice of $n$, $V_n$ separates jets at both $u$ and $v$ and hence the general $H \in V_n$ is also smooth at $u$ and $v$. 
\end{proof}

\begin{theorem}[Main decomposition theorem -- general version]
\label{thm:decomposition_theorem_alternative}
Let $(X,\Delta)$ be a Cohen--Macaulay pair $(X,\Delta)$ such that  either
\begin{enumerate}[label=(\alph*)]
\item \label{situation:decomposition_theorem_alternative:semi_ample} $-K_X-\Delta$ is a semi-ample $\bQ$-Cartier divisor, or
\item \label{situation:decomposition_theorem_alternative:finite_field} $k \subseteq \obF_p$ and $-K_X - \Delta$ is a nef $\bQ$-Cartier divisor.
\end{enumerate}  
Assume that $A$ admits a $k$-point $0 \in A(k)$ and the Albanese morphism $f \colon X \to A$ is a surjective fibration and that the general fiber $(X_t,\Delta_t)$ is strongly $F$-regular.  Then,  
\begin{equation}\label{eq:isom_main_thm}
(X,\Delta) \times_A I \cong_I  (X_0,\Delta_0) \times_k I, 
\end{equation}
where 
\begin{enumerate}
\item \label{itm:decomposition_theorem_alternative:I_def}
\begin{equation*}
I := \Isom_A((X,\Delta;L_X),(X_0 \times_k A, \Delta_0 \times_k A; L_0 \times_k A)),
\end{equation*}
for some relatively ample line bundle $L_X$ on $X$,  $X_0 := f^{-1}(\{0\})$, $\Delta_0:= \Delta|_{X_0}$ and $L_0 := L_{X}|_{X_0}$.
\item \label{itm:decomposition_theorem_alternative:surjective} $I \to A$ is surjective and flat. 
\end{enumerate}
In particular,  all fibers $(X_t, \Delta_t)$ of $f$ over $t \in A\big(\ok\big)$ are isomorphic.  Moreover, the map $I \to A$ has a structure of a torsor under $G = \Isom((X_0, \Delta_0, L_0))$ in the flat topology.

\end{theorem}

\begin{proof}
Let us start the proof by reducing the statement to point \autoref{itm:decomposition_theorem_alternative:surjective}. 
First, note that $\pr_2 \colon I \times_A I \to I$ has the natural diagonal section. Hence, points \autoref{itm:base_change_for_isom_scheme:base_change} and \autoref{itm:base_change_for_isom_scheme:sections} of \autoref{prop:base_change_for_isom_scheme} yield automatically point \autoref{itm:decomposition_theorem_alternative:I_def} of the present theorem as soon as point \autoref{itm:decomposition_theorem_alternative:surjective} is shown.  Second, the as soon as we know point \autoref{prop:base_change_for_isom_scheme}, for an $t \in A(k)$, the fiber $(X_t, \Delta_t)$ is also a fiber of $(X, \Delta) \times_A I$. However, all the fibers of the latter over $k$ points are isomorphic to

Second, note that as the base field $k$ is perfect, by base-changing to the algebraic closure, we may assume that in fact $k$ is algebraically closed. Here we are using \autoref{prop:base_change_for_isom_scheme}. 

To sum up, we assume from now that $k$ is algebraically closed and we only have to prove that that $I \to A$ is surjective and flat.

 By applying \autoref{thm:nef_anti_rel_canonical_flat} we obtain that $f$ is equidimensional and hence flat because $X$ is Cohen--Macaulay. 
Let us  fix then the following notation and the basic setting.  We assume that the dimension of $X$ is equal to $n$, and the relative dimension of $f$ is equal to $d$.  We take a very ample divisor $\wt{L}$ on $X$ satisfying the conditions:
\begin{align*}
R^i f_*\cO_X\left(j \wt{L}\right) = 0, \text{ for $i>0$ and $j>0$,} \\
\Sym^m f_*\cO_X\left(\wt{L}\right) \to f_*\cO_X\left(m \wt{L}\right) \text{ is surjective for every $m>0$.}
\end{align*}
This ensures first that the sheaves $f_*\cO_X\left(j \wt{L}\right)$ are vector bundles the formation of which commutes with every base change. Second, it ensures that the relative section ring 
\[
\bigoplus_{m \geq 0} f_*\cO_X\left(m \wt{L}\right),
\] 
and hence the whole morphism $f$ as well are determined by $f_*\cO_X\left(\wt{L}\right)$ and the natural relations given by the kernel of the multiplication map $\Sym^d f_*\cO_X\left(\wt{L}\right) \to f_*\cO_X\left(d \wt{L}\right)$.

Furthermore we choose a very ample divisor $D$ on $A$, such that for any two closed points $x, y \in A$, and for general $D_1, \dots, D_{n-d-1} \in |D|$ through $x$ and $y$, the curve $C_{x,y}:=\bigcap_{i=1}^{n-d-1} D_i$ is  smooth (by \autoref{lem:Bertini}) and irreducible. 

Set $H = f^*D$.  Along the lines of \autoref{notation:triviality_K_trivial} and \autoref{lem:relative_ample_top_self_intersection}, we observe that $H^{n-d+1} = \left(f^*D\right)^{n-d+1} = 0$, and consequently we may find positive integers $m$ and $n$ such that the divisor $L_X = m\tilde L - nH$ is $f$-ample and satisfies the condition 
\begin{equation}
L_X^{d+1} \cdot H^{n-d-1} = 0 \label{eq:condition_intersection_relative}
\end{equation}  
We remark that $L_X$ is clearly $f$-ample since $H$ is a pullback of a divisor from the base $A$.  By the projection formula it also satisfies the conditions explained above:
\begin{align}
R^i f_*\cO_X\left(j L_X\right) = 0, \text{ for $i>0$ and $j>0$,} \\
\Sym^m f_*\cO_X\left(L_X\right) \to f_*\cO_X\left(m L_X\right) \text{ is surjective for every $d > 0$.}
\end{align}
We prove the surjectivity and the flatness of $I \to A$ using the results of \autoref{sec:numerical_flatness} by restricting to the curves $C_{x,0} \sim D^{n-d-1}$ introduced above, where $x$ is an arbitrary  closed point.  For the ease of notation, fix this $x \in X$, set $C := C_{x,0} \subset A$, and consider the pullback pair $\left(Y,\Delta_Y\right) = \left(X \times_A C,\Delta_Y\right)$ and the Cartier divisor $L = L_{X}|_Y$. To explain what $\Delta_Y$ is, note first that as $X$ is Cohen-Macaulay and $C$ is a smooth curve in $A$, $Y$ is also Cohen-Macaulay. Furthermore, as $C$ goes through a general point of $A$, the general fibers of $Y \to C$ are normal. Taking into account that the special fibers are reduced, we obtain that $Y$ is normal (e.g., \cite[Lemma 6.2]{Codogni_Patakfalvi_Positivity_of_the_CM_line_bundle_for_families_of_K-stable_klt_Fanos}), as well as by \cite[Proposition 2.1]{Codogni_Patakfalvi_Positivity_of_the_CM_line_bundle_for_families_of_K-stable_klt_Fanos} there is a $\Delta_Y$ on $Y$ such that 
\begin{equation}
\label{eq:decomposition_theorem_alternative:adjunction}
K_{X/A} + \Delta|_Y \sim_{\bQ} K_{Y/C} + \Delta_Y.
\end{equation}
 By the base change properties of the log-isomorphism scheme presented in \autoref{sec:Isom_scheme}, and by flattening decomposition \cite[Lecture 8]{Mumford_Lectures_on_curves_on_an_algebraic_surface} it suffices prove that the morphism
\[
\Isom_C((Y,\Delta_Y;L),(X_0 \times C, \Delta_0 \times C;L_0 \times_A C)) \explparshift{160pt}{-50pt}{\isom}{$\Isom$ scheme behaves well with respect to base change (see \autoref{prop:base_change_for_isom_scheme})} I \times_A C \to C
\] 
is surjective. 

To this aim, we want to apply the results of \autoref{sec:isotriviality_finite_fields}.  For this purpose, we verify that the morphism $g \colon (Y,\Delta_Y) \to C$ and the Cartier divisor $L$ satisfy all the requirements described in \autoref{notation:triviality_K_trivial}. 
First, the normality of $Y$ was shown above. Second, if $t \in C$ is general, by \cite[Lemma 6.2]{Codogni_Patakfalvi_Positivity_of_the_CM_line_bundle_for_families_of_K-stable_klt_Fanos}, $\left(Y_t, \left( \Delta_Y\right)_t \right)= (X_t, \Delta_t)$.   In particular, as we have seen above that the latter is strongly $F$-regular, so is the former. 
Third, the condition $L^{d+1} = 0$ is shown by the following computation:
\[
L^{d+1} = L_X^{d+1} \cdot Y 
\explshift{-50pt}{=}{$C \in D^{n-d-1}$ and $Y \in H^{n-d-1}$} 
L_X^{d+1} \cdot H^{n-d-1}
\expl{=}{\autoref{eq:condition_intersection_relative}} 0.
\]
Finally, 
we have the following identification
\[
-K_{Y/C} - \Delta_Y \explshift{-30pt}{\sim_{\bQ}}{\autoref{eq:decomposition_theorem_alternative:adjunction}} (-K_{X/A} - \Delta)|_{Y} \explshift{30pt}{\sim}{$K_{X/A}+\Delta = K_X - f^*K_A + \Delta = K_X + \Delta$} (-K_X - \Delta)|_{Y},
\]
which implies that the divisor $-K_{Y/C}-\Delta_Y$  is nef in the case of assumption \ref{situation:decomposition_theorem_alternative:finite_field} and semi-ample in the case of assumption \ref{situation:decomposition_theorem_alternative:semi_ample}. In either, case it is   nef as required.  

In order to finish the proof of the isomorphism \autoref{eq:isom_main_thm}, we may now directly apply \autoref{lemma:semi_ample_lemma} in the case of assumption \ref{situation:decomposition_theorem_alternative:semi_ample} and \autoref{cor:finite_field} in the case of assumption \ref{situation:decomposition_theorem_alternative:finite_field}.  

We are now ready to finish the proof by getting the torsor structure of $I \to A$. We note that using \autoref{remark:isom_of_two_is_a_torsor} in both cases $I \to A$ is acted on by the $k$-group scheme $G = \Isom((X_0, \Delta_0, L_0))$ and therefore using \cite[Proposition 4.1 a)]{Milne_Etale_cohomology} it suffices to prove that the natural map $G \times I \to I \times_A I$ given by the formula:
\[
(g,\sigma) \mapsto (\sigma,\sigma \circ g)
\] 
is an isomorphism, where:
\begin{align*}
g \in G(T) & = \Isom((X_0, \Delta_0, L_0))(T) \\ 
\sigma \in I(T) & = \Isom_A((X,\Delta;L_X),(X_0 \times_k A, \Delta_0 \times_k A; L_0 \times_k A))(T) 
\end{align*}
for some scheme $T$.  Since $I \to A$ is flat and the relevant isomorphism schemes are compatible with base change this can be verified fiberwise.  This is now clear because $(X,\Delta) \to A$ is already proven to be isotrivial.
\end{proof}

\begin{theorem}[Main decomposition theorem]
\label{thm:decomposition_theorem}
Let $(X,\Delta)$ be a strongly $F$-regular and globally $F$-split projective pair over $k$ such that $K_X + \Delta$ is $\QQ$-Cartier of index coprime to $p$ and
either
\begin{enumerate}[label=(\alph*)]
\item \label{situation:decomposition_theorem:semi_ample} $-K_X-\Delta$ is semi-ample, or
\item \label{situation:decomposition_theorem:finite_field} $k \subseteq \obF_p$ and $-K_X - \Delta$ is nef.
\end{enumerate}  
Let $X \to A$ be the Albanese morphism of $X$, and assume that $A$ admits a $k$-point $0 \in A(k)$.  Then,  
\begin{equation*}
(X,\Delta) \times_A I \cong_I (X_0,\Delta_0) \times_k  I, 
\end{equation*}
where 
\begin{enumerate}
\item \label{itm:decomposition_theorem_alternative:I_def}
\begin{equation*}
I := \Isom_A((X,\Delta;L_X),(X_0 \times_k A, \Delta_0 \times_k A; L_0 \times_k A)),
\end{equation*}
for some relatively ample line bundle $L_X$ on $X$,  $X_0 := f^{-1}(\{0\})$, $\Delta_0:= \Delta|_{X_0}$ and $L_0 := L_{X}|_{X_0}$.
\item \label{itm:decomposition_theorem:surjective} $I \to A$ is surjective and flat.
\item \label{itm:decomposition_theorem:singularities}
$(X_0, \Delta_0)$ is strongly $F$-regular and globally $F$-split such that  $-(K_X + \Delta)$ has Cartier index coprime to $p$.
 \end{enumerate}
In particular,  all fibers $(X_t, \Delta_t)$ of $f$ over $t \in A\big(\ok\big)$ are isomorphic.  Moreover, the map $I \to A$ has a structure of a torsor under $G = \Isom((X_0, \Delta_0, L_0))$ in the flat topology.
\end{theorem}

\begin{proof}
First, just as at the beginning of the proof of \autoref{thm:decomposition_theorem_alternative} we may assume that $k$ is algebraically closed.  Moreover, since $(X,\Delta)$ is $F$-split and the Cartier index of $K_X + \Delta$ is coprime to $p$, we may apply \autoref{prop:albanese_of_fsplit} to see that $f$ is a surjective fibration $F$-split relative to $\Delta$ (see \autoref{sec:relatively_fsplit} or \cite[Def 5.1]{Ejiri_When_is_the_Albanese_morphism_an_algebraic_fiber_space_in_positive___characteristic?} for the definition of $f$ being $F$-split).

Second, \emph{we claim that the general fiber $(X_t,\Delta_t)$ of $f$ satisfies the same assumptions as $(X,\Delta)$.} Since $f$ is $F$-split relative to $\Delta$, all the fibers $(X_t,\Delta_t)$ are $F$-split, and hence reduced.  Consequently, using the assumption on the Cartier index of $-K_X - \Delta$, we may use \autoref{prop:fibres_are_regular_codimension_one} to deduce that the general fiber is regular in codimension one.  Since $X$ is strongly $F$-regular, it is also Cohen--Macaulay (see \cite[Corollary 2.5]{Hochster_Huneke_Tight_closure_and_strong_F-regularity}) and thus the general fiber is Cohen--Macaulay too.  This means that the general fiber satisfies the condition in Serre's criterion for normality and is therefore normal.  This allows us to apply \autoref{thm:general_Fiber_SFR2} to see that the general fiber $(X_t,\Delta_t)$ is strongly $F$-regular.  Finally, we observe that $-K_{X_t} - \Delta_t$ is semi-ample using the adjunction formula described in \autoref{sec:relative_canonical_bundle_adjunction_smooth_base}. This concludes our claim. 

By the above claim,  \autoref{thm:decomposition_theorem_alternative} can be applied to $(X,\Delta)$, yielding our result.
\end{proof}


\section{Finite automorphisms}
\label{sec:finite_automorphisms}

The first author has learned the idea of \autoref{prop:finite_automorphism} and in particular the reference \cite{Rosenlicht_Some_basic_theorems_on_algebraic_groups} from J\'anos Koll\'ar.  Note that the proof works in any characteristic.

\begin{proposition}\label{prop:finite_automorphism}
If $(X,\Delta)$ is a projective klt pair over an algebraically closed field $k$ such that $K_X + \Delta$ is pseudo-effective and $\sL$ is an ample line bundle, then 
\begin{equation*}
\Aut(X,\Delta; \sL):=\{\  \sigma \in \Aut(X,\Delta) \ |\  \sigma^* \sL \cong \sL \ \}
\end{equation*}
is finite. 
\end{proposition}

\begin{proof}
 {\scshape Assume now that $\Aut(X, \Delta ;\sL)$ is not finite.} Write $\Delta= \sum_{s=1}^r c_s \Delta_s$, where $c_i$ are distinct and $\Delta_i$ are reduced. Since $\sL^{\otimes q}$ is very ample, for some   divisible enough integer $q>0$, $\Aut(X,\Delta ;\sL)$ can be identified with the $k$-points of the linear algebraic group 
 \begin{equation*}
 G:= \left\{ \  \left. \alpha \in PGL\left(h^0\left( \sL^{\otimes q}\right), k\right) \ \right| \  \alpha(X) = X,\  \forall s : \alpha(\Delta_s) = \Delta_s \ \right\}.
 \end{equation*}
 In particular, $G$ has infinitely many $k$-points, and hence it is positive dimensional. Hence, it contains an algebraic sub-group $H$ isomorphic either to $\bG_m$ or to $\bG_a$. Hence we have an action $\sigma : H \times X \to X$. The stabilizer sub-group scheme $S \to X$ of $\sigma$ is a closed, proper (so not the whole $H \times X$) sub-group scheme of the group scheme $H \times X \to X$ over $X$. We may assume that $S$ contains no sub-group scheme of the form $H' \times X \to X$ for some finite $k$-subgroup scheme of $H'$ of $H$, as then we could replace $H$ by $H/H'$ (which is still a $1$-dimensional linear algebraic group and hence isomorphic to $\bG_a$ or $\bG_m$).
 
 {\scshape We claim that there is an open set of $X$, where the stabilizer of $\sigma$ is trivial}, that is, the group scheme $S$ introduced above is trivial over this open set.   As $H \times X$ is irreducible, the only way $S$ can be not equal to it is if $\dim S < \dim H \times X = \dim X +1$. Hence, there is an open set, $U_{\fin} \subseteq X$ over which $S$ is finite. However, then over $U_{\fin}$ there is an integer $m$, such that $S|_{U_{\fin}}$ is $m$-torsion. So, $S|_{U_{\fin}}$ is a sub-group scheme of $H[m] \times U_{\fin} \to U_{\fin}$, where $H[m]$ denotes the $m$-torsion sub-group scheme of $H$. In particular, we obtain a morphism $\phi : U_{\fin} \to \Hilb_{k}(H[m] )$.
 Assume now that we are able to prove that there are only finitely many $k$-subgroup schemes $K$ of $H[m]$. Then, $\Hilb_{k}(H[m])$ would be of dimension $0$, and hence, as $U_{\fin}$ is irreducible,  $\phi$ would map it to a single point of $\Hilb_{k}(H[m] )$. That is, $S|_{U_{\fin}} = K \times U_{\fin}$ would hold for some $k$-subgroup scheme $K$ of $H[m]$. Then, using that $S$ is closed in $H \times X$, we would obtain that $K \times X$ is a sub-group scheme of $S$ over $X$ too. This would be a contradiction. 
 
 Hence, to prove the above claim, we only have to prove that $H[m]$ has only finitely many $k$-sub-group schemes.  We remark that the assumption $\dim H = 1$ is necessary (take a look below for an example of infinitely many subgroups of $\alpha_p \times \alpha_p = (\bG_a \times \bG_a)[p]$).  For the proof, we write $m= i p^j$, where $\gcd(i,p)=1$ if $\charfield k =p>0$ and $i:=m$ otherwise. Let $\tH$ be an arbitrary sub-group scheme of $H[m]$. As we are over a prefect field (in fact, even algebraically closed), $\tH$ splits as $\tH_{\et} \times \tH_{\inf}$, where $\tH_{\et}$ is \'etale and hence a subgroup(scheme) of $H[i]$ and $\tH_{\inf}$ is connected and hence  a sub-group scheme of $H[p^j]$.  So, it is enough to prove that $H[i]$ and $H[p^j]$ have finitely $k$-subgroup schemes separately. This is immediate for the former, however one has to be slightly careful with the latter as there are finite connected group schemes over $k$ that have infinitely many subgroup schemes, e.g., $\alpha_p \times \alpha_p$ admits a family of distinct embeddings $\alpha_p \to \alpha_p \times \alpha_p$ parametrized by the projective space of its Lie algebra. Luckily, $H[p^j] \cong \alpha_{p^j}$ or $\mu_{p^j}$, and then it has a unique chain of subgroup schemes given by the $p^l$-torsion subgroup schemes for $1 \leq l \leq j$.  To make a precise argument, we proceed by induction on $j$  proving that $H[p^j]$ has only finitely may subgroup schemes. For $j=1$ this is immediate by dimension reasons. Then, it is enough to prove that every non-zero sub-group scheme $\tH \subseteq H[p^j]$ contains $H[p]$, as then $\tH/H[p]$ becomes a subgroup scheme of $H[p^j]/H[p] \cong H[p^{j-1}]$ (the isomorphism given by multiplication by $p$),  for which we already know the statement. However, the latter statement is quite straightforward: let $l>0$ be the smallest integer such that $p^l \cdot \tH = \{0\}$. Then, $p^{l-1} \cdot \tH$ is a non-zero subgroup scheme of both $\tH$ and of $H[p]$, and by dimension reasons then it has to be equal to the latter. This finishes the proof of our claim.
 
Note now that $U:=X \setminus \Supp \Delta$ is necessarily $H$ invariant. By intersecting it with the open set found in the above claim,  we may also assume that $H$ acts freely on an open set $U$ contained in $X \setminus \Supp \Delta$. Further, by \cite[Theorem 2 and 10]{Rosenlicht_Some_basic_theorems_on_algebraic_groups} or
 \cite{Rosenlicht_Another_proof_of_a_theorem_on_rational_cross_sections}, we may assume that $V:=U/H$ exists as a scheme and that there is a section $s \colon V \to U$. However, then $U \cong V \times H$ necessarily. This yields a birational morphism $X \dashrightarrow \bP^1 \times Y$, where $Y$ is a normal compactification of $V$, $U$ is in the domain of $f$, and $U$ is mapped isomorphically to $H \times V$ for some fixed embedding $H \subseteq \bP^1$. Let $Z$ be the normalization of the graph of this map. That is, we have the following commutative diagram, where $f$ and $g$ are birational, proper and the other morphisms are open embeddings, over the images of which $f$ and $g$ are isomorphisms:
 \begin{equation*}
 \xymatrix{
& Z \ar[dl]_f \ar[dr]^g \\
 X & \ar@{_(->}[l]^{\iota_X} \ar@{_(->}[u]^{\iota_Z}  U \cong  H \times V \ar@{^(->}[r]_j & \bP^1 \times Y
 }
 \end{equation*}
 One can write
 \begin{equation}
 \label{eq:finite_automorphism:log_canonical}
  f^* (K_X + \Delta)  = K_Z + \Gamma
 \end{equation}
 for compatible choices of $K_X$ and $K_Z$. Let $C$ be a general $\bP^1$ in $\bP^1 \times Y$.
Then:
 \begin{multline}
 \label{eq:finite_automorphism:Kodaira_dimension}
 K_X + \Delta \textrm{ is pseudo-effective }
 \expl{\Rightarrow}{\autoref{eq:finite_automorphism:log_canonical}}
 K_Z + \Gamma \textrm{ is pseudo-effective }
\\ \explshift{-40pt}{\Rightarrow}{$K_{\bP^1 \times Y} + g_* \Gamma= g_* (K_Z + \Gamma)$}
 K_{\bP^1 \times Y} + g_* \Gamma \textrm{ is pseudo-effective }
\explparshift{120pt}{-120pt}{\Rightarrow}{$C$ is the general element of a moving family of curves}
( K_{\bP^1 \times Y} + g_* \Gamma) 
\explparshift{165pt}{30pt}{\cdot}{as $C$ is a general elment of the ruling of $\bP^1 \times Y$, it avoids every codimension $2$ subset of $\bP^1 \times Y$, in particular, it avoids the non $\bQ$-Cartier locus of $K_{\bP^1 \times Y} + g_* \Gamma$ }
 C \geq 0 
 \hspace{60pt}
 \end{multline}
  As  $f$ and $g$ are  isomorphisms over $U$ and $U \cap \Supp \Delta= \emptyset$, we obtain that
 \begin{equation}
 \label{eq:finite_automorphism:avoid}
 g(\supp \Gamma) \subseteq \bP^1 \times Y \setminus j(H \times V).
 \end{equation}
 As $(X, \Delta)$ is klt, all coefficients of $\Gamma$ are smaller than $1$. By \autoref{eq:finite_automorphism:avoid}, the only possible components of $g_* \Gamma$ that $C$ meets are the ones that are the components of $(\bP^1 \setminus H) \times Y$. As $C$ meets there are at most two such components, and $C$ meets them in multiplicity $1$, we obtain that
 \begin{equation*}
  ( K_{\bP^1 \times Y} + g_* \Gamma) \cdot C = \deg (K_{\bP^1} + g_* \Gamma) \cdot C < -2 + 2 \cdot 1 = 0.
  \end{equation*}
This contradicts \autoref{eq:finite_automorphism:Kodaira_dimension}. Hence, our assumption was false, and $\Aut(X, \Delta; \sL)$ is finite. 
\end{proof}

\section{Weak Bogomolov--Beauville decomposition}
\label{sec:bogomolov_beauville_actual_proofs}

In this chapter we apply the results from two previous chapters to get the weak Beauville--Bogomolov decomposition for weakly ordinary varieties with trivial canonical class.  First, we generalize to the setting of \autoref{thm:decomposition_theorem} the necessary definition which appeared before in \cite{Greb_Kebekus_Petternel_Singular_Spaces_with_Trivial_Canonical_Class}.  We say that a morphism is \emph{quasi-\'etale} if it is \'etale in codimension one.

\subsection{Quasi-\'etale maps and augmented irregularity}

Before defining the augmented irregularity in our situation, we need a lemma:

\begin{lemma}
\label{lemma:properties_preserved_etale_codim_one}
Let $(X,\Delta)$ be a normal pair, and let $f \colon X' \to X$ be a quasi-\'etale morphism of normal varieties.  Assume that $(X,\Delta)$ satisfies one of the properties:
\begin{enumerate}[1)]
\item sharply $F$-pure,
\item globally $F$-split,
\item strongly $F$-regular,
\item globally  $F$-regular.
\end{enumerate}
Then the pair $\left(X',\Delta'\right)$, where $\Delta' := f^*\Delta$, satisfies the same property.
\end{lemma}

\begin{proof}

We first prove the results concerning sharp $F$-purity and global $F$-splitting.  Since sharp $F$-purity is just a local version of global $F$-splitting, we may focus on latter notion.  We take a Frobenius splitting $s \colon F^e_*\cO_X(\lceil (p^e - 1) \Delta \rceil) \to \cO_X$ of the pair $(X,\Delta)$.  Let $f_U \colon U' \to U$ be a restriction of $f$ to the intersection of the regular locus of $X$ and the \'etale locus of $f$.  By restricting the splitting $s$ to $U$ and taking the pullback along $f_U^*$ we obtain a Frobenius splitting 
\[
f_U^*(s_{|U}) \colon F^e_*\cO_{U'}\left(\left\lceil (p^e - 1)\Delta' \right\rceil\right) \expl{\isom}{flat base change applies because $f_U$ is \'etale over a regular scheme} f_U^*F^e_*\cO_{U}(\lceil (p^e - 1)\Delta \rceil) \xrightarrow{f_U^*(s_{|U})} f_U^*\cO_U \isom \cO_{U'}
\] defined on the open subset $U'$.  Since the relevant sheaves are reflexive, we may extend $s_{U'}$ to a splitting on the whole $X'$, which concludes this part of the proof.

We now proceed to the statements concerning $F$-regularity.  As above it suffices to approach the global case.  In this situation, we take a Cartier divisor $C \subset X$ containing the non-\'etale locus of $f$ and such that $X \setminus C$ is smooth and affine (and hence globally strongly $F$-regular).  We denote by $C'$ the preimage $f^*C$.  By \cite[Theorem 3.9]{Schwede_Smith_Globally_F_regular_an_log_Fano_varieties} it suffices to prove that $X' \setminus C'$ is globally strongly $F$-regular and the map 
\begin{equation}
\label{eq:properties_preserved_etale_codim_one:upstairs}
\cO_{X'} \to F^e_*\cO_{X'}\left(\left\lceil (p^e - 1)\Delta' \right\rceil + C'\right)  
\end{equation}
 splits.   The former statement is clear since $X' \setminus C'$ is affine and smooth.  In order to see the latter, by using \cite[Theorem 3.9]{Schwede_Smith_Globally_F_regular_an_log_Fano_varieties} again we obtain that $F^e_*\cO_{X}(\left\lceil (p^e - 1)\Delta \right\rceil + C) \to \cO_{X}$ splits.  Then as above by  taking pullbacks of this splitting  along the  \'etale locus we obtain the splitting of \autoref{eq:properties_preserved_etale_codim_one:upstairs}.
\end{proof}

We recall that using \cite[Theorem 1.1]{Ejiri_When_is_the_Albanese_morphism_an_algebraic_fiber_space_in_positive___characteristic?} the Albanese morphism of a pair $(X,\Delta)$ satisfying the conditions of \autoref{thm:decomposition_theorem} is surjective. Furthermore, by \autoref{lemma:properties_preserved_etale_codim_one}, these conditions hold also for quasi-\'etale covers.  This justifies the existence of the maximum in the following definition, which is in fact bounded from above by $\dim X$.

\begin{definition}[Augmented irregularity]
\label{def:augmented_irregularity}
Let $(X,\Delta)$ be  a strongly $F$-regular and globally $F$-split projective pair over $k$ such that the Cartier index of $K_X+\Delta$  coprime to $p$.  We define the \emph{augmented irregularity} $\wh{q}(X)$ by the formula:
\[
\wh{q}(X) = \max \left\{ \left. \dim \Alb_{X'} \right|  X' \to X \text{ is a finite quasi-\'etale morphism of normal varieties}\right\}.
\] 
Note that if $X$ happens to be smooth then it is enough to consider only \'etale morphisms in the above definition by the purity of branch locus. 
\end{definition}

\begin{proposition}
\label{lem:invariance_augmented_irregularity}
Let $X$ and $Y$ be normal projective varieties.  Assume that $X \to Y$ is either a finite quasi-\'etale map or a universal homeomorphism.  Then $\wh{q}(X) = \wh{q}(Y)$.
\end{proposition}
\begin{proof}
The statement concerning quasi-\'etale maps is immediate from \autoref{def:augmented_irregularity}.  We therefore proceed to the case when $X \to Y$ is a universal homeomorphism.  First we claim that the normalized pullback map induces an equivalence between the categories of finite quasi-\'etale normal coverings of $X$ and $Y$.  Indeed, let $Y'$ be a normal variety and let $Y' \to Y$ be a finite quasi-\'etale map.  The scheme $X'$ defined as a normalization of the pullback $X \times_Y Y'$ is normal and admits a finite quasi-\'etale map to $X$.  This yields a functor in one direction.  In order to describe an inverse, we take a finite quasi-\'etale covering $\pi \colon X' \to X$.  By \stacksproj{04DY} we see that $\pi$ can be uniquely descended to a finite \'etale covering $\pi_U \colon U' \to U$ of a big open subset $U \subset Y$.  The morphism $\pi_U$ can be uniquely extended to a covering of $Y$ by taking the normalization of $Y$ inside the fraction field of $U$.  This yields the required inverse functor.  We leave the details to the reader.  We conclude the proof of the proposition by observing that for a pair of finite quasi-\'etale maps $X' \to X$ and $Y' \to Y$ associated to each other via the above equivalence the natural map $X' \to Y'$ is a universal homeomorphism, and hence $\dim \Alb_{X'} = \dim \Alb_{Y'}$ by \autoref{cor:univ_homeo_alb_dim}.
\end{proof}

\subsection{The proof of the weak Beauville--Bogomolov decomposition}

\noindent In the proof of the weak Beauville--Bogomolov decomposition for $-K_X-\Delta$ numerically trivial we use the following result which implies that the torsors arising from \autoref{thm:decomposition_theorem} are particularly simple, that is, they are dominated by torsors induced by $n$-torsion points, for $n \in \NN$.

\begin{proposition}[{\cite[Proposition]{Nori_Fundamental_Group_Scheme_Of_An_Abelian_Variety}}]
\label{prop:torsors_on_abelian_vars}
Let $A$ be an abelian variety over a field $k$, let $G$ be a finite group scheme over $k$, and let $P \to A$ be a $G$-torsor with a $k$-rational point over $0 \in A$.  Then there exists an integer $n \in \NN$, a homomorphism of group schemes $\phi : A[n] \to G$, and a morphism $A \to P$,  which is $A[n]$-equivariant with the $A[n]$ action on $P$ induced by $\phi$.
\end{proposition}

The statement of \autoref{prop:torsors_on_abelian_vars} is summarized in the following diagram:
\begin{equation*}
\xymatrix@C=90pt@R=8pt{
A \ar[dr]_{n_A} \ar[r]^{A[n] \textrm{ equivariant}} & P \ar[d] \\
& A
}
\end{equation*}
Before stating \autoref{thm:beauville_bogomolov_general_logcy}, we note the following lemma stating that the pairs in \autoref{thm:beauville_bogomolov_general_logcy} automatically have log-canonical divisors with Cartier index prime-to-$p$. 

\begin{lemma}
\label{lem:index_prime_to_p}
Let $(X, \Delta)$ be a globally $F$-pure pair such that $K_X + \Delta \equiv 0$, then the Cartier index of $K_X + \Delta$ is prime-to-$p$.
\end{lemma}

\begin{proof}
By \autoref{prop:technical_F_splittings_and_section}, a global Frobenius splitting $F^e_*\cO_X(\lceil (p^e-1)\Delta \rceil) \to \cO_X$ of $(X,\Delta)$ corresponds to a $\Gamma$ such that
\begin{equation}
\label{eq:index_prime_to_p}
0 \leq \Gamma \sim  \lfloor(1-p^e)(K_X + \Delta ) \rfloor =  \underbrace{(1-p^e)(K_X + \Delta)}_{\equiv 0} + \Big(\underbrace{\lfloor (1-p^e) \Delta \rfloor - (1-p^e) \Delta}_{=: D, \textrm{ where } D \leq 0}\Big) 
\end{equation} 
Hence, in the sum $ \Gamma + (-D)$ both summands are effective divisors, and the sum itself is numerically trivial. This implies that $\Gamma= D=0$, and then by \autoref{eq:index_prime_to_p}, $(1- p^e) (K_X + \Delta) \sim 0$.
\end{proof}

\begin{theorem}[Weak Beauville--Bogomolov decomposition]
\label{thm:beauville_bogomolov_general_logcy}
Let $(X,\Delta)$ be a globally $F$-split pair with strongly $F$-regular singularities  such that $-K_X-\Delta \equiv 0$.  Then there exists a composition 
\begin{equation*}
Y \to W \to X
\end{equation*}
such that
\begin{enumerate}
\item $W \to X$ is  a finite quasi-\'etale morphism,
\item $Y \to W$ is a finite infinitesimal torsor  under $G^0=\prod_{i=1}^{\wh{q}(X)} \mu_{p^{j_i}}$ for some integers $j_i \geq 0$, and
\item $(Y,\Delta_Y) \isom (Z_0,\Delta_0) \times B$ where
\begin{enumerate}
\item  $B$ is an abelian variety of dimension $\wh{q}(X)$,and
\item $(Z_0,\Delta_0)$ is a  globally $F$-split pair with  strongly $F$-regular singularities such that $K_{Z_0}+\Delta_{Z_0} \equiv 0$ and $\wh{q}(Z_0) = 0$.
\end{enumerate}
\end{enumerate}
Additionally we can assume that $G^0$ acts diagonally on the factors and that the action is faithful on $Z_0$ and free on $B$.
\end{theorem}

\begin{proof}
We show the statement using \autoref{thm:decomposition_theorem}.  First,  take a quasi-\'etale morphism $\pi \colon Z \to X$ such that $Z$ is normal and $\dim \Alb_Z = \wh{q}(X)$. Note that by \autoref{lemma:properties_preserved_etale_codim_one}, $(Z, \pi^* \Delta)$ is globally $F$-split. Hence, by \autoref{prop:albanese_of_fsplit}, $\Alb_Z$ is ordinary. In particular,  by \cite[Prop 2.3.2]{Hacon_Patakfalvi_A_generic_vanishing_in_positive_characteristic},  $\Alb_Z$ has $p$-rank $g:=\dim \Alb_Z $. Hence, by
\cite[``The $p$-rank'', p. 146]{Mumford_Abelian_varieties} $\left(\Alb_Z \otimes_k \overline{k}\right)[p^r]=\left(\factor{\bZ}{p^r\bZ}\right)^{\times g} \times \mu_{p^r}^{\times g}$ (the third factor mentioned in \cite{Mumford_Abelian_varieties} has height $1$ when the $p$-rank equals $g$ \cite[page 147, line -7]{Mumford_Abelian_varieties}). In particular $\left(\left(\Alb_Z \otimes_k \overline{k}\right)[p^r]\right)^0\cong \mu_{p^r}^{\times g}$.

As $k$ is perfect, we can perform any additional finite base-extensions, and we can count them into the quasi-\'etale part of the cover we construct. Hence, we may assume that $Z$ has a rational point (and hence the Albanese variety has one) and we may also assume that for the identity component of the $p^r$ torsion subgroup of $\Alb_Z$ we have $(\Alb_Z[p^r])^0 \cong  \mu_{p^r}^{\oplus \dim \Alb_Z}$, for some integer $r>0$ specified later. 

We consider the Albanese morphism $f \colon Z \to \Alb_Z$, where $\Alb_Z$ is an abelian variety.   By \autoref{lemma:properties_preserved_etale_codim_one} we see that the assumptions concerning $(X,\Delta) $ are also satisfied for $(Z,\Delta_Z)$, where $\Delta_Z = \pi^*\Delta$, and therefore by \autoref{thm:decomposition_theorem}, there exists a torsor in the flat topology $I \to \Alb_Z$ such that 
\[
(Z,\Delta_Z) \times_{\Alb_Z} I \isom I \times (Z_0,\Delta_0),
\] 
where $(Z_0,\Delta_0)$  possesses all the required properties.  More precisely, the map $I \to \Alb_Z$ is in fact a torsor under $\Isom\left(Z_0,\Delta_{Z_0}, L_0\right)$, for some ample polarization $L$ on $Z$, which is a finite group scheme by \autoref{prop:finite_automorphism}. Additionally, from the precise description of $I$ as an $\Isom$ scheme in \autoref{thm:decomposition_theorem}, using the base-change properties of $\Isom$ stated in \autoref{prop:base_change_for_isom_scheme}, we see that $I$ has a $k$-rational point over $0 \in A(k)$ given by the identity automorphism.
Therefore, by \autoref{prop:torsors_on_abelian_vars} there exists a torsor $B \to \Alb_Z$ under a finite abelian group-scheme $G$  and a morphism $B \to I$ that is $G$-equivariant via a group-scheme homomorphism $G \to \Isom\left(Z_0,\Delta_{Z_0}, L_0\right)$. Additionally, the total space of $B$ is also an abelian variety and $G=B[r]$ for some $r>0$, which is the integer $r$ that was promised at the beginning of the proof to get fixed later. In particular, the identity component of $G$ is isomorphic to $\prod_{i=1}^{\wh{q}(X)} \mu_{p^{j_i}}$ for some integers $j_i \geq 0$, which is the assumption we keep instead of insisting on $G= B[r]$.  Replacing now both $B$ and $G$ by the quotient with  $\Ker \big(G \to \Isom\left(Z_0,\Delta_{Z_0}, L_0\right)\big)$ we may also assume that $G \to \Isom\left(Z_0,\Delta_{Z_0}, L_0\right)$ is injective.   Set now $Y:= Z \times_{\Alb_Z} B$. As $B \to \Alb Z$ factors through $I \to \Alb Z$, $I \times_{\Alb_Z} B \to B$ has a section. Then, using \autoref{prop:base_change_for_isom_scheme} and the precise description of $I$ as an $\Isom$ scheme from \autoref{thm:decomposition_theorem}, we obtain that $(Y, \Delta_Y) \cong_B \left( Z_0, \Delta_{Z_0} \right) \times B$.  Let us denote this isomorphism by $\phi : Y \to Z_0 \times B$.

Our next task is to define $W$.
The group scheme $G$ is abelian and can be decomposed as a direct sum $G =  G^0 \oplus G^{\et}$, where $G^0$ is infinitesimal and $G^{\et}$ is \'etale.     We set $W:= Y /G^0$. As $W \to Z$ is a quotient by $G^{\et}$, $W \to Z$ is \'etale, and hence $W \to X$ is quasi-\'etale. We also have $\dim B = \dim \Alb_Z = \wh{q}(X)$.

To show $\wh{q}(Z_0) = 0$, consider first the following line of equalities 
\begin{equation}
\label{eq:beauville_bogomolov_general_logcy:augmented_irregularity}
\dim B = \wh{q}(X)
\explparshift{180pt}{-175pt}{ =}{\autoref{lem:invariance_augmented_irregularity}, as $W \to X$ is quasi-\'etale}
\wh{q}(W) 
\explparshift{150pt}{140pt}{ =}{\autoref{lem:invariance_augmented_irregularity}, as $Z_0 \times 
B \to W$ is a universal homeomorphism}
 \wh{q}(Z_0 \times B)
\end{equation}
Assume now that $\wh{q}(Z_0)>0$ holds. Then, there is a quasi-\'etale cover $V \to Z_0$ with $q(V)>0$. Hence,  $V \times B$ is a quasi-\'etale cover of $Z_0 \times B$ with $q(V \times B) > \dim B$. This contradicts \autoref{eq:beauville_bogomolov_general_logcy:augmented_irregularity}.

Concerning the addendum: the $G$-action on  $Z_0 \times B$ is defined by passing to $Z \times_{\Alb_Z} B$ via $\phi^{-1}$, then applying the $G$-action there on the $B$-coordinate, and then returning to $Z_0 \times B$ via $\phi$. So let us start with $(z,b) \in Z_0 \times B$ and trace its image through this process (more precisely we should consider also a test-scheme $S$ with which we take product of every space, but let us omit that for simplicity of notation). Note that $\phi$ is an isomorphism over $B$, so it restricts to an isomorphism $\phi_b$ on the fibers over every $b \in B$. By applying $
\phi^{-1}$ to $(z,b)$  we get $\Big(\phi^{-1}_b(z),b\Big)$. Then, $g\in G$ takes $\Big(\phi^{-1}_b(z),b\Big)$ to $\Big(\phi^{-1}_b(z),gb\Big)$. Getting back to the $Z_0 \times B$ side, we see that $(z,b)$ is taken by $g$ during the overall process to $\bigg(\phi_{gb}\Big(\phi^{-1}_b(z)\Big),gb\bigg)$. The key point is that $\phi_{gb}= g \circ \phi_b$, because the trivialization of $I \times_{\Alb_Z} B$ is based on the diagonal section. Hence, $g \in G$ takes $(z,b)$ to $(gz,gb)$. So, the action of $G$ is indeed diagonal on $Z_0 \times B$. Hence, the same is true for the action of $G^0 \subseteq G$. 

Additionaly, $G^0$ acts freely on $B$ as $G^0 \subseteq G$, and $G$ acts freely on $B$.  Furthermore, $G^0$ acts faithfully on $Z_0$ as it is a subgroup of $\Isom\left(Z_0,\Delta_{Z_0}, L_0\right)$.
\end{proof}

\begin{proof}[Proof of \autoref{cor:beauville_bogomolov}]
By \autoref{prop:ordinary_is_fsplit} weakly ordinary varieties with trivial canonical class are globally $F$-split, and hence \autoref{thm:beauville_bogomolov_general_logcy} can be directly applied for $\Delta = 0$.
\end{proof}

\begin{remark}
\label{rem:domination}
We claim that every two decompositions as in \autoref{thm:beauville_bogomolov_general_logcy} can be dominated via a composition of a quasi-\'etale map and an abelian infinitesimal torsor by a third one, if we relax the condition stated in the addendum of \autoref{thm:beauville_bogomolov_general_logcy} about the free/faithful action on the respective factors. So, take two such decomposition $Z_{0}^1 \times B^1$ and $Z_{0}^2 \times B^2$ as in the following diagram. The other parts of the diagram are constructed after the diagram itself, and the numbers on the arrows denote the number of the step of this construction in which the given arrow is constructed.
\begin{equation}
\label{eq:domination}
\xymatrix@C=45pt{
B^1 \ar[d]|{\textrm{\autoref{itm:domination:initial_Albanese}}} 
& & B \ar[rr]|{\textrm{\autoref{itm:domination:from_middle_Albanese_top}}} \ar[ll]|{\textrm{\autoref{itm:domination:from_middle_Albanese_top}}} \ar[d]|{\textrm{\autoref{itm:domination:middle_Albanese}}}  & 
& B^2 \ar[d]|{\textrm{\autoref{itm:domination:initial_Albanese}}} \\
%
%
\Alb_{W^1} & & \Alb_W \ar[ll]|{\textrm{\autoref{itm:domination:from_middle_Albanese}}}  \ar[rr]|{\textrm{\autoref{itm:domination:from_middle_Albanese}}}   & & \Alb_{W^2} \\
%
%
 & \ar[uul]|(0.3){\textrm{\autoref{itm:domination:initial_Albanese}}} Z_{0}^1 \times B^1 \ar[d]_{\rm inf. ab}^{\alpha^1} 
& Z_0 \times B \ar@/_3pc/[uu]|(0.7){\textrm{\autoref{itm:domination:middle_Albanese}}}  \ar[l]|{\textrm{\autoref{itm:domination:product}}}_(0.4){\gamma^1} \ar[r]|(0.6){\textrm{\autoref{itm:domination:product}}}^(0.7){\gamma^2} \ar[d]_{\rm inf. ab \quad}|{\textrm{\autoref{itm:domination:middle_covers}}}  
& Z_{0}^2 \times B^2 \ar[d]^{\rm inf. ab}_{\alpha^2} \ar[uur]|(0.3){\textrm{\autoref{itm:domination:initial_Albanese}}} &  \\
%
%
& W^1  \ar[uul]|(0.7){\textrm{\autoref{itm:domination:initial_Albanese}}} \ar[dr]_{q.-\et}^{\beta^1} 
& W \ar[d]|{\textrm{\autoref{itm:domination:middle_covers}}} \ar[r]|{\textrm{\autoref{itm:domination:middle_covers}}} \ar[l]|{\textrm{\autoref{itm:domination:middle_covers}}}  \ar@/_3pc/[uu]|(0.3){\textrm{\autoref{itm:domination:middle_Albanese}}}  & W^2 \ar[ld]^{q.-\et}_{\beta^2}  \ar[uur]|(0.7){\textrm{\autoref{itm:domination:initial_Albanese}}} \\
%
%
& & X.
} 
\end{equation}
We construct the commutative diagram \autoref{eq:domination} through the following steps:
\begin{enumerate}
\item
\label{itm:domination:initial_Albanese}
As $\alpha_i$ are inseparable, $\dim \Alb_{W^i}= \dim \Alb_{Z^i_0 \times B^i} = \dim B_i$. Hence, using the functoriality of the Albanese morphisms we obtain the isogenies $B^i \to \Alb_{W^i}$ as in \autoref{eq:domination}.
\item \label{itm:domination:middle_covers} As the proof of \autoref{thm:beauville_bogomolov_general_logcy} starts by taking a high enough quasi-\'etale cover, we may find a composition $Z_0 \times B \to W \to X$ as in  \autoref{thm:beauville_bogomolov_general_logcy}, such that $W$ dominates both $W^1$ and $W^2$.
\item \label{itm:domination:middle_Albanese} As in point \autoref{itm:domination:initial_Albanese} we obtain the isogeny $B \to \Alb_W$.
\item \label{itm:domination:from_middle_Albanese} By the universal property of Albanese morphisms we obtain the isogenies $\Alb_W \to \Alb_{W^1}$ and $\Alb_W \to \Alb_{W^2}$.
\item \label{itm:domination:from_middle_Albanese_top} By possibly increasing $B \to \Alb_W$ by a higher inseparable isogeny, we obtain also the isogenies $B \to B^1$ and $B \to B^2$.
\item \label{itm:domination:product} As $Z_0^i \times B^i \cong W^i \times_{\Alb_{W^i}} B^i$, we obtain the morphisms $\gamma^i$ .
\item It is immediate from the construction, that $\gamma^i$ map the fibers of the projection $Z_0 \times B \to B$ to the fibers of the projection $Z_0^i \times B \to B^i$, as these projections, $\gamma_i$ and the isogenies $B \to B^i$ commute. However, it is not clear from the construction if there are similar morphisms $Z_0 \to Z_0^i$ with which $\gamma^i$ and the corresponding projections commute. 
\end{enumerate}
\end{remark}


\section{Corollaries}

We start with a short summary of the classification of smooth surfaces with trivial canonical class defined over an algebraically field of characteristic $p>0$ provided by Bombieri and Mumford,  with a special emphasis on the relation between augmented irregularity and $F$-splitting (see \cite{Liedtke_Algebraic_surfaces_in_positive_characteristic,Bombieri_Mumford_Enriques_classification_of_surfaces_in_char_p_II,Bombieri_Mumford_Enriques_classification_of_surfaces_in_char_p_III}).

\begin{proposition}\label{prop:classification_ordinary_cy}
Let $S$ be a smooth surface defined over an algebraically closed field of characteristic $p>0$ with trivial canonical class.  Then $S$ belongs to one of the following classes of surfaces having the described properties:
\begin{enumerate}
    \item abelian surfaces,
    \item K3 surfaces are simply connected and hence satisfying $\wh{q}(S) = 0$, with the invariants:
        \begin{align*}
            h^1(S,\cO_S) = 0, \quad b_1(S) = 0, \quad b_2(S) = 22,
        \end{align*}
    \item \label{itm:classification_ordinary_cy:Enriques} non-classical Enriques surfaces in characteristic $p = 2$, with invariants:
            \begin{align*}
                h^1(S,\cO_S) = h^2(S,\cO_S) = 1, \quad b_1(S) = 0, \quad b_2(S) = 10
            \end{align*}
        the surfaces are called singular (resp. supersingular) if equivalently the Frobenius action on $H^2(S,\cO_S)$ is bijective (resp. zero) or the surface is globally $F$-split (resp. not globally $F$-split).  The fundamental group of an Enriques surface is non-trivial if and only the surface is singular, it is then isomorphic to $\ZZ/2\ZZ$.
    \item hyperelliptic surfaces, given as a quotient of a product of two elliptic curves by an action of a potentially non-reduced abelian group scheme, with invariants:
        \begin{align*}
            h^1(S,\cO_S) = 2, \quad p_g(S) = 1, \quad b_1(S) = b_2(S) = 2,
        \end{align*}
        and hence satisfying with $\wh{q}(S) \geq 1$.
        \item a quasi-hyperelliptic surfaces  in characteristic $p=2,3$, given as a quotient of a product of $\bP^1$ and an elliptic curve by a finite group scheme. We have $b_1(S) = 2$, $b_2(S) =2$, $\wh{q}(S)=1$, and $S$ admits a fibration $S \to E$ to an elliptic curve such that the fibers are all rational curves with a single cusp singularity. As general fibers of fibrations of globally $F$-split varieties are globally $F$-split \cite[Cor 2.5]{Gongyo_Li_Patakfalvi_Schwede_Tanaka_Zong_On_rational_connectedness_of_globally_F_regular_threefolds}, and since cusps are not $F$-pure, we see that $S$ is not globally $F$-split. 
\end{enumerate}
\end{proposition}

\begin{lemma}
\label{lemma:decomp_rel_dim_2}
Let $X$ be a smooth projective $F$-split variety with numerically trivial canonical class and $\wh{q}(X) = \dim X - 2$ (e.g. a threefold with $\wh{q}(X) = 1$).  Then there exists a generically finite map $A \times Z \to X$ where $A$ is an abelian variety of dimension $\wh{q}(X)$ and $Z$ is a globally $F$-split K3 surface.
\end{lemma}

\begin{proof}
We use \autoref{cor:beauville_bogomolov} to see that there exists a composition of an \'etale map and a finite universal homeomorphism $X' \to X$ such that $X' \isom A \times Y$ where $A$ is an abelian variety of dimension $\wh{q}(X)$ and $Y$ is an $F$-split Gorenstein strongly $F$-regular surface with trivial canonical class and $\wh{q}(Y) = 0$.  In particular, $Y$ has DuVal singularities. Let $Y' \to Y$ be the minimal resolution. As the exceptional divisors of $Y' \to Y$ are all smooth rational curves, we show that $\wh{q}(Y') = \wh{q}(Y)=0$: Let $W'$ be a quasi-\'etale cover of $Y'$. By purity of branch locus $W' \to Y'$ is \'etale, and hence the preimages of the exceptional rational curves are disjoint unions of smooth rational curves. Let $W$ be the normalization of $Y$ in the function field of $W'$. This is a finite cover $Y$ which agrees with $W'$ over the isomorphism locus of $Y' \to Y$. Hence 
$W' \to W$  contracts exactly the above mentioned rational curves. Therefore, the birational contraction $W' \to W$ factors the morphism $\alb_{W'}$. This implies that we have equality for the spaces $\Alb_{W'} = \Alb_W$. However, by the $\wh{q}(Y)=0$ assumption the latter is a point, and therefore $\Alb_{W'}$ is also a point. This concludes that $\wh{q}(Y')=0$. 

By \autoref{prop:classification_ordinary_cy}, $Y'$ is either a K3 surface or a singular Enriques. In the first case, set $Z:=Y'$ and the second one set $Z \to Y'$ to be the $\bZ/2 \bZ$ Galois \'etale cover stated to exist in \autoref{prop:classification_ordinary_cy},. Note that in the latter case, $Z$ is again a K3 surface, and it is globally $F$-split by \autoref{lemma:properties_preserved_etale_codim_one}.\autoref{itm:classification_ordinary_cy:Enriques}.   This concludes our proof.
\end{proof}

\subsection{Fundamental groups} 

In characteristic zero, as direct corollaries of the Beauville--Bogomolov decomposition, one obtains the following:
\begin{corollary}[of the Beauville--Bogomolov decomposition in char. zero]
\label{cor:virtually_abelian}
Let $X$ be a smooth projective variety with numerically trivial canonical class over $\bC$.  Then
\begin{enumerate}
\item $\pi_1(X)$ and $\piet(X)$  are virtually abelian -- the abelian part comes from the abelian variety in the decomposition,
\item and consequently the condition $\wh{q}(X) \neq 0$ is equivalent to $|\pi_1(X)|=|\piet(X)| = \infty$.
\end{enumerate}
\end{corollary}

\noindent It is natural to expect that similar behaviour can be exhibited in characteristic $p>0$ under some favorable arithmetic assumptions, for example some version of ordinarity.  Unfortunately, we have not been able to get any unconditional results in this direction.  However, using the corollaries of \autoref{thm:beauville_bogomolov_general_logcy}, we can show that the above conditions hold true for weakly ordinary varieties with numerically trivial class and the augmented irregularity satisfying $\wh{q}(X) \geq \dim X - 2$.  In the proof we need the following:

\begin{lemma}[{\cite[Lemme 4.4.17]{Deligne_Theorie_de_Hodge_II}}]
\label{lem:fundamental_group_resolutions}
Let $f \colon Y \to X$ be a generically finite morphism of normal varieties.  Then the image of $f_* \colon \piet(Y) \to \piet(X)$ is of finite index.
\end{lemma}

\begin{corollary}
\label{cor:fundamental_group_virtually_abelian}
Let $X$ be a smooth projective $F$-split variety with numerically trivial canonical class and $\wh{q}(X) \geq \dim X - 2$ (e.g. a threefold with $\wh{q}(X) \neq 0$) defined over an algebraically closed field $k$.  Then the fundamental group $\piet(X)$ is virtually abelian.
\end{corollary}

\begin{proof}
The only non-trivial case is when $\wh{q}(X) = \dim X - 2$.  In this situation, by \autoref{lemma:decomp_rel_dim_2} there exists a generically finite map $A \times Z \to X$, where $A$ is an abelian variety and $Z$ is a K3 surface.  This gives the statement by \autoref{lem:fundamental_group_resolutions}.
\end{proof}

It would be interesting to obtain \autoref{cor:virtually_abelian} in characteristic $p>0$ (or even only for ordinary varieties) without appealing to any Beauville--Bogomolov type results.
\begin{question}
Let $X$ be a smooth projective variety with trivial canonical class defined over $k$.  Is $\piet(X)$ virtually abelian? 
\end{question}
As a particular case, one would like to exclude the $p$-completion of a free group on $n \geq 2$ generators as the fundamental group of a $K$-trivial variety.

\subsection{Betti numbers}

In this section we prove that the Betti numbers of $F$-split threefolds with $\hat q(X) \geq 1$ and with numerically trivial canonical class are bounded from above.  Before proceeding with the actual argument we recall the following classical result.

\begin{lemma}[{\cite[Proposition 1.2.4]{Kleiman_Algebraic_cycles_and_the_Weil_conjectures}}]\label{lem:generically_finite_etale_cohomology_inj}
Let $f \colon X \to Y$ be a surjective morphism of smooth varieties defined over an algebraically closed field $k$.  Then for every $i \geq 0$ the natural maps $f^* \colon H^i_c(Y_{\et},\QQ_{\ell}) \to H^i_c(X_{\et},\QQ_{\ell})$ and $f^* \colon H^i(Y_{\et},\QQ_{\ell}) \to H^i(X_{\et},\QQ_{\ell})$ are injective.
\end{lemma}

\begin{theorem}
\label{thm:Betti_numbers}
Let $X$ be a smooth projective $F$-split threefold with numerically trivial canonical class and $\wh{q}(X) \geq 1$.  Then the Betti numbers satisfy the following inequalities:
\[
b_{1}(X) = b_{5}(X) \leq 6, \quad b_{2}(X) = b_{4}(X) \leq 23, \quad b_3(X) \leq 44.
\]
\end{theorem}

\begin{proof}
In the case of $\wh{q}(X)=1$, by \autoref{lemma:decomp_rel_dim_2} there exists a generically finite morphism $E \times Z \to X$, where $E$ is an elliptic curve and $Z$ is a K3 surface. In the case of $\wh{q}(X) \geq 2$, by \autoref{thm:beauville_bogomolov_general_logcy}, there exists a finite morphism $A \to X$ where $A$ is an abelian $3$-fold. By \autoref{lem:generically_finite_etale_cohomology_inj} this implies that $b_i(X) \leq \max\{b_i(E \times Z), b_i (A)\}$.  The claim now follows from the K\"unneth formula:
\begin{gather*}
b_0(E) = b_2(E) = 1, \quad b_1(E) = 2, \\
b_0(Z) = b_4(Z) = 1, \quad b_1(Z) = b_3(Z) = 0, \quad b_2(Z) = 22,
\end{gather*}
and the general formula for the Betti number of abelian threefolds:
\begin{equation*}
b_1(A) =6, \quad b_2(A) =15, \quad b_3 = 20.
\end{equation*}
\end{proof}

\begin{corollary}
\label{cor:rational_points}
If $X$ is a smooth projective globally $F$-split and geometrically irreducible threefold with numerically trivial canonical class and $\wh{q}(X) \geq 1$ over a finite field $\bF_q$, then we have
\begin{equation}
\label{eq:rational_points:absolute}
X(\bF_q) \leq 1 + 6 q^{1/2}  + 23 q + 44 q^{3/2} + 23 q^2 + 6 q^{5/2} + q^3.
\end{equation}
and
\begin{equation}
\label{eq:rational_points:difference}
|X(\bF_q) - 1 - q^3| \leq 6 q^{1/2}  + 23 q + 44 q^{3/2} + 23 q^2 + 6 q^{5/2}.
\end{equation}
In particular, 
\begin{equation}
\label{eq:rational_points:non_empty}
X(\bF_q) \neq \emptyset \textrm{ for } q \geq 83.
\end{equation}
\end{corollary}

\begin{proof}
Items \autoref{eq:rational_points:absolute} and \autoref{eq:rational_points:difference} are direct consequence of \autoref{thm:Betti_numbers}, and of the Weil conjectures. Then  \autoref{eq:rational_points:non_empty} follows from \autoref{eq:rational_points:difference} by directly computing that for every integer $q \geq 83$,
\begin{equation*} 
1 + q^3 > 6 q^{1/2}  + 23 q + 44 q^{3/2} + 23 q^2 + 6 q^{5/2}.
\end{equation*}
\end{proof}

\section{Examples}
\label{sec:examples}

\emph{For simplicity in this section the base field $k$ is algebraically closed, and as usual of characteristic $p>0$.}

In \autoref{rem:differences_to_char_zero}, we explained that the two main differences between our theorems and the characteristic zero statements are the following relaxations:
\begin{enumerate}
\item \label{itm:inseparability} In our case, the cover $Z \to X$, for which $Z$ splits as a product, can have inseparability.
\item \label{itm:singularity} The stably Albanese trivial factor of $Z$ might be singular even if we start with a smooth variety.
\end{enumerate}

When the first version of the present paper was made available we only knew singular, however strongly $F$-regular, examples of $K$-trivial varieties really requiring presence of the above.  In this section we describe a simple construction which combined with an example provided by Lam \cite[Thm 2.11]{Lam_On_a_question_of_Patakfalvi_Zdanowicz}, in response to one of the questions we raised, leads to a smooth example of dimension four.  Below, we also provide our older examples based on the work of Matsumoto \cite{Matsumoto_mu_p_and_alpha_p_actions_on_K3_surfaces}.  Those, though singular, have an advantage of being of dimension three.

\begin{construction}\label{const:quotient_of_product}
\label{constr:diag_quotient}
Let $S$ be a projective variety such that 
\begin{enumerate}
\item \label{itm:diag_quotient:K_trivial} $K_S \sim 0$,
\item \label{itm:diag_quotient:GFP} $S$ is globally $F$-split, 
\item \label{itm:diag_quotient:stab_Alb_dim} $\hat q (S)=0$, and
\item \label{itm:diag_quotient:inf_auto} $\mu_p \subseteq \Aut(S)$. 
\end{enumerate}
Consider also the condition:
\begin{enumerate}[resume]
\item \label{itm:diag_quotient:smooth} the quotient $S/\mu_p$ is smooth.
\end{enumerate}
Let $E$ be an ordinary elliptic curve, and let $\mu_p \subseteq E$ be the kernel of the Frobenius endomorphism of $E$. We set then $Y:= S \times E/\mu_p$, where $\mu_p$ acts diagonally. We note that the fibers of $f \colon Y \to E':=E/\mu_p$ are all isomorphic to $S$. In particular, the natural morphism $S \times E  \to Y \times_{E'} E$ is a morphism of fibrations over $E$ mapping the fibers (which are isomorphic to $S$) isomorphically. In particular, this morphism is an isomorphism. 

\begin{proposition}
Let us consider classes of singularities which are inherited by the total space of a fibration if both the base and all the fibres belong to the class.  Then the variety $Y$ inherits the milder class of singularities among those appearing on $S$ or $S/\mu_p$.
\end{proposition}
\begin{proof}
Since $f \colon Y \to E'$ has fibers isomorphic to $S$ it follows that $Y$ has at most as bad singularities as $S$.  Symmetrically, under condition \autoref{itm:diag_quotient:smooth}, the map $g \colon Y \to S/\mu_p$ has fibers isomorphic to $E$ and therefore $Y$ inherits singularities of $S/\mu_p$.  In particular, if $S/\mu_p$ is smooth then $Y$ is smooth too.
\end{proof}

 As $\mu_p$-quotients preserve global $F$-splitting \cite[Thm 4.1 \& Rem 4.9]{Achinger_Ilten_Suss_F-split_and_F-regular_varieties_with_a_diagonalizable_group_action}, $Y$ is globally $F$-split. We claim that since $\mu_p$ acts freely on $E$ it also does on $S \times E$.  For this purpose, we recall that an action of $G \times X \to X$ of a group scheme $G$ on a scheme $X$ is free if and only if the corresponding set theoretic action $G(T) \times X(T) \to X(T)$ is free for every test scheme $T$.  This is preserved under taking products with any $G$-scheme.  Using \cite[Thm A]{Carvajal-Rojas_Finite_torsors_over_strongly_F-regular_singularities} we now see that $K_Y \sim 0$ . Additionally, $f: Y \to E'$ is the Albanese morphism by using that the fibers of $f$ are isomorphic to $S$ in conjunction with  assumption \autoref{itm:diag_quotient:stab_Alb_dim} above.

Let $L$ be an ample Cartier divisor on $S$ and let $p_S : S \times E \to S$ be the first projection. Then $p_S^* (pL)$ descends to a Cartier divisor $M$ on $Y$. We set $I := \Isom_{E'}\Big(\big(S \times E', p_S^*(pL)\big), \big(Y, M\big)\Big)$. We note that 
\begin{equation}
\label{eq:diag_quotient:isom}
I = \Isom (S, pL) \times E/ \mu_p, 
\end{equation} 
where the action is again diagonal, and via post-composition on the left factor. 
\end{construction}

\begin{proposition}
\label{prop:diag_action_tivializing}
In the situation of \autoref{constr:diag_quotient}, if $\Isom(S, pL)$ is a finite group scheme such that $\mu_p \subseteq \Isom(S, pL)$ is a closed subgroup scheme, then as schemes over $E'$ we have an isomorphism $I_{\red}$ is a disjoint union of copies of $E$ with the structure map  on each copy  being the quotient map of the construction of $E'$.  

In particular, a base-change $T \to E'$ trivializes $Y \to E'$ if and only if it factors through $E$. 
\end{proposition}

\begin{proof}
For ease of notation, set $H:= \Isom(S, pL)$, and let $G$ be the identity component of $H$.
Note that as schemes over $E$ we have a $\mu_p$-equivariant isomorphism $H \times E \cong_E \bigcup G \times E$, such that the $\mu_p$-action is diagonal on each of the $G \times E$ factors on the right. Hence, by \autoref{eq:diag_quotient:isom}, it is enough to show that $\left(G \times E/\mu_p \right)_{\red} \cong_{E'} E$. As $\mu_p\hookrightarrow G$ is a closed embedding, so is $\mu_p \times E/\mu_p  \to G \times E/\mu_p $. Hence, it is enough to show that $\mu_p \times E/\mu_p  \cong E$. However, this is immediate by identifying this quotient with the surjective group homomorphism $\mu_p \times E \to E$ given by $(x,y) \mapsto x - y$. 
\end{proof}

We may sum up the above consideration by:

\begin{proposition}\label{prop:contstruction_gives_example}
If a singular variety $S$ satisfies the conditions of \autoref{const:quotient_of_product} then the resulting variety $Y$ inherits the class of singularities of $S$ (or is smooth if condition \autoref{itm:diag_quotient:smooth} is satisfied) and is also $F$-split $K$-trivial and its Albanese morphism, whose fibers are singular, becomes a product only after an inseparable base change.
\end{proposition}

\subsection{Smooth examples}

We now describe a construction by Lam (see \cite[Section 2.4 \& Appendix A]{Lam_On_a_question_of_Patakfalvi_Zdanowicz}) of a threefold $S$ which satisfies all six conditions in \autoref{const:quotient_of_product} and hence leads to a smooth example of a $K$-trivial variety requiring all of our relaxations in its weak Bogomolov--Beauville decomposition.  

\begin{example}
\label{ex:Lam}
We set $p = 5$ and consider the action of $\mu_5$ on $\PP^4_k$ given by all different characters of $\mu_5$ (note that the fixed locus is a set of isolated points):
\[
\mu_5 \times \PP^4_k \to \PP^4_k \quad \xi \cdot (x_0,x_1,x_2,x_3,x_4) = (x_0,\xi x_1,\xi^2 x_2, \xi^3 x_3, \xi^4 x_4).
\]
We also consider the map $\PP^4_k \to \PP^{25}_k$ given by the linear system of all the invariant quintics. The map factors via the quotient map $\pi \colon \PP^4_k \to \PP^4_k/\mu_5$ and it turns out (see \cite[Appendix A]{Lam_On_a_question_of_Patakfalvi_Zdanowicz}) that the induced morphism $\PP^4_k/\mu_5 \to \PP^{25}_k$ is an embedding (here we use the assumption that $p=5$).  We consider the variety $S' \subset \PP^4_k/\mu_5$ given as a general hyperplane section under this embedding.  

\begin{proposition}
The variety $S = \pi^{-1}(S')$ satisfies the conditions from \autoref{const:quotient_of_product} necessary for application of \autoref{prop:contstruction_gives_example}. 
\end{proposition}

\begin{proof}
Clearly, by the construction $S$ admits a $\mu_5$-action and hence \autoref{itm:diag_quotient:inf_auto} holds.  Moreover, we observe that $S' = S/\mu_5$ is smooth, and hence \autoref{itm:diag_quotient:smooth} is satisfied, because a general hyperplane section does not pass through finitely many singular points of $\PP^4_k/\mu_5$.   To see that \autoref{itm:diag_quotient:GFP} is true, that is $S$ is $F$-split, we remark that $F$-splitting is an open condition in a linear system and $S$ is in fact a general $\mu_5$-invariant quintic which degenerates to the sum of hyperplanes $\{x_0x_1x_2x_3x_4 = 0\}$.  This is however $F$-split by a direct computation.  

In order to see that $K_S \sim 0$ we can use adjunction once we know that $S$ is in fact normal. However, this follows from \autoref{prop:fibres_are_regular_codimension_one} since $S$ is the general fibre of a morphism $Y \to E'$ from an $F$-split variety as in the description of the construction above.  This is admittedly a slightly peculiar way of proving that $S$ is not too singular, but we could not easily utilize the fact that $S$ is a total space of $\mu_5$-torsor over a smooth variety $S'$.

We also need to prove that $\hat q(S)=0$.  By \autoref{lem:invariance_augmented_irregularity} it suffices to show that $\hat q(S') = 0$.  By the singular version of the Lefschetz hyperplane theorem $\piet(S) = \{1\}$ (see \stacksproj{0ELE}) and hence using invariance of the \'etale site $\piet(S') = \{1\}$.  Consequently, we need to show that $\dim \Alb S' = 0$, since $S'$ is smooth and therefore every quasi-\'etale cover is \'etale.  However, if the fundamental group is trivial then the Albanese is zero dimensional, because the prime to $p$ abelian part of the fundamental group of $X$ is isomorphic to prime to $p$ torsion on the Albanese.

Finally, to see that $S$ is singular we reason by contradiction. We first observe that it admits a $\mu_5$-action and hence $ H^3(S,\Omega^1_S) = H^0(S,T_S) \neq 0$ where the equality follows from Serre duality and the condition $\omega_S \simeq \cO_S$ which is clearly true for a smooth quintic.  This yields a contradiction because Hodge numbers are invariant (see \cite[Théorème 1.5]{Seminaire_de_Geometrie_Algebrique_du_Bois_Marie_7_2_Groupes_de_monodromie_en_geometrie_algebrique}) under liftings for hypersurfaces and $h^{1,3} = 0$ for quintic threefolds in characteristic zero.  
\end{proof}

\end{example}

\subsection{Singular examples of smallest possible dimension}

\begin{example}
\label{ex:inseparable_base_change_necessary}
Consider the example constructed by Matsumoto \cite[Example 10.5]{Matsumoto_mu_p_and_alpha_p_actions_on_K3_surfaces}, that is a surface:
\[
S := \left\{w^2 + x^5 y + y^5 z + z^5 x + a x^2y^2 z^2=0\right\} \subseteq P:=\bP_k(3,1,1,1)_{w,x,y,z},
\] 
where $k$ is an algebraically closed field of characteristic $p=7$ and $a^3 \neq 1$. This surface $S$ is a K3 surface with only Du-Val singularities and a $\mu_p$-action. In particular, all the assumptions of \autoref{constr:diag_quotient} are automatically satisfied \cite{Hara_Classification_of_two_dimensional_F_regular_and_F_pure_singularities}, except possibly that $S$ is globally $F$-split. 

\emph{We claim that in fact $S$ is also globally $F$-split.} Let $R$ be the graded ring $k[w,x,y,z]$, where $\deg w= 3$, and $\deg x= \deg y = \deg z = 1$, and let $k[w,z,y,z]_d$ denote the homogeneous part of degree $d$. Then, $H^0\left(P, \omega_P^{1-p}\right) \to H^0(P, \sO_P)=k$ can be identified with the $k$-linear homomorphism $k[w,x,y,z]_{6(p-1)} \to k$ sending $(wxyz)^{p-1}$ to $1$ and all the other monomials to $0$. Consider now the following diagram:
\begin{equation*}
\xymatrix{
k\cong H^0\left(P, \omega_P(S)^{1-p} \right)  \ar[d]^{\cong} \ar[r] &  H^0\left(P, \omega_P^{1-p}\right) \ar[r] &  H^0(P, \sO_P) \cong k \ar[d]^{\cong} \\
k \cong H^0\left(S, \omega_S^{1-p}\right) \ar[rr] & & H^0(S, \sO_S) \cong k
}
\end{equation*}
The diagram shows that the one can show that $X$ is globally $F$-split by computing that the composition of the top horizontal maps is bijective. Hence, one needs to compute whether the coefficient of $(wxyz)^{p-1} = (wxyz)^{6}$ in  
\begin{equation*}
\left(w^2 + x^5 y + y^5 z + z^5 x + a x^2y^2 z^2\right)^{p-1} = \left(w^2 + x^5 y + y^5 z + z^5 x + a x^2y^2 z^2\right)^{6}
\end{equation*}
 is non-zero. However, this coefficient is \Big(accounting for $\left(w^2\right)^3 \cdot \left(x^5 y\right) \cdot \left(y^5 z\right) \cdot \left(z^5x \right)$ and for $\left(w^2\right)^3 \cdot \left( a x^2 y^2 z^2\right)^3$\Big):
\begin{equation*}
\frac{6!}{3!} + a^3\frac{6!}{3! 3!} = 120 + a^3 20 
\explshift{-30pt}{=}{$\characteristic k = 7$}
1 - a^3 
\explshift{60pt}{\neq}{$a^3 \neq 1$ by assumption} 0 
\end{equation*}
This concludes our claim above. 

Hence, all conditions of \autoref{constr:diag_quotient} are satisfied, and according to \autoref{prop:diag_action_tivializing} $Y \to E'$ constructed using $S$ of the present example yields an example where in \autoref{thm:beauville_bogomolov_general_logcy} we must have an inseparable base-change. 
\end{example}

The topic of the remainder of \autoref{sec:examples} is to prove that in dimension two such examples cannot be non-singular.  We learned the argument from Yuya Matsumoto whose paper \cite{Matsumoto_mu_p_and_alpha_p_actions_on_K3_surfaces} gives a detailed description of $\mu_p$ and $\alpha_p$ quotients of K3 surfaces with canonical singularities.  We start with a lemma describing the behaviour of Betti numbers under resolutions.

\begin{lemma}
\label{lem:betti_numbers_resolution}
Let $f \colon Y \to X$ be a minimal log resolution of a proper canonical surface defined over an algebraically closed field.  Then $b_1(Y) = b_1(X)$ and $b_2(Y) = b_2(X) + \sum_{x \in X_{\rm sing}} d_x$, where $d_x$ is the number of irreducible components in the minimal resolution of the germ of $X$ at $x \in X_{\rm sing}$. 
\end{lemma}

\begin{proposition}[Matsumoto, private communication]
Let $X \to Y$ be a generically finite map of normal Gorenstein surfaces with trivial canonical class and at most canonical singularities defined over an algebraically closed field of characteristic $p>0$.  Assume that $Y$ is smooth and weakly ordinary (equiv.  globally $F$-split).  Then $X$ is smooth as well.
\end{proposition}

\begin{proof}
Assume by contradiction that $X$ is not smooth.  Let $X' \to X$ be the minimal resolution of singularities.  Since $X$ was canonical the resolution is crepant and therefore $X'$ is a smooth Calabi--Yau variety.  The map $X' \to Y$ is generically finite, and therefore by \autoref{lem:generically_finite_etale_cohomology_inj} and \autoref{lem:betti_numbers_resolution} the Betti numbers  satisfy the inequalities $b_1(X') \geq b_1(Y)$, $b_2(X') > b_2(Y)$.  Both surfaces $X'$ and $Y$ are Calabi--Yau and therefore by the classification of surfaces (see \autoref{prop:classification_ordinary_cy}) we obtain a table where the entries contain the contradictory inequalities of Betti numbers for the associated possibilities for types of $X'$ and $Y$ (respectively in rows and columns).

{\scriptsize
\begin{center}
\begin{tabular}{ |c|c|c|c|c| }
\cline{1-5}
& K3  & abelian & Enriques & hyperelliptic \\
\cline{1-5}
K3 & $22 = b_2(X') \not> b_2(Y) = 22$ & $0 = b_1(X') \not> b_1(Y) = 4$ & No contradiction & $0 = b_1(X') \not> b_1(Y) = 2$ \\
\cline{1-5}
abelian & $2 = b_2(X') \not> b_2(Y) = 22$ & $6 = b_2(X') \not> b_2(Y) = 6$ & $2 = b_2(X') \not> b_2(Y) = 10$ &  No contradiction \\
\cline{1-5}
Enriques & $10 = b_2(X') \not> b_2(Y) = 22$ & $0 = b_1(X') \not> b_1(Y) = 4$ & $10 = b_2(X') \not> b_2(Y) = 10$ & $0 = b_1(X') \not> b_1(Y) = 2$ \\
\cline{1-5}
hyperelliptic & $2 = b_2(X') \not> b_2(Y) = 22$ & $2 = b_1(X') \not> b_1(Y) = 4$ & $2 = b_2(X') \not> b_2(Y) = 10$ & $2 = b_2(X') \not> b_2(Y) = 2$ \\
\cline{1-5}
\cline{1-5}
\end{tabular}
\end{center}
}

This leaves us with the following cases:

\begin{enumerate}
\item either $X'$ is abelian and $Y$ is hyperelliptic, which gives a contradiction because $X$ was singular canonical, and hence $X'$ contains a rational curve, or
\item $X'$ is a K3 surface and $Y$ is an Enriques surface.
\end{enumerate}
In the latter case, we first observe that if $p > 2$ we may consider the cartesian diagram
\[
\xymatrix{
    X' \times_Y Y' \ar[r]\ar[d]_{\text{\'et}} & Y' \ar[d]^{\text{\'et}} \\
    X' \ar@{-->}[ru]\ar[r] & Y.
}
\]
where $Y' \to Y$ is the canonical \'etale K3 cover of $Y$.  Since $X'$ is a K3 surface, it is simply connected, and hence the map $X' \times_Y Y' \to X'$ admits a section yielding a generically finite morphism $X' \to Y'$.  The map $X' \to Y'$ contracts a rational curve and therefore $b_2(X') > b_2(Y')$ which is a contradiction, because both $X'$ and $Y'$ are K3 surfaces.  If $p = 2$, we use the weakly ordinary assumption to see that $Y$ is in fact a singular Enriques surface (see \autoref{prop:classification_ordinary_cy} for the definition of singular Enriques surfaces), and hence we may repeat the above argument because a required $Y'$ also exists.
\end{proof}

\begin{remark}
Without the weak ordinarity assumption the last part of the above argument breaks down.  There are examples Enriques surfaces in characteristic two (non $F$-split and called supersingular) which are covered by singular varieties with trivial canonical class.
\end{remark}


\section{Counterexample to the decomposition in the non weakly ordinary case}\label{sec:moret_bailly}

In this section we recall the examples of Calabi--Yau varieties arising from Moret-Bailly's pencils of abelian varieties, and utilize them to show that the weakly ordinary/globally $F$-split assumption cannot be dropped from our Beauville-Bogomolov type decomposition statements.  We refer to \cite[Section 3]{Schroer_Some_CY_threefolds_with_obstructed_deformations} for a detailed exposition of three dimensional version of the the construction.  In what follows we provide a slightly more general family of examples which in particular work in arbitrarily large characteristic, although we need to increase the dimension.

\begin{example}
\label{ex:Moret_Bailly}
Let us assume that $k$ is algebraically closed and let $A = E \times \ldots \times E$ be a product of $n+1$ supersingular elliptic curves.  Observe that $E$ contains a subgroup scheme isomorphic to $\alpha_p$, and hence $A$ admits a subgroup scheme isomorphic to $\alpha_p \times \ldots \times \alpha_p$, which in turn contains a family of subgroup schemes isomorphic to $\alpha_p$ parametrized by $\bP^n$.  More precisely, every line $\ell \in \PP({\rm Lie}_{\alpha_p \times \ldots \times \alpha_p})$ induces a sub-Lie algebra ${\rm Lie}_\ell \subset {\rm Lie}_{\alpha_p \times \ldots \times \alpha_p}$ which in turn gives a subgroup scheme $G_\ell$ isomorphic to $\alpha_p$.  In fact, all such group schemes can be combined into a tautological family $G_{\PP^{n}} \subset A \times \PP^{n}$ whose relative Lie algebra is isomorphic to the pullback $\cO_{A \times \PP^{n}}(-1)$ of $\cO_{\PP^{n}}(-1)$ under the projection $A \times \PP^{n} \to \PP^{n}$.  Taking a base change of the family under any degree $r$ self-map of $\PP^{n}$ we obtain a relative group scheme $G_\ast$ with Lie algebra $\cO_{A \times \PP^{n}}(-r)$.  

We set $X = A \times \PP^{n} / G_{\ast}$, and $\pi \colon A \times \PP^{n} \to X$ the natural projection.  As $G_{\ast}$ acts freely on $A \times \PP^{n}$, $X$ is smooth.  In order to compute the canonical divisor of $X$ we use the theory of foliations (see \cite{Ekedahl_Foliations_and_inseparable_morphisms}).  The Lie algebra of $G_{\ast}$ yields a foliations $\cO_{A \times \PP^{n}}(-r) \subset \cT_{A \times \PP^{n}}$ which in turn implies using \cite[Proposition 1.7]{Ekedahl_Foliations_and_inseparable_morphisms} that
\[
\pi^*\omega_X = \cO_{A \times \PP^{n}}\left(r(p-1) - n -1 \right).
\]
Since $\pi$ is a universal homeomorphism, factorizing Frobenius morphism of $X$, this formula implies that canonical divisor of $X$ is numerically trivial if and only if 
\[
r(p-1) = n+1.
\]
The dimension of the variety $X$ is then equal to $\dim X = 2n+1$, which agrees with $2p-3$ if $r=1$.  However, we obtain many other examples for $r \neq 1$. For example, the list of the low dimensional examples is the following:
\begin{itemize}
    \item $\dim X = 3$: $(n,p,r) = (1,2,2)$ or $(n,p,r) = (1,3,1)$
    \item $\dim X = 5$: $(n,p,r) = (2,2,3)$
    \item $\dim X = 7$: $(n,p,r) = (3,2,4)$ or $(n,p,r) = (3,3,2)$, or $(n,p,r) = (3,5,1)$
\end{itemize}

We set $E' = E/\alpha_p$.  We first note that since $X$ maps to a supersingular abelian variety $A':=E' \times \ldots \times E'$ its reduced Picard scheme is supersingular and hence admits no non-trivial $p$-torsion points.  This implies that $\omega_X$ is not only numerically trivial but actually trivial.
\end{example}

We note that a construction similar to \autoref{ex:Moret_Bailly} was used in  \cite{Roessler_Schroer_Moret-Bailly_families_and_non-liftable_schemes} to construct non-liftable Calabi-Yau varieties in arbitrary high characteristics.

\begin{proposition}
\label{prop:Moret_Bailly_Albanese}
Let $X$ be one of the $K$-trivial varieties constructed in \autoref{ex:Moret_Bailly}.  Then, the Albanese morphism of $X$ is equal to the natural projection $X \to A'$.
\end{proposition}

\begin{proof}
We consider the following diagram, where the vertical arrows are Albanese morphisms:
\begin{equation}
\label{eq:Moret_Bailly_Albanese}
\xymatrix{
E \times \ldots \times E \times \PP^{n} \ar[r]^(0.7){\tau} \ar[d] & X \ar[r]\ar[d] & E' \times \ldots \times E' \times \PP^{n} \ar[d] \\
A=E \times \ldots \times E \ar[r]^(0.6){\psi} & \Alb X \ar[r] & E' \times \ldots \times E' = A',
}
\end{equation}
As $\tau$ is inseparable, the images of the $\bP^{n}$ fibers via $\tau$ are rationally chain connected, and hence, they must be contracted by $X \to \Alb X$. It follows then that all the arrows in diagram \autoref{eq:Moret_Bailly_Albanese} are surjective morphism of varieties.  In particular, the left bottom arrow $\psi$ is an isogeny whose kernel satisfies the condition $\Ker(\psi) \subset \alpha_p \times \ldots \times \alpha_p$.  Restricting the composition $E \times \ldots \times E \times \PP^{n} \to \Alb X$, to the fibers over different  closed points of $\PP^{n}$ we observe that $\Ker(\psi)$ contains $G_\ell$, for every $\ell \in \PP^{n}$, and so we claim that $\Ker(\psi) = \alpha_p \times \ldots \times \alpha_p$.  Indeed, $\Ker(\psi)$ is a height one group scheme, consequently it is determined by its Lie algebra, which clearly contains the Lie algebra of $\alpha_p \times \ldots \times \alpha_p$.
\end{proof}

\autoref{prop:Moret_Bailly_Albanese} tells us already  that bad things happen in the non weakly ordinary case. For example, \cite[Thm 1.1]{Ejiri_When_is_the_Albanese_morphism_an_algebraic_fiber_space_in_positive___characteristic?}
fails in the non-weakly ordinary case, which  in the weakly ordinary and $K$-trivial case tells us that the general fiber of the Albanese morphism is $F$-pure and non-reduced. Indeed, in the case of \autoref{prop:Moret_Bailly_Albanese}, the Albanese morphism has non-reduced fibers the reduced subscheme of which is $\bP^n$.

The following proposition takes this further, by showing that  the weakly ordinary/globally $F$-split assumption cannot be dropped from our Beauville--Bogomolov type decomposition statements. 

\begin{proposition}
\label{prop:bb_decomp_counterexample_moret_bailly}
Let $X$ be one of the $K$-trivial varieties constructed in \autoref{ex:Moret_Bailly}, and let $Y \xrightarrow{\rho} Z \xrightarrow{\gamma} X$ be a composition of finite morphisms such that
\begin{enumerate}
\item $Y$ is integral, 
\item $\gamma$ is \'etale, and 
\item $\rho$ is a torsor under an infinitesimal group scheme over $k$.
\end{enumerate}
Then, the Albanese morphism $Y \to \Alb Y$ is not split. 
\end{proposition}

\begin{proof}
Assume that $\alb_Y \colon Y \to \Alb(Y)$ is split. That is, $Y = \Alb(Y) \times V$ for some integral projective variety $V$. In particular, we have $V \cong \alb_Y^{-1}(0)$. 

First note that it is enough to show that $V \cong \bP^{n}$. Indeed, as torsors are relatively $K$-trivial, the following line yields a contradiction
\begin{equation*}
0 \sim \tau^* K_X \sim K_Y 
\expl{\sim}{$Y \cong \Alb(Y) \times V$} 
\pr_{\Alb_Y}^* K_{\Alb Y}  + \pr_{V}^* K_V
\expl{\not\sim}{$V \cong \bP^{n}$}
 0. 
\end{equation*}
So, \emph{in the remaining part of the proof we show that $V \cong \bP^{n}$}.

As $ \gamma \times_X \left(A \times \bP^{n}\right) $ is \'etale, and as $\bP^{n}$ is \'etale simply connected, this morphism can be identified with $D \times \bP^{n} \to A \times \bP^{n}$ induced by an \'etale isogeny $D \to A$. 

Now, let us analyze the morphism $\rho \times_ Z \left( D \times \bP^{n}\right) $. This is a torsor under an infinitesimal group scheme. Since the fundamental group scheme functor is trivial for the projective space \cite[Corollary, p. 93]{Nori_Fundamental_Group_Scheme} and commutes with products of proper varieties \cite[Proposition 2.1]{Mehta_Subramanian_On_the_fundamental_group_scheme}, we obtain that this morphism can be identified with $W \times \bP^{n} \to D \times \bP^{n}$, induced by an infinitesimal torsor $W \to D$. According to  \autoref{prop:torsors_on_abelian_vars}, $W_{\red} \to D$ is an isogeny of abelian varieties. Set $B:= W_{\red}$. 

Now, consider the following commutative diagram, consisting of three parallel commutative squares. The right one, drawn with dotted arrows, consists of the Albanese varieties of the varieties in the middle square. The left square, drawn with dashed arrows, then consists of the reduced subschemes of the  fibers over $0$ of the corresponding Albanese morphisms. Note that as three of these fibers are reduced either by construction or by the splitting assumption on $Y \to \Alb Y$, we have to actually take the reduced subscheme only for $\alb_X^{-1}(0)$.
\begin{equation*}
\xymatrix@C=30pt{
& & \ar@{-->}[dll]_{\id_{\bP^{n}}} \ar@{-->}[dd]|{\parbox{10pt}{$\hole$  \\ $\hole$ }} \bP^{n} \cong \{0\} \times \bP^{n} \ar@{^(->}[r] 
%
%
& B \times \bP^{n} \ar[r] \ar[dd] \ar[dll]
%
%
& B \ar@{..>}[dd] \ar@{..>}[dll] \\
%
%
\bP^{n} \cong \{0\} \times \bP^{n} \ar@{^(->}[r] \ar@{-->}@/_1pc/[dd] &
%
%
A \times \bP^{n} \ar@/_1pc/[dd] \ar[r] & 
%
%
A \ar|{\parbox{10pt}{$\hole$  \\ $\hole$ }}@/_1pc/@{..>}[dd] \\ 
%
%
& & V = \alb_Y^{-1}(0)  \ar@{-->}[dll] \ar@{^(->}[r]
%
%
& Y \ar[r]^{\alb_Y} \ar[dll]
%
%
&  \Alb(Y) \ar@{..>}[dll] \\
%
%
\alb_X^{-1}(0)_{\red} \ar@{^(->}[r]
%
%
& X \ar[d] \ar[r]
%
%
& 
A' \expl{=}{ \autoref{prop:Moret_Bailly_Albanese}} \Alb(X) \\
%
%
& \bP^{n} \ar@/^8pc/[uuul]^{\cong}
}
\end{equation*}
Let $\pr_{\bP^{n}}$ be the projection $A \times \bP^{n} \to \bP^{n}$. We see from the diagram that $\alb_X^{-1}(0)_{\red}$ factors the morphism $\left. \pr_{ \bP^{n}}\right|_{0 \times \bP^{n}}: \{0\} \times \bP^{n}  \to \bP^{n}$, and it is irreducible of dimension $n$. Hence, we have $\bP^{n} \cong \alb_X^{-1}(0)_{\red}$, and the morphism $\{0\} \times \bP^{n} \to  \alb_X^{-1}(0)_{\red}$ identifies with $\id_{\bP^{n}}$. But, then, $V$ also factors the identity morphism of $\bP^{n}$, which implies that $V \cong \bP^{n}$.

\end{proof}

\appendix

\section{Isotriviality of the Albanese morphism of klt varieties with nef anticanonical divisor in characteristic zero -- jointly with Giulio Codogni}
\label{sec:appendix}

The aim of this appendix is to generalize the result of \cite[Theorem 1.2]{Cao_Albanese_maps_of_projective_manifolds_with_nef_anticanonical_bundles} 
\cite[Proposition 2.8]{Cao_Albanese_maps_of_projective_manifolds_with_nef_anticanonical_bundles} concerning the structure of the Albanese morphism of varieties with nef anticanonical class to the boundary and singular setting.   The boundary free problems were already posed in \cite[Section 4]{Demailly_Peternell_Schneider_Compact_Kahler_Manifolds_with_Hermitian_semipositive}.  Additionally, the $K$-trivial version of our results was shown in \cite{Xu_Homogeneous_fibrations_on_log_Calabi_Yau_varieties}.

We note that  after the present paper was made public, th main statements of this section were reproved with other methods in \cite{Wang_Structure_of_projective_varieties_with_nef_anticanonical_divisor_the_case_of_log_terminal_singularities},
and they were also generalized to similar statement for the MRC fibration. 

\emph{In this section we work over the field of complex numbers denoted by $\CC$.}

\begin{definition}
\label{defin:isotrivial_etale}
Let $f \colon X \to T$ be a flat projective morphism of varieties defined over $\CC$.  We say that $f$ is \emph{isotrivial} if it is locally trivial (isomorphic to the product family) in the \'etale topology.  We say that a morphism $(X,\Delta) \to T$ of a pair is \emph{isotrivial} if it is locally trivial in the \'etale topology as a morphism of pairs.
\end{definition}

We first recall the necessary isotriviality statement in the analytic topology.  We say that a flat fibration $f \colon X \to T$ of projective varieties is \emph{analytically isotrivial relative to a divisor $L$ on $X$} if there exists a trivializing cover $\pi \colon U \to X$ in the analytic topology such that the pair $(X_U,\pi^*L)$ is isomorphic to $(X_t \times U,L_t \times U)$ for a closed point $t \in T$.  The definition clearly extends to the setting with boundary.

We shall need some results concerning numerically flat bundles in characteristic zero.  We recall that by definition, a vector bundle $\cE$ on a smooth projective variety $X$ is \emph{numerically flat} if both $\cE$ and $\cE^\vee$ are nef vector bundles.  In characteristic zero this condition turns out to be equivalent to $H$-semistability and the vanishing 
\[
c_1(\cE).H^{n-1} = c_2(\cE).H^{n-2} = 0,
\] for every (equiv. any) ample line bundle $H$ on $X$. Additionally, in characteristic zero, it is also equivalent to the condition that $\cE$ comes from an extension of irreducible unitary flat bundles (equiv.~irreducible unitary representations of $\pi_1(X)$).  In the proof of the following proposition, we summarize the known results and justify some of the above claims.  We refer to \cite[\S 1.3 and \S 5]{Langer_On_the_S-fundamental_group_scheme_I} for a detailed discussion.

\begin{proposition}[Numerically flat bundles and local systems]
\label{prop:numerically_trivial_are_local_systems}
Let $T$ be a smooth projective variety.  Then for every numerically flat bundle $\cE$ there exists a canonical choice of a holomorphic integrable connection on $\cE$.  Moreover, every $\cO_X$-linear homomorphism of numerically flat bundles is parallel with respect to the canonical connections.
\end{proposition}

\begin{proof}
First of all, we observe using \cite[Thm 1.18 \& Cor 1.19]{Demailly_Peternell_Schneider_Compact_Complex_manifolds_with_numerically_effective_tangent_bundles} that every numerically flat bundle is an extension of unitary local systems and hence also a semistable bundle with vanishing Chern classes.  The association $\cE \to (\cE,0)$, where $0$ denotes the zero morphism $\cE \to \cE \otimes \Omega^1_T$, yields a fully faithful functor from the category of numerically flat bundles to the category of semistable Higgs sheaves $\cE$ satisfying $c_1(\cE) = c_2(\cE) = 0$.  By \cite[Corollary 3.10]{Simpson_Higgs_bundles_and_local_systems}, the latter category is in fact equivalent to the category of local systems.  By the remark at the end of subsection ''Examples`` of  \cite[Section 3]{Simpson_Higgs_bundles_and_local_systems}, the equivalence does not change the holomorphic structure of the bundles in our case which consequently gives the claim of the proposition.
\end{proof}

\begin{proposition}[{\cite[Proposition 2.8]{Cao_Albanese_maps_of_projective_manifolds_with_nef_anticanonical_bundles}}]
\label{prop:isotriviality_criterion_cao}
Let $f \colon X \to T$ be a flat projective fibration of a normal variety $X$ onto a smooth projective variety $T$, and let $L$ be an $f$-ample line bundle.  Assume that $f_*\cO_X(nL)$ is a numerically flat vector bundle for every integer $n \geq 1$.  Then, the morphism $f$ is analytically isotrivial relative to some multiple of the divisor $L$.  

Assume further that $X$ is $\QQ$-factorial and equipped with a boundary divisor $\Delta$ whose components do not contain any fibers of $f$ and such that $(X,\Delta)$ is klt.  Then the morphism $f \colon (X,\Delta) \to T$ is analytically isotrivial relative to some multiple of the divisor $L$ if the sheaves $f_*\cO_X(nL-D)$ are numerically flat (in particular locally free), for every irreducible component $D$ of the boundary $\Delta$. 
\end{proposition}

\begin{proof}
Since $L$ is $f$-ample, we observe that for $m>0$ sufficiently large the formation of $f_*\cO_X(nmL)$ commutes with base change and that the family $f \colon X \to T$ is determined by the sheaf of rings $\cR_{X/T} = \bigoplus_{n \in \NN} f_*\cO_X(nmL)$.  Under our assumptions the bundles $f_*\cO_X(nmL)$ forming $\cR_{X/T}$ are numerically flat.  By \autoref{prop:numerically_trivial_are_local_systems} 
this implies that they admit a structure of local systems which in turn means they are trivialized on the universal cover.  Taking a sufficiently large integer $m > 0$ we may assume that the natural multiplication maps
\[
\Sym^n f_*\cO_X(mL) \to f_*\cO_X(nmL) 
\]
are surjective.  Using 
\cite[Theorem 6.2.12 (iii)]{Lazarsfeld_Positivity_in_algebraic_geometry_II}, we see that the symmetric power bundles are also numerically flat, and consequently the multiplication is a surjective morphism of numerically flat bundles.  Consequently, by \autoref{prop:numerically_trivial_are_local_systems} the homomorphism is parallel with respect to canonical connections.  In particular, the kernels of multiplications maps, describing the family entirely, are also local systems and hence are trivialized on the universal cover. Hence, $\sR_{X/T}$ itself is a sheaf of rings trivialized on the universal cover.  This concludes the proof of the first part of the statement. 

For the second part of the statement, our assumptions yield that the sheaf of ideals $\bigoplus_{n \in \bN} f_* \sO_X(nmL - D) =: \sI_{D/X} \subseteq \sR_{X/T}$ is also numerically flat. Hence as above, $\sI_{D/X}$  trivializes on the universal cover with $\sI_{D/X} \hookrightarrow \sR_{X/T}$ being horizontal as well. So, it is enough to show that for every $t \in T$ we have $\sI_{D/X}|_{X_t} \cong \sI_{D_t/X_t}$, as then not only the fibers of $X \to T$ are identified on the universal cover, but also the ideals of the fibers of the components of the boundary.  We defer the proof of this statement to the following \autoref{lem:base_change_ideals} since the result will actually be used in a few other places. 

We also remark that we crucially take advantage of numerical flatness of the bundles in question.  A morphism of arbitrary local systems is not necessarily flat.  For examples on varieties admitting non-trivial global differential forms there are many non-isomorphic flat structures on the trivial bundle, the identity map is the not parallel.
\end{proof}

\begin{lemma}\label{lem:base_change_ideals} 
Let $f \colon X \to T$ be a flat projective fibration of a klt pair $(X,\Delta)$ onto a smooth projective variety $T$, and let $D$ be a $\QQ$-Cartier divisor which does not contain any irreducible components of the fibers of $f$.  Suppose that $L$ is an $f$-ample line bundle.  Then for every regular base change $S \to T$  (here regular means that $S$ is regular) the natural morphism 
\[
\cO_X(nL - D)_{|X_S}  \to \cO_X(nL_S - D_S)
\]
is an isomorphism for every $n \in \ZZ$.  Moreover, for $n \gg 0$ the higher direct images satisfy the vanishing
\[
R^i f_*\cO_X(nL - D) = 0
\]
for $i>0$ and universally, that is, after any base change.  For all such $n$ the sheaves $f_*\cO_X(nL - D)$ are locally free with formation commuting with any regular base change.
\end{lemma}
\begin{proof}
First, we claim that for $n$ sufficiently large we have the appropriate base change property for $\cO_X(nL - D)$.  We recall that $X_S$ denotes the associated base change $X \times_T S$.  Since $(X,\Delta)$ is klt and $D$ is $\QQ$-Cartier, for every $n>0$ the sheaves $\cO_X(nL - D)$ are Cohen--Macaulay (see \cite[Corollary 5.25]{Kollar_Mori_Birational_geometry_of_algebraic_varieties}).  This property is preserved under regular base change, and therefore the natural map 
\begin{equation}
\label{eq:isotriviality_criterion_cao:ideals}
\cO_X(nL - D)_{|X_S} \explpar{250pt}{\isom}{The sheaf $\cO_X(nL - D)$ is Cohen--Macaulay, and hence reflexive of rank one after every base change.  The associated divisor is equal to $nL_S - D_S$ because $D$ does not contain any fibers.} \cO_{X_S}(nL_S - D_S)
\end{equation}
is an isomorphism.  By \cite[Corollary 2.14]{Bhatt_Ho_Patakfalvi_Schnell_Moduli_of_products_of_stable_varieties} the sheaves $\cO_X(nL - D)$ are flat. Consequently using \autoref{eq:isotriviality_criterion_cao:ideals} in conjunction with the $f$-ampleness of $L$ and the cohomology and base change theorem, for $n \gg 0$, we have the necessary vanishing of higher direct images and the identification $f_*\cO_X(nL-D)_{|X_S} \cong f_{S*}\cO_{X_S}(nL_S - D_S)$.
\end{proof}

\begin{proposition}[{Algebraization of \autoref{prop:isotriviality_criterion_cao}}]
\label{prop:isotriviality_criterion_cao_algebraized}
Let $f \colon X \to T$ be a flat projective fibration of a normal variety $X$ onto a smooth projective variety, and let $L$ be an $f$-ample line bundle.  Assume that the bundles $f_*\cO_X(nL)$ are numerically flat, for every $n>0$. Then the morphism $f$ is isotrivial over $T$ relative to some multiple of $L$.  

The same statement holds true for a morphism of pairs $(X,\Delta) \to T$ provided that $X$ is $\QQ$-factorial, no irreducible components of the boundary $\Delta$ contain any fibers of $f$, the pair $(X,\Delta)$ is klt and the sheaves $f_*\cO_X(nL - D)$ are numerically flat vector bundles for every integer $n>0$.
\end{proposition}

\begin{proof}
We begin by using \autoref{prop:isotriviality_criterion_cao} to see that the family is analytically isotrivial relative to to some multiple of $L$.  Using GAGA this means that there exists an integer $m$ such that for every two points $t,t' \in T$ the pairs $(X_t,m L_t)$ and $(X_{t'},m L_{t'})$ are isomorphic.  Consequently, the relative isomorphism scheme $\Isom_T((X,m L),(X_t \times T,m L_t \times T)$ is a torsor over $T$ under a finite type smooth (we work in characteristic zero) group scheme $\Aut(X_t,m L_t)$.  The map $\Isom_T((X,m L),(X_t \times T,m L_t \times T) \to T$ is therefore smooth and hence admits sections \'etale locally.  This finishes the proof of the first part of the proposition.  The same argument works in the boundary setting by substituting the relevant isomorphism scheme with its boundary version.  See \autoref{constr:log_isom} for the construction, and \autoref{prop:base_change_for_isom_scheme} for the appropriate base change properties.
\end{proof}

\begin{remark}
\label{remark:char_zero_multiple_equal_one}
From the proof (see the final remarks in the proof of \autoref{prop:isotriviality_criterion_cao}) we directly see that the multiplying factor $m$ necessary to get analytic isotriviality with respect to line bundle $mL$ is equal to one if the natural multiplication map $\Sym^n \cO_X(L) \to \cO_X(nL)$ is surjective for every $n \geq 1$.
\end{remark}

Since we develop tools to get numerical flatness of relative section bundles only over curves, we need the following lemma for bases of higher dimension.

\begin{lemma}[{Algebraization of \autoref{prop:isotriviality_criterion_cao}}]
\label{lemma:isotriviality_criterion_cao_algebraizedII}
Let $f \colon X \to T$ be a flat projective fibration defined over $\CC$, and let $L$ be an $f$-ample line bundle.  Suppose that direct images of $\cO_X(nL)$  satisfy the following conditions:
\begin{itemize}
    \item $f_*\cO_X(nL)$, are locally free with formation commuting with any base change,
    \item the multiplication map $\Sym^m\cO_X(L) \to \cO_X(mL)$ is surjective, for every $m \geq 1$.
\end{itemize} 
Moreover, assume that for a covering family of smooth projective curves $\{C_s \subset X\}_{s \in S}$ the vector bundle 
\[
f_*\cO_X(nL)_{|C_s} \expl{\isom}{formation of pushforward commutes with base change} f_{C_s*}\cO_{X_{C_s}}(nL_{|C_s})
\] 
is numerically flat, for every $n>0$ and for every $s \in S$.  Then the morphism $f$ is isotrivial over $T$ relative to $L$.  
\end{lemma}

\begin{proof}
As in the proof of \autoref{prop:isotriviality_criterion_cao_algebraized} we consider the relative isomorphism scheme $\Isom_T((X,L),(X_t \times T,L_t \times T))$.  In order to get our claim, it suffices to prove that the natural map from the isomorphism scheme is surjective onto the base $T$.  We denote by $X_{C_s}$ the base change $X \times_T C_s$, and by $f_{C_s}$ the natural projection.  Since the family $\{C_s\}_{s \in S}$ covers $T$, the bundles $f_{C_s*}\cO_{X_{C_s}}(nL_{C_s})$ are isomorphic to $f_*\cO_X(nL)_{|C_s}$ and the surjectivity of the multiplication map which is inherited by the bundles $f_{C_s*}\cO_{X_{C_s}}(nL_{|C_s})$ (by the base change and the assumptions), using base change theorem for the isomorphism scheme (see \autoref{prop:base_change_for_isom_scheme}) we may reduce to the case of $T = C_s$.  This finishes the proof using \autoref{prop:isotriviality_criterion_cao_algebraized} along with \autoref{remark:char_zero_multiple_equal_one}.
\end{proof}

\begin{remark}
\label{remark:numerical_flatness_divisor}
We note that the cohomological criteria necessary in \autoref{lemma:isotriviality_criterion_cao_algebraizedII} can be guaranteed by taking a sufficient power of the line bundle $L$.  Moreover, the statement of \autoref{lemma:isotriviality_criterion_cao_algebraizedII} holds true for a morphism of pairs $(X,\Delta) \to T$ provided that every irreducible component $D$ of the boundary $\Delta$ intersect every fiber of $f$ properly and the formation of $f_*\cO_X(nL - D)$ commutes with the base change $C_s \to T$ yielding numerically flat bundles
\[
f_*\cO_X(nL-D)_{|C_s} \isom f_{C_s*}\cO_{X_{C_s}}\left(nL_{C_s} - D_{C_s}\right)
\] 
for every sufficiently large integer $n>0$, where $f_{C_s} \colon X_{C_s} \to C_s$ is the base change of $f$ along the inclusion $C_s \to T$.  This happens for example if the divisor $D$ is $\QQ$-Cartier and the pair $(X,\Delta)$ is klt.  The proof is provided in \autoref{lem:base_change_ideals}.
\end{remark}

\subsection{Semi-positivity in characteristic zero}

Here we restate and provide the reference for the characteristic zero version of the main result of \autoref{sec:semi-positivity_engine}.

\begin{theorem}
\label{lem:semipositivity_char_zero}
Assume we are in the following situation:
\begin{enumerate}
\item $f \colon X \to T$ is an equidimensional fibration between normal projective varieties with normal general fibers,
\item $U \subseteq T$ is a non-empty open set, 
\item $\Gamma \geq 0$ be a $\bQ$-divisor such that $K_{X/T} + \Gamma|_{f^{-1}U}$ is $\bQ$-Cartier,
\item  $L$ is a nef  $\bQ$-Cartier $\bQ$-divisor such that $K_{X/T} + \Gamma + L$ is $f$-nef over $U$, and
\item     
$(X_t, \Gamma_t)$ is klt for every closed point $t \in U$
\end{enumerate}
 Then $K_{X/T} + \Gamma+L$ is pseudo-effective.

{\scshape Case $(*)$:}
If $T$ is a curve, and $K_{X/T} +  \Gamma+L$ is $f$-nef $\bQ$-Cartier (so globally, not only over $U$), then $K_{X/T} + \Gamma + L$ is not only pseudo-effective but also nef. 
\end{theorem}

\begin{proof}
Fix an ample divisor $H$ on $X$. As in the proof of \autoref{lem:semipositivity}, \emph{it is enough to show that for every rational number $\varepsilon>0$ and every divisible enough integer $m>0$ (where divisible enough depends on $\varepsilon$),   $f_* \sO_X(m (K_{X/T} + \Gamma + L + \varepsilon H ) ) $ is weakly positive}. So, we prove this statement in the rest of the proof, for which we fix a rational number $\varepsilon>0$. Note that by fixing a general effective $\bQ$-divisor $\Gamma' \sim_{\bQ} L + \varepsilon H$ and setting $\Delta:= \Gamma + \Gamma'$, it is enough to show that  $f_* \sO_X(m (K_{X/T} + \Delta ) )$ is weakly positive. Additionally,  $(X_t, \Delta_t)$ is klt for every $t \in U$, after possibly shrinking $U$.

Let $\pi : Y \to X$ be a log resolution of $(Y, \Delta)$, and set $g:= f \circ \pi$. Let $\Delta_Y$ be the crepant boundary on $Y$. Then, it is enough to show that 
  $g_* \sO_Y(m (K_{Y/T} + \{\Delta_Y\} ) ) $ is weakly positive for every integer $m>0$ divisible enough:
\begin{equation*}
K_Y + \{ \Delta_Y\}  = \pi^* (K_X + \Delta ) - 
\explparshift{420pt}{-150pt}{\lfloor \Delta_Y \rfloor}{ $-\lfloor \Delta_Y \rfloor|_{g^{-1}U}$ is effective, and $\pi$-exceptional $\Rightarrow$ the natural homomorphism  $g_* \sO_Y(m (K_{Y/T} + \{\Delta_Y\}  ) ) \cong f_* \pi_*  \sO_Y(m (K_{Y/T} + \{\Delta_Y\}  ) ) \hookrightarrow f_* \sO_X(m (K_{X/T} + \Delta ) )$ is generically an isomorphism}
\end{equation*}
However, the semi-positivity of   $g_* \sO_Y(m (K_{Y/T} + \{\Delta_Y\} ) ) $ is exactly the statement of \cite[Thm 1.1]{Fujino_Notes_on_the_weak_positivity_theorems}. 

\end{proof}

\subsection{Flatness and reducedness of fibers}

In this subsection we state the generalization of the result of \cite{Lu_Tu_Zhang_Zheng_On_semistability_of_Albanese_maps} to the case of an arbitrary map $(X,\Delta) \to T$ where $(X,\Delta)$ is klt and $-K_{X/T}-\Delta$ is nef.  As described in \autoref{sec:flatness} our argument is very similar to the one given in \emph{loc.cit}.

\begin{proposition}
\label{prop:flatness_and_reducedness_char_zero}
Let $f \colon (X,\Delta) \to T$ be a morphism from a normal, projective klt pair to a smooth variety such that $-K_{X/T}-\Delta$ is nef.  Then $\Delta$ does not contain any irreducible components of the fibers of $f$, the morphism $f$ is flat and all its fibers are reduced.
\end{proposition}

\begin{proof}

We can prove that $f$ is equidimensional, and hence flat (underlying varieties of klt pairs are Cohen--Macaulay), arguing as in \autoref{thm:nef_anti_rel_canonical_flat} replacing the semi-positivity theorem used at the end of the proof with the analogous characteristic zero result provided in \autoref{lem:semipositivity_char_zero}.  Similarly, we can obtain the non-reducedness of the fibers and the statement concerning the containment of irreducible components of fibers of $f$ by using the line of reasoning from \autoref{prop:reduced_fibres} with the substitution of the semi-positivity result.
\end{proof}

\subsection{Numerical flatness in characteristic zero}

\begin{proposition}
\label{thm:numerical_flatness_char_zero}
Let $f \colon (X,\Delta) \to T$ be a fibration from a normal, projective pair to a smooth curve.  Assume that $-K_{X/T}-\Delta$ is ($\QQ$-Cartier) nef, and that the general fiber $(X_t,\Delta_t)$ is klt.  Let $L$ be an $f$-ample line bundle on $X$ satisfying $L^{\dim X} = 0$.  Then the following statements hold true. 
\begin{enumerate}
\item\label{numerical_flatness_item1} For every $m>0$ the vector bundle $f_*\cO_X(mL)$ is numerically flat.
\item\label{numerical_flatness_item2} Let $D$ be an irreducible component of the divisor $\Supp \Delta$.  Assume that $D$ is $\QQ$-Cartier.  Then for every $m>0$ such that $mL - D$ is $f$-ample, the sheaf $f_*\cO_X(mL - D)$ is a numerically flat vector bundle.
\end{enumerate}
\end{proposition}

\begin{proof}

The following argument is a characteristic zero combination of \autoref{thm:nef} and \autoref{thm:K_trivial_numerically_flat}.  We recall the argument for the sake of completeness and because some characteristic $p>0$ subtleties are not relevant in the present proof.

\emph{We first prove statement \autoref{numerical_flatness_item1}.} We observe that by \autoref{prop:flatness_and_reducedness_char_zero} the map $f$ is flat with reduced fibers, and hence we can use the semi-positivity result from \autoref{lem:semipositivity_char_zero}.  We now proceed to the proof of the fact that $f_*\cO_X(mL)$ is nef.  We begin by showing that $L$ is nef.  We fix $\varepsilon>0$. Consider $L + \varepsilon f^* H$. It is enough to prove, by limiting $\varepsilon \to 0$, that $L+ \varepsilon f^* H$ is nef. Set $N:= q(L + \varepsilon f^* H)$, for $q \gg 0$. By \autoref{lem:asymptotic_RR_relative_ample}, we know that \begin{equation*}
h^0(X, N)  \geq  q^{d+1} \frac{(L + \varepsilon f^* H)^{d+1}}{(d+1)!} + O(q^d).
\end{equation*}
Since
\begin{equation*}
(L + \varepsilon f^* H)^{d+1} 
\expl{=}{$L^{n+1} =0$}
\varepsilon (d+1) (\deg H)  \left(L_t^n \right)    >0,
\end{equation*}
we obtain that $h^0(X,N) \neq 0$ (using $q \gg 0$). So, choose  $\Gamma  \in |N|$. For $\varepsilon'$ small enough,  $\left(X_t, \Delta_t+ \varepsilon' \Gamma_t \right)$ is klt for $t \in T$ general. Hence, according to \autoref{lem:semipositivity_char_zero}, the following divisor is nef  
\begin{equation*}
\underbrace{\underbrace{K_{X/T} +\Delta + \varepsilon' \Gamma}_{\parbox{92pt}{\scriptsize $\left(X_t, \Delta_t+ \varepsilon' \Gamma_t \right)$ is klt for $t \in T$ general} } + (\underbrace{ - (K_{X/T} + \Delta)}_{\textrm{nef}})}_{\textrm{$f$-ample}}
 \equiv 
 \varepsilon' q(L + \varepsilon f^* H) 
\end{equation*}
Hence $L + \varepsilon f^* H$  is nef, which concludes our statement.  

Now we deduce that $f_*\cO_X(mL)$ is nef.  For this purpose, we observe that
\[
mL = K_{X/T} + \Delta + \underbrace{mL + (-K_{X/T}-\Delta)}_{\textrm{nef and $f$-ample}}
\]
and hence by \cite[Proposition 6.3]{Codogni_Patakfalvi_Positivity_of_the_CM_line_bundle_for_families_of_K-stable_klt_Fanos} the pushforward $f_*\cO_X(mL)$ is nef.

We now show that $f_*\cO_X(mL)^*$ is nef.  Assume for the sake of contradiction that this is not the case. 
According to the characteristic zero version of \autoref{lem:nef_L_min_max}, this means that 
\begin{equation*}
\mu_{\max} \left( f_* \sO_X(aL) \right) > 0 
\end{equation*}
We take the maximal destabilizing sheaf $\cE$ of $f_* \sO_X(aL)$.  According to \cite[Corollary 6.4.14]{Lazarsfeld_Positivity_in_algebraic_geometry_II} , $\sE^{\otimes r}$ is semi-stable with slope $r\mu (\sE)$. In particular, we may fix $r >0$ such that $\mu\left(\sE^{\otimes r}\right) > 2 g(T) +1$. Consider then the composition
\begin{equation*}
\sE^{\otimes r} \hookrightarrow \left( f_* \sO_X(aL) \right)^{\otimes r} \to f_*\sO_X(raL) .
\end{equation*}
This is non-zero (diagonal tensors do not go to zero). So, we obtain that 
\begin{equation*}
\mu_{\max}\left(f_*\sO_X(raL)\right) > 2 g(T) +1.
\end{equation*}
Hence for any closed point $t \in T$, by \cite[Prop 5.7]{Codogni_Patakfalvi_Positivity_of_the_CM_line_bundle_for_families_of_K-stable_klt_Fanos}, we have 
$h^0\left(T, ( f_*\sO_X(raL))(-t) \right)\neq 0$. 

Choose $0 \neq \Gamma \in \left| raL - X_t\right| $.  Then for any $0 < \varepsilon \ll 1$, $\left(X_t, \Delta_t + \varepsilon \Gamma_t\right)$ is  a klt pair for $t \in T$ general, where the notion of generality depends on $\varepsilon$.  In particular, according to \cite[Corollary 6.4]{Codogni_Patakfalvi_Positivity_of_the_CM_line_bundle_for_families_of_K-stable_klt_Fanos}, the following divisor is nef:
\begin{equation*}
\varepsilon \Gamma \sim_{\bQ} 
\underbrace{
\underbrace{K_{X/T} + \Delta +   \varepsilon \Gamma}_{\parbox{118pt}{\scriptsize $\left(X_t,\Delta_t + \varepsilon \Gamma_t\right)$  is klt for $t \in T$ general  }} 
+(\underbrace{-(K_{X/T} + \Delta)}_{\textrm{nef}}) 
}_{\textrm{$f$-ample}}
\end{equation*}
 However,  then
\begin{equation*}
\Gamma^{d+1} =\left(raL - X_t\right)^{d+1} 
\expl{=}{$L^{d+1}=0$} 
- (n+1) ra L_t^d <0,
\end{equation*}
which contradicts nefness of $\Gamma$.  Hence, our initial assumption was false, which concludes the proof of \autoref{numerical_flatness_item1}.  

\emph{We now proceed to the proof of \autoref{numerical_flatness_item2}}.  

First, we prove that for $m \gg 0$ sufficiently large the divisor $mL - D$ is nef.  Let $\varepsilon = 1/m$ and choose $m \gg 0$ such that $L - \varepsilon D$ is $f$-ample.  The claim follows from the semipositivity result and the relation
\[
L - \varepsilon D = 
\underbrace{
\underbrace{K_{Y/T} + \Delta - \varepsilon D}_{\parbox{95pt}{\scriptsize $\left(X_t,\Delta_t - \varepsilon D_t\right)$ is klt for $t \in T$ general, and $m \gg 0$}} 
+ \underbrace{(-K_{Y/T}-\Delta) + L.}_{\text{nef}}
}_{\text{$f$-ample}}
\]
In order to see that $f_*\cO_X(mL - D)$ is nef, we consider the relation
\[
mL - D = K_{X/T} + \Delta + \underbrace{(-K_{X/T}-\Delta) + mL - D}_{\textrm{nef and $f$-ample}}
\]
and apply \cite[Proposition 6.3]{Codogni_Patakfalvi_Positivity_of_the_CM_line_bundle_for_families_of_K-stable_klt_Fanos} as above.  However, this already implies that $f_*\cO_X(mL - D)$ is numerically flat.  Indeed, the vector bundle $f_*\cO_X(mL - D)$ is nef, and hence of non-negative degree.  It is also a torsion-free subsheaf of $f_*\cO_X(mL)$ which is numerically flat, and hence semistable of degree zero.  This implies that $f_*\cO_X(mL - D)$ is nef of degree zero and hence numerically flat by the result of Langer (see \cite[Proposition 5.1]{Langer_On_the_S-fundamental_group_scheme_I}) recalled in \autoref{prop:num_flat_equivalent_defs}.
\end{proof}

\subsection{The proof of the main result}

\begin{theorem}
\label{thm:char_zero_isotrivial}
Let $f \colon (X,\Delta) \to T$ be a fibration of a normal projective $\bQ$-factorial pair onto a smooth projective variety.  Assume that the $-K_{X/T} - \Delta$ is nef and that $(X,\Delta)$ is klt.  Then $f$ is isotrivial, and in particular all the fibers $f$ are isomorphic. 
\end{theorem}

\begin{proof}
By \autoref{prop:flatness_and_reducedness_char_zero}, $f$ is flat with reduced fibers,  and $\Delta$ does not contain any fiber.  Let $d$ be the relative dimension of $f$.  As in the proof of \autoref{thm:decomposition_theorem} we find an $f$-ample divisor $\wt L$ such that $\wt L^{d+1} \cdot H^{n-d-1} = 0$, for a given very ample divisor $H$ on $T$, and the following conditions are satisfied 
\begin{enumerate}[(a)]
\item the bundles $f_*\cO_X(n \wt L)$ are locally free with formation commuting with any base change,
\item the multiplication map $\Sym^n \cO_X(n \wt L) \to \cO_X(n\wt L)$ is surjective, for every $n > 0$,
\item the line bundle $L$ is sufficiently $f$-ample (so that $nL-D$ is $f$-ample for every $n>0$) and satisfies the vanishing $R^if_*\cO_X(nL - D) = 0$ for every irreducible component $D$ of $\Supp \Delta$ and $i>0$.
\end{enumerate}
The intersection condition implies that for a general smooth curve $C$ in the intersection $H^{n-d-1}$ the line bundle $L = \wt L_{|X_C}$, where $X_C = X \times_T C$, satisfies the equality $L^{d+1} = 0$.  Here general means that it is not contained in the proper closed subvariety of the base over which the fibers are not normal. 

We claim that the morphism $f_C \colon X_C \to C$ satisfies the assumptions of \autoref{thm:numerical_flatness_char_zero}.  Indeed, we see that the fibers of $f$ are reduced by \autoref{prop:flatness_and_reducedness_char_zero}.  Since $X$ is klt it is also Cohen--Macaulay and therefore $X_C$ is $S_2$.  The reducedness of all the fibers and normality of the general fibre imply that $X_C$ is also $R_1$ and consequently it is normal by Serre's criterion.  Moreover the general fiber of $(X_C,\Delta_C)\to C$ is klt, and by \cite[Proposition 2.1]{Codogni_Patakfalvi_Positivity_of_the_CM_line_bundle_for_families_of_K-stable_klt_Fanos} we have the base change formula $-K_{X_C/C}-\Delta_C = (-K_{X/T}-\Delta)|_{X_C}$ yielding nefness of $-K_{X/C} - \Delta_C$.  As a result, we see using the base change and \autoref{thm:numerical_flatness_char_zero} that the line bundle $\wt L$ and the family of smooth curves in $H^{n-d-1}$ satisfies the assumption of \autoref{lemma:isotriviality_criterion_cao_algebraizedII} which implies that the morphism $f \colon X \to T$ is isotrivial over $T$.  In order to extend this to the boundary setting, we observe that the conclusion of \autoref{remark:numerical_flatness_divisor} and \autoref{thm:numerical_flatness_char_zero} could be used in the boundary setting, that is, the required sufficient $f$-ampleness property holds uniformly for a family of curves covering a big open subset of $T$.   More precisely, the condition (c) above and the $(X,\Delta)$ klt assumption ensure that both the base change property for the relative section rings (see \autoref{lem:base_change_ideals}) and the $f$-ampleness of $nL-D$ necessary in \autoref{thm:numerical_flatness_char_zero} are satisfied.  We may extend isotriviality from the open subset covered by a family of smooth complete intersection to the whole $T$ by taking a closure, because the divisors are uniquely determined in codimension one.  
\end{proof}

\begin{theorem}
\label{thm:char_zero_isotrivial_generically}
Let $f \colon (X,\Delta) \to T$ be a fibration of a normal $\bQ$-factorial projective pair onto a normal projective variety.  Assume that the $-K_{X/T} - \Delta$ is nef and that $(X,\Delta)$ is klt.  Then $f$ is generically locally isotrivial.
\end{theorem}
\begin{proof}
The proof is very similar to the argument given above.  We claim that the fibrations is isotrivial over the smooth locus $U$ of $T$.  Indeed, as in the proof of \autoref{thm:char_zero_isotrivial}, we construct a divisor $\wt L$ and a family of complete intersection curves covering $T$.  Since $T$ is normal, its singular locus is of codimension at least two, and therefore there exists a subfamily of smooth curves $\{C_s\}_{s \in S}$ covering $U$.  We conclude as above.
\end{proof}

\begin{corollary}
\label{cor:dps_projective_singular}
Let $(X,\Delta)$ be a projective $\bQ$-factorial klt pair such that $-K_X-\Delta$ is nef.  Then the Albanese morphism $\pi \colon (X,\Delta) \to \Alb X$ is an isotrivial fibration.
\end{corollary}

\begin{proof}
First, using the result of \autoref{sec:relative_canonical_bundle_normal_spaces}, we see that $-K_{X/T} - \Delta$ is nef.  Since $(X, \Delta)$ is klt, the general fiber $(X_t,\Delta_t)$ of $\pi$ is klt as well, and we may therefore apply \autoref{thm:char_zero_isotrivial} to conclude.
\end{proof}

\bibliographystyle{skalpha} 
\bibliography{includeNice}
\end{document}